\documentclass[11pt,a4paper]{article}

\usepackage{auxhook}
\usepackage{latexsym}
\usepackage{blindtext}

\PassOptionsToPackage{force}{filehook}

\usepackage[T1]{fontenc}
\usepackage[utf8]{inputenc}
\usepackage[main=english,french]{babel}
\usepackage{graphicx}
\usepackage{graphics,subcaption}
\usepackage{fancyhdr}
\usepackage{lmodern}
\usepackage{color}
\usepackage{transparent}
\usepackage{float}
\usepackage{pythontex}
\usepackage{tikz-cd}
\usepackage{listings}
\usepackage[linktocpage=true,citecolor=green, urlcolor=black]{hyperref}
\usepackage{tabularx}
\usepackage[title,titletoc]{appendix}
\usepackage{xcolor}
\usepackage[mathscr]{eucal}
\usepackage{calrsfs}
\usepackage[skins]{tcolorbox}

\usetikzlibrary{shapes.misc}
\usepackage{pgf,tikz}
\usepackage{caption}
\usetikzlibrary{decorations.markings}
\usetikzlibrary{arrows}
\usetikzlibrary{patterns}
\usepackage{tikz-3dplot}
\usetikzlibrary{calc}
\usepackage{amsfonts}
\usepackage{amsmath,amssymb,amsthm}
\usepackage{mathtools}
\usepackage{relsize}
\usepackage{titlesec}
\usepackage{tcolorbox}
\usepackage{upgreek}

\newcolumntype{P}[1]{>{\centering\arraybackslash}p{#1}}

\usepackage{stmaryrd}
\usepackage{trimclip}

\usepackage[shortlabels]{enumitem}

\usepackage{wasysym}

\usepackage{hyperref}
\usepackage[nameinlink]{cleveref}
\crefname{section}{Appendix}{Appendices}

\makeatletter
\DeclareRobustCommand{\shortto}{%
	\mathrel{\mathpalette\short@to\relax}%
}

\newcommand{\short@to}[2]{%
	\mkern2mu
	\clipbox{{.5\width} 0 0 0}{$\m@th#1\vphantom{+}{\shortrightarrow}$}%
}
\makeatother

\theoremstyle{definition}
\newtheorem{definition}{Definition}[section]

\newtheorem{remark}[definition]{Remark}

\newtheorem{lemma}[definition]{Lemma}

\newtheorem{Proposition}[definition]{Proposition}
\newtheorem*{proposition*}{Proposition}
\newtheorem*{conjecture*}{Conjecture}

\titleformat{\subsubsection}
{\normalfont\normalsize\bfseries}{\thesubsubsection}{1em}{}
\titlespacing*{\subsubsection}
{0pt}{3.25ex plus 1ex minus .2ex}{1.5ex plus .2ex}

\usepackage[margin=2.5cm]{geometry}

\usepackage{mathptmx}

\usepackage{authblk}

\usepackage{comment}
\usepackage[strict]{csquotes}
\usepackage[
style=numeric-comp,
autolang=other,
backend=biber,
sorting=nyt,
citereset=none,
maxnames=10,
natbib=true,
hyperref=true,
url=true,
giveninits=true,
]{biblatex}

\title{Moser Regularization of the CR3BP for Symmetric Orbits: Bifurcation Networks of Lunar Orbits}
\author[1]{Cengiz Aydin}
\author[2]{Chankyu Joung}
\affil[1]{Institut für Mathematik, Universität Heidelberg, Germany}
\affil[ ]{\textit {E-mail address: \href{mailto:cengiz.aydin@hotmail.de}{cengiz.aydin@hotmail.de}}}
\affil[2]{Department of Mathematical Science, Seoul National University, South Korea}
\affil[ ]{\textit {E-mail address: \href{mailto:chankyu.joung@gmail.com}{chankyu.joung@gmail.com}}}
\date{}

\begin{document}

\setcounter{page}{1}
\pagenumbering{arabic}	
\maketitle

\begin{abstract}
	We develop a computational framework based on Moser regularization for the analysis of symmetric periodic orbits in the spatial Circular Restricted Three-Body Problem (CR3BP).\ The formulation enables the numerical continuation of periodic orbit families through collision with either primary, together with the computation of Floquet multipliers and stability indices entirely in regularized coordinates.\	We further provide a characterization and numerical detection scheme for symmetric periodic collision orbits.\ We apply this framework to the Earth--Moon CR3BP and uncover bifurcation networks connecting several classical lunar periodic orbit families.\ Higher-period vertical bifurcations of the planar $L_1$ Lyapunov family generate spatial bridges reconnecting different regions of the same family across its collision with the Earth, while corresponding bifurcations of the $L_2$ Lyapunov family form bridges to planar families associated with retrograde comet orbits.\ Furthermore, period-multiplying bifurcations of Near-Rectilinear Halo Orbit (NRHO) families produce a systematic network connecting to planar prograde and retrograde families.\ These bifurcation connections in the near-Moon dynamics retain the characteristic organization by covering multiplicities of the corresponding bifurcation network in the limiting Hill problem.\  Conley--Zehnder indices and related symplectic invariants complement the numerical continuation and provide a topological organization of the results.\
\end{abstract}

\begin{center}
\begin{tabular}{ll}
	\textbf{Keywords} & $\quad$ Moser regularization $\cdot$ restricted three-body problem $\cdot$ symmetric orbits $\cdot$\\
	& $\quad$ Earth--Moon CR3BP $\cdot$ lunar orbits $\cdot$ bifurcation network \\
	\textbf{MSC 2020} & $\quad$ 70F07 $\cdot$ 70G45 $\cdot$ 70H12
\end{tabular}
\end{center}

\setcounter{tocdepth}{5}
\tableofcontents

\section{Introduction}

The Circular Restricted Three-Body Problem (CR3BP) is one of the fundamental problems in celestial mechanics which concerns the motion of an infinitesimal mass moving under the gravitational influence of the Newtonian attraction of two bodies (the primaries) that revolve around each other in circular orbits with uniform velocity.\ The infinitesimal mass may approach a collision with the primaries, causing the equations of motion to become singular.\ Understanding the behavior of trajectories near collisions has motivated the development of a variety of regularization theories \cite{birkhoff,szebehely,celletti,roa}, whose purpose is to transform the singular differential equations into regular ones.\ There are plenty of results on many topics related to the applications of regularization techniques, such as
\begin{itemize}[noitemsep]
	\item[(a)] computation of periodic orbits \cite{broucke,howell_breakwell,belbruno_2019_family},
	\item[(b)] sharper analysis of bridges and bifurcation relationships among orbit families \cite{schmidt,aydin_batkhin,joung_2026_bifurcations}, 
	\item[(c)] study of quasi-periodic motions \cite{palacian} and close encounters \cite{lega,cardin},
	\item[(d)] design of transfer trajectories~\cite{oshima} and trajectory optimization \cite{korneev}, and
	\item[(e)] topology of regularized energy level sets \cite{albers,aydin_contact,frauenfelder_koert,joung_koert}.
\end{itemize}

While regularization in the topics (a)--(d) is essential for the continuation of solutions through collision and the better control of the numerical error near collisions, in the topic (e) regularization aims to complete the flow such that energy level sets become compact, to which topological and symplectic methods can be applied.\

In many celestial mechanical systems, such as in the CR3BP, the most common type of collision singularity is a collision of two bodies, whose mutual gravitational interaction dominates the dynamics near collision.\ Consequently, the local near-collision dynamics is governed by a perturbed Kepler problem.\ This makes classical regularization methods developed for the Kepler problem natural tools to apply within the topics (a)--(e), which include the works of Levi-Civita (LC) \cite{levi_civita}, Kustaanheimo--Stiefel (KS)~\cite{kustaanheimo_stiefel}, Moser \cite{moser_1970_regularization} and Ligon--Schaaf \cite{ligon_schaaf}.\ Historically, according to \cite[Section 3.12]{szebehely}, Euler appears to have been the first to address the problem of regularization when studying the collision of two bodies moving on the same straight line.\ 

Moreover, among the most representative regularization methods within the computational results in (a)--(d) are the LC and KS regularizations, which locally regularize one singularity in planar and spatial problems, respectively.\ We notice that Birkhoff \cite{birkhoff} also proposed a global transformation for the CR3BP to obtain a simultaneous regularization with both primaries for the planar system, which was generalized by Waldvogel \cite{waldvogel} for the spatial system.\ These regularization transformations (except Ligon--Schaaf, see below) rely on a common procedure characterized by three features:\
\begin{itemize}[noitemsep]
	\item[(1)] Rescaling of time (to control the increase of speed at collision),
	\item[(2)] performing a suitable change of coordinates, and
	\item[(3)] extending or blowing up the phase space by using the preservation of energy.
\end{itemize} 

The beginning part (1) may be attributed to Sundman \cite{sundman}, while the main ideas of (2) and (3) were contributed by LC for the planar Kepler problem, by using a conformal squaring in the physical plane, which suggests using complex notation.\ This was generalized by KS for the spatial system by using an elegant quaternion formalism.\ Another, more geometric way to regularize the Kepler problem was proposed by Moser \cite{moser_1970_regularization}, which is valid in every dimension.\ By exchanging position and momentum variables, rescaling of time, and using the stereographic projection, Moser realized each fixed negative energy level as the unit cotangent bundle of the $n$-sphere, transforming symplectically Keplerian motion into the geodesic flow of $S^n$.\ Indeed, the collision singularity, which corresponds to unbounded momentum, is mapped to a regular point on the $n$-sphere.\ Moser's construction is level-wise, meaning that the symplectomorphism is defined separately on each fixed negative energy level set.\ We notice that Ligon--Schaaf provides a global symplectomorphism, without reparametrizing time, that is valid for the whole negative-energy region at once, such that the Kepler flow becomes a reparametrized geodesic flow.\

However, the Moser regularization has been used more often in theoretical symplectic research in the topic (e) so far, but is rarely applied in the computational topics (a)--(d).\ To our knowledge, only the works \cite{belbruno_2019_family,joung_2026_bifurcations} utilize the Moser regularization in numerical computations.\ Therefore, the purpose of this paper is to provide a Moser regularization scheme for the numerical computation of periodic orbits of the CR3BP, with a focus on symmetric periodic orbits.\ Even though the Moser and KS regularizations are related~\cite{kummer}, a key benefit of the Moser regularization is that the transformation is symplectic.\ This allows us to perform all relevant calculations directly in Moser regularized coordinates, including Floquet multipliers, stability indices, and symplectic invariants of the linearized flow, rather than evaluating through a chain of transformations as is done in KS coordinates~\cite{howell_breakwell,aydin_batkhin}.\

In this paper, we develop a numerical framework for symmetric periodic orbit families of the spatial CR3BP in Moser regularized coordinates.\ Building on the constrained Hamiltonian formulation for numerical integration and stability analysis developed in~\cite{belbruno_2019_family,joung_2026_bifurcations}, we extend the methods in the following directions.\ First, we adapt the symmetry reduction of periodic orbits~\cite{bray_goudas} to the ambient space of the Moser regularized energy surface.\ This yields four classes of symmetric periodic orbits with explicit boundary conditions.\ We derive symmetry reduced formulas for the monodromy matrix and stability indices, and formulate single and multiple shooting correction schemes over half or quarter period segments, entirely within regularized coordinates.\ Second, we characterize the geometry of periodic collision orbits.\ On each regularized energy surface, the collision states associated with a primary form a $2$-sphere, whose intersections with the fixed point sets of the two reversing symmetries determine the collision configurations.\ We establish a criterion for detecting collisions along periodic orbit families and develop a Newton correction scheme.\ Finally, we provide a publicly available implementation of the numerical algorithms, together with test scripts that reproduce all numerical tables, including stability indices, Krein signatures, and Conley--Zehnder indices, available at \url{https://github.com/ckjoung/moser}.\

We apply this framework to lunar orbit families in the Earth--Moon CR3BP, focusing on the planar Lyapunov families about $L_1$ and $L_2$, as well as the halo families that bifurcate vertically from them.\ Lunar orbits have played a central role in space mission design for decades~\cite{howell_barden,koon,gomez_et_al}, and studies have revealed a rich network of bifurcations and connections among these and related periodic orbit families~\cite{doedel,aydin_batkhin}.\ However, the complete bifurcation network, including connections through collision, has not yet been studied systematically, and several important features of this network remain unknown.\ Upon further continuation, the $L_1$ and $L_2$ Lyapunov families approach collision with the primaries and undergo higher-period vertical bifurcations, producing spatial families with periods integer multiples of that of the corresponding Lyapunov orbits.\ Such bifurcations have been studied in Hill's modification of the CR3BP~\cite{kalantonis} and its radiation pressure extension~\cite{ragos_kalantonis}, but have been less explored in the Earth--Moon system.\ Moreover, collision is commonly treated as an endpoint of numerical continuation~\cite{doedel,kalantonis}, whereas regularized continuation reveals families to typically pass through collision and connect to other branches.\ This issue is particularly relevant to Near-Rectilinear Halo Orbits (NRHOs).\ In the Hill problem, the NRHO family relates to the vertical collision orbit, which together with planar retrograde and direct families form a structured bifurcation network organized by covering multiplicities~\cite{aydin_batkhin}.\ Here, we use the Moser regularized framework to explore the corresponding network in the Earth--Moon CR3BP.\

\textit{Structure of the paper.}\ In Section \ref{sec:2} we briefly recall the CR3BP dynamics.\ Section \ref{sec:3} is devoted to the Moser regularization scheme for the CR3BP.\ We first derive the regularized equations of motion.\ Then, we discuss the analysis of symmetric periodic orbits, computation of Floquet multipliers and stability indices, continuation, Newton correction, and multiple shooting method within Moser regularized coordinates.\ In Section \ref{sec:4} we describe a numerical scheme for detecting symmetric periodic collision orbits in Moser regularization of the CR3BP.\ In Section \ref{sec:5}, by using the techniques developed in Section \ref{sec:3} and Section \ref{sec:4}, we present the numerical results in the Earth--Moon CR3BP.\ For a systematic organization of the bifurcation results, we use symplectic methods and homological invariants (mostly related to the Conley--Zehnder index) that are summarized in Appendix~\ref{appendix_a}.\ Finally, Appendix~\ref{appendix_b} contains tables of data.\

\section{Circular Restricted Three-Body Problem (CR3BP)}
\label{sec:2}

\begin{figure}[t]
	\centering
	\includegraphics[scale=1]{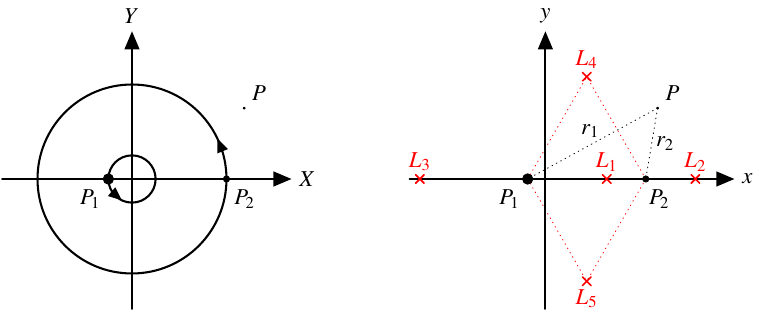}
	\caption{CR3BP (from \cite{aydin_dro}) in an inertial reference frame $(X,Y)$ and in a rotating reference frame $(x,y)$, with primaries $P_1$ and $P_2$, infinitesimal body $P$ and libration points $L_i$ ($i=1,2,3,4,5$).}
	\label{fig:cr3bp}
\end{figure}

The Circular Restricted Three-Body Problem (CR3BP) describes the orbital motion of a particle $P$ with infinitesimal mass under the gravitational influences of two point masses, a heavy primary $P_1$ and a light primary $P_2$ of masses $m_1$ and $m_2$ $(m_1 \geq m_2)$, respectively.\ This model assumes that the two point masses move along circles in the same plane with constant angular velocity around their common center of mass (see Figure \ref{fig:cr3bp}).\ It is convenient to choose the units of mass, distance and time such that the gravitational constant is~1, and the period of the circular orbits is $2\pi$.\ In these units, both the distance between the two primaries and their total mass are normalized to one.\ The mass parameter is defined as $\mu = m_2 / (m_1 + m_2)$, so that, under the normalization $m_1 + m_2 = 1$, the masses of $P_1$ and $P_2$ are $1 - \mu$ and~$\mu$, respectively, with $0 < \mu \leq 1/2$.\ It is also convenient to study the dynamics within a coordinate system that rotates with an angular velocity of the orbital angular rate of the primaries.\ Then, the positions of the primaries $P_1$ and $P_2$ relative to the rotating frame are fixed correspondingly at $(-\mu,0,0)$ and at $(1-\mu,0,0)$ (see Figure \ref{fig:cr3bp}).\ Within these assumptions, the dynamics of the CR3BP is expressed via the Hamiltonian function \cite[Chapter 10]{szebehely}:\
\begin{align} \label{hamiltonian}
	H(x,y,z,p_x,p_y,p_z) = \frac{1}{2} \left( p_x^2 + p_y^2 + p_z^2 \right) - \frac{1-\mu}{r_1} - \frac{\mu}{r_2} + p_x y - p_y x,
\end{align}
where $r_1 = \left( (x+\mu)^2 + y^2 + z^2 \right)^{\frac{1}{2}}$ and $r_2 = \left( (x - 1 + \mu)^2 + y^2 + z^2 \right)^{\frac{1}{2}}$ indicate the distances from the infinitesimal body to the primaries.\ The equations of motion are
\begin{align} \label{equations}
	\left\{ \dot{k} = \frac{\partial H}{\partial p_k},\quad \dot{p}_k = - \frac{\partial H}{\partial k} \right\},\quad \text{with } k = x,y,z.
\end{align}

The CR3BP framework has five points of equilibrium, also called libration points or Lagrange points \cite[Chapter 4]{szebehely}:\ Three collinear points $L_1$ (located within the primaries), $L_2$ and $L_3$ (outside the interval joining the primaries), and two equilateral triangular points $L_4$ and $L_5$, as illustrated in Figure \ref{fig:cr3bp}.\ Furthermore, since $H$ is a first integral of the system, periodic solutions come in a smooth one-parameter family, parameterized by the energy.\ The phase space of the system is 6-dimensional and endowed with the standard symplectic form $\omega_0 = \sum dp_k \wedge dk$ ($k=x,y,z$).\ The flow of the Hamiltonian vector field $X_H$, that is uniquely determined by $dH(\cdot) = \omega_0( \cdot, X_H )$, is equivalent to the equations~\eqref{equations}.\

\section{Moser Regularization of the CR3BP}
\label{sec:3}

\subsection{Moser Regularized Equations of Motion}

We describe the regularization of collision with the light primary located at $\mathbf{m} = (1-\mu,0,0)$.\ The regularization of collisions with the heavy primary at $\mathbf{e} = (-\mu, 0, 0)$ follows analogously by replacing $\mu$ by $1 - \mu$ and reflecting the dynamics by the symmetry $(x, y, z, p_x, p_y, p_z) \mapsto -(x, y, z, p_x, p_y, p_z)$.\

\begin{figure}[t]
	\centering
	\includegraphics[width=0.5\textwidth]{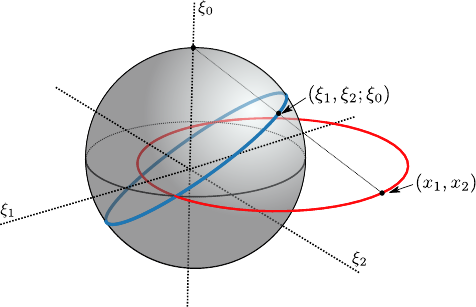}
	\caption{
		Illustration of the Moser regularization (from~\cite{joung_2026_bifurcations}), that is shown for the planar case where the momentum variables $(x_1,x_2)$ are mapped to $(\xi_0,\xi_1,\xi_2)$ on the unit 2-sphere by the inverse stereographic projection.\ The collision orbits pass through the North Pole $(\xi_0,\xi_1,\xi_2)=(1,0,0)$.}
	\label{fig:moser-reg}
\end{figure}

First we shift the origin to $\mathbf{m}$ and exchange symplectically the position and momentum coordinates:
\begin{equation}\label{change_pos_mom}
	\begin{aligned}
		\mathbf{x} &= (x_1, x_2, x_3) = -(p_x,\, p_y - (1 - \mu),\, p_z), \\
		\mathbf{y} &= (y_1, y_2, y_3) = (x - (1 - \mu),\, y,\, z).
	\end{aligned}
\end{equation}
Then, we apply the Moser transformation~\cite{moser_1970_regularization}, which maps symplectically the Euclidean phase space
\[
(\mathbf{x},\mathbf{y}) \in \mathbb{R}^3 \times \mathbb{R}^3
\]
to the cotangent bundle of the unit 3-sphere,
\[
T^*S^3
=
\left\{
(\boldsymbol{\xi},\boldsymbol{\eta}) \in \mathbb{R}^4 \times \mathbb{R}^4
\,\middle|\,
|\boldsymbol{\xi}|^2 = 1,\,
\boldsymbol{\xi}\cdot\boldsymbol{\eta} = 0
\right\} \subset T^* \mathbb{R}^4,
\]
that we identify with the tangent bundle of the unit 3-sphere $T S^3$ by using the standard metric on $S^3$.\ Here, $\boldsymbol{\xi}=(\xi_0, \xi_1, \xi_2, \xi_3)$ denotes the coordinates on the sphere $S^3\subset \mathbb{R}^4$, and $\boldsymbol{\eta}=(\eta_0, \eta_1, \eta_2, \eta_3) \in T_{\xi}S^3$ is a tangent vector.\ The key idea of the Moser regularization is to transform the momentum coordinates~$\mathbf{x}$ to regularized momentum $\boldsymbol{\xi}$ by utilizing the inverse stereographic projection (see Figure~\ref{fig:moser-reg}), so that the collision locus $|\mathbf{x}| = \infty$ corresponds to an ordinary point, namely the North Pole $\boldsymbol{\xi}=(1,0,0,0)$ of the sphere.\ The transformation, that is induced by the stereographic projection, is given by
\begin{gather}\label{eqn:moser-trans}
	x_1 = \frac{\xi_1}{1-\xi_0}, \quad x_2 = \frac{\xi_2}{1-\xi_0}, \quad x_3 = \frac{\xi_3}{1-\xi_0}, \\
	y_1 = \eta_0 \xi_1 + (1-\xi_0)\eta_1, \quad y_2 = \eta_0 \xi_2 + (1-\xi_0)\eta_2, \quad y_3 = \eta_0 \xi_3 + (1-\xi_0)\eta_3,\nonumber
\end{gather}
with the inverse transformation that is presented by
\begin{equation}
	\begin{aligned}
		\xi_0 &= \frac{|\mathbf{x}|^2 - 1}{|\mathbf{x}|^2 + 1}, \qquad \xi_k = \frac{2 x_k}{|\mathbf{x}|^2 + 1}, \qquad k=1,2,3, \\
		\eta_0 &= \mathbf{x}\cdot\mathbf{y}, \qquad \eta_k = \frac{|\mathbf{x}|^2 + 1}{2}\, y_k
		- (\mathbf{x}\cdot\mathbf{y}) x_k,
		\qquad k=1,2,3. \\
	\end{aligned}
	\label{eqn:moser-transform}
\end{equation}
By construction, these transformations are symplectic, as they preserve the canonical 1-form,
$ \sum_{j = 1}^3 y_j dx_j = \sum_{j = 0}^3 \eta_j d \xi_j. $
Furthermore, it holds that
\begin{align}\label{relation}
	|\boldsymbol{\eta}| = \frac{1}{2} \left( |\mathbf{x}|^2 + 1 \right) |\mathbf{y}| = \frac{\mathbf{|y|}}{1 - \xi_0}.
\end{align}

\noindent
Now let $c$ be the energy level of the original Hamiltonian.\ After the transformation (\ref{change_pos_mom}), we use a new time parameter $s$ and define a new Hamiltonian by
\begin{align}\label{transformed_ham}
	ds = \frac{dt}{|\mathbf{y}|},\quad K_c(\mathbf{x},\mathbf{y}) = |\mathbf{y}| ( H - c).
\end{align}
By utilizing the transformation (\ref{eqn:moser-trans}) and the identities (\ref{relation}) into (\ref{transformed_ham}), we obtain a transformed Hamiltonian on $T^* S^3$, which we denote by the same letter,
$$K_c (\boldsymbol{\xi},\boldsymbol{\eta}) = f(\boldsymbol{\xi},\boldsymbol{\eta}) |\boldsymbol{\eta}| - \mu,$$
where
\begin{equation}
	\label{eqn:energy_f}
	\begin{split}
		f(\boldsymbol{\xi},\boldsymbol{\eta}) &= 1 - (1 - \xi_0) \left( c + \frac{1}{2} \right) + (1 - \xi_0)(\eta_1 \xi_2 - \eta_2 \xi_1)\\
		&\quad - (1 - \mu) (1 - \xi_0) \left( \frac{1}{| \eta_0 \vec \xi + (1 - \xi_0) \vec \eta - (-1,0,0) |} + \eta_0 \xi_1 + (1 - \xi_0)\eta_1 + \frac{1}{2}(1 - \mu) \right),
	\end{split}
\end{equation}
with $\vec{\xi} = (\xi_1,\xi_2,\xi_3)$, and $\vec \eta = (\eta_1,\eta_2,\eta_3)$.\ We further consider the Hamiltonian $Q_c$ on $T^* S^3$ that is defined by
\begin{equation}
	\label{eqn:energy_Q}
	Q_c (\boldsymbol{\xi},\boldsymbol{\eta}) = \frac{1}{2 \mu} \left(  K_c (\boldsymbol{\xi},\boldsymbol{\eta}) + \mu \right)^2 = \frac{1}{2 \mu} f(\boldsymbol{\xi},\boldsymbol{\eta})^2 |\boldsymbol{\eta}|^2. 
\end{equation}
The level set $Q_c = \mu/2$ consists of two components $f(\boldsymbol{\xi},\boldsymbol{\eta})|\eta| = \mu$ and $f(\boldsymbol{\xi},\boldsymbol{\eta})|\eta| = -\mu$.\ The first corresponds to $K_c = 0$, on which the dynamics of $Q_c$ are a reparametrization of the dynamics of $K_c$, which in turn corresponds to the dynamics of $H$ on $H = c$.\ To keep the flow on the constraint surface $T^* S^3$, we define the two constraints
\begin{equation} 
	\label{eqn:constraints_g}
	g_1(\boldsymbol{\xi},\boldsymbol{\eta}) = |\boldsymbol{\xi}|^2 - 1,\qquad g_2(\boldsymbol{\xi},\boldsymbol{\eta}) =  \boldsymbol{\xi} \cdot \boldsymbol{\eta}, 
\end{equation}
and determine the constrained Hamiltonian equations with respect to the Lagrange multiplier approach, that is characterized by
\begin{align*}
	\begin{dcases}
		\dot{\xi}_k = \frac{\partial Q_c}{\partial \eta_k} + \lambda_1 \frac{\partial g_1}{\partial \eta_k} + \lambda_2 \frac{\partial g_2}{\partial \eta_k}  = \frac{\partial Q_c}{\partial \eta_k} + \lambda_2 \xi_k,\qquad k=0,1,2,3\\
		\dot{\eta}_k = - \frac{\partial Q_c}{\partial \xi_k} - \lambda_1 \frac{\partial g_1}{\partial \xi_k} - \lambda_2 \frac{\partial g_2}{\partial \xi_k} = - \frac{\partial Q_c}{\partial \xi_k} - 2 \lambda_1 \xi_k - \lambda_2 \eta_k,\qquad k=0,1,2,3,
	\end{dcases}
\end{align*}
with Lagrange multipliers $\lambda_1 \equiv \lambda_1(\boldsymbol{\xi},\boldsymbol{\eta})$, and $\lambda_2 \equiv \lambda_2(\boldsymbol{\xi},\boldsymbol{\eta})$.\ The preservation of the two constraints implies that
$$ \sum_{j = 0}^{3} \xi_j \dot{\xi}_j = 0, \qquad \sum_{j = 0}^{3} \dot{\xi}_j \eta_j + \sum_{j = 0}^{3} \xi_j \dot{\eta}_j = 0, $$
from which the Lagrange multipliers are determined to be
$$ 
\lambda_1 =  \frac{1}{2|\boldsymbol{\xi}|^2} \left( \sum_{j = 0}^{3} \eta_j \frac{\partial Q_c}{\partial \eta_j} - \sum_{j = 0}^{3} \xi_j \frac{\partial Q_c}{\partial \xi_j} \right),
\qquad 
\lambda_2 = - \frac{1}{|\boldsymbol{\xi}|^2} \sum_{j = 0}^{3} \xi_j \frac{\partial Q_c}{\partial \eta_j} .
$$
Finally, the constrained Hamiltonian equations that are associated to $Q_c | _{T^* S^3}$ are given by
\begin{equation}
	\begin{aligned}
		\begin{dcases}
			\dot{\xi}_k = \frac{\partial Q_c}{\partial \eta_k} - \frac{\xi_k}{|\boldsymbol{\xi}|^2} \left( \sum_{j = 0}^{3} \xi_j \frac{\partial Q_c}{\partial \eta_j} \right),\qquad k=0,1,2,3,\\
			\dot{\eta}_k = - \frac{\partial Q_c}{\partial \xi_k} - \frac{\xi_k}{|\boldsymbol{\xi}|^2} \left( \sum_{j = 0}^{3} \eta_j \frac{\partial Q_c}{\partial \eta_j} - \sum_{j = 0}^{3} \xi_j \frac{\partial Q_c}{\partial \xi_j} \right) + \frac{\eta_k}{|\boldsymbol{\xi}|^2} \left( \sum_{j = 0}^{3} \xi_j \frac{\partial Q_c}{\partial \eta_j} \right),\qquad k=0,1,2,3,
		\end{dcases}
	\end{aligned}
	\label{eqn:ode-moser}
\end{equation}
which coincide with the equations provided in \cite{belbruno_2019_family,joung_2026_bifurcations}.\ The normalization by $|\boldsymbol{\xi}|^2$ leaves the flow on $T^*S^3$ unchanged, but is required to make $g_1$ and $g_2$ first integrals of~\eqref{eqn:ode-moser}.\ Consequently, the linearized flow preserves the variations of $g_1$ and $g_2$, which will be used in Section~\ref{sec:floquet} to determine the Floquet multipliers.\ The following summarizes the correspondence between unregularized and regularized trajectories.
\begin{Proposition}
	A trajectory $\boldsymbol{\gamma}(t)$ of the unregularized CR3BP on the energy level $H=c$ corresponds to a trajectory $\tilde{\gamma}(s)$ of the Moser regularized system~\eqref{eqn:ode-moser} on the component $K_c = 0$ ($f(\boldsymbol{\xi},\boldsymbol{\eta})|\eta| = \mu$) of the energy level $Q_c=\mu/2$ via the transformations \eqref{eqn:moser-trans} and \eqref{eqn:moser-transform}, and the time reparametrization $dt / ds = |\mathbf{y}| = (1 - \xi_0) |\boldsymbol{\eta}|$.
\end{Proposition}

\subsection{Analysis of Periodic Orbits with respect to Symmetries}

Symmetries are symplectic or anti-symplectic involutions of the phase space that leave the equations invariant.\ While the symplectic symmetries denote the time-preserving involutions, the anti-symplectic symmetries denote the time-reversal ones in the Hamiltonian context~\cite{lamb_roberts}.\ From the inspection of the Hamiltonian (\ref{hamiltonian}), two anti-symplectic and one symplectic symmetries within the unregularized CR3BP are typically found, that are correspondingly given by:\
\begin{align*}
	\rho_x (x,y,z,p_x,p_y,p_z) = &(x,-y,-z,-p_x,p_y,p_z) \text{ ($\pi$-rotation around the $x$-axis)},\\
	\rho_{xz} (x,y,z,p_x,p_y,p_z) = &(x,-y,z,-p_x,p_y,-p_z) \text{ (reflection at the $xz$-plane)},\\
	\sigma (x,y,z,p_x,p_y,p_z) = &(x,y,-z,p_x,p_y,-p_z) \text{ \phantom{$-$}(reflection at the plane $\{z=0\}$)}.
\end{align*}
These satisfy $\rho_x \circ \rho_{xz} = \rho_{xz} \circ \rho_{x} = \sigma$.\ The restriction of the spatial system to the fixed point set of $\sigma$,
$$ \text{Fix}(\sigma) = \{ (x,y,0,p_x,p_y,0) \}, $$
describes the planar CR3BP, in which an inherent anti-symplectic symmetry is characterized by
\begin{equation*}
	\rho : (x,y,0,p_x,p_y,0) \mapsto (x,-y,0,-p_x,p_y,0),
\end{equation*}
corresponding to the reflection at the $x$-axis.\ These classical symmetries of the spatial CR3BP take the following form in the Moser regularized coordinates:\
\begin{equation}\label{symmetries_moser}
	\begin{aligned}
		\tilde{\rho}_x(\boldsymbol{\xi},\boldsymbol{\eta})
		&=
		(\xi_0, -\xi_1, \xi_2, \xi_3,
		-\eta_0, \eta_1, -\eta_2, -\eta_3), \\
		\tilde{\rho}_{xz}(\boldsymbol{\xi},\boldsymbol{\eta})
		&=
		(\xi_0, -\xi_1, \xi_2, -\xi_3,
		-\eta_0, \eta_1, -\eta_2, \eta_3), \\
		\tilde{\sigma}(\boldsymbol{\xi},\boldsymbol{\eta})
		&=
		(\xi_0, \xi_1, \xi_2, -\xi_3,
		\eta_0, \eta_1, \eta_2, -\eta_3).
	\end{aligned}
\end{equation}
These leave the Hamiltonian $Q_c$ invariant, and since the Moser transformations are symplectic, they have the same properties as the corresponding underlying symmetry in the unregularized system:\ Both $\tilde{\rho}_x$ and~$\tilde{\rho}_{xz}$ are anti-symplectic, while $\tilde{\sigma}$ is symplectic.\ They also satisfy $\tilde{\rho}_x \circ \tilde{\rho}_{xz} = \tilde{\rho}_{xz} \circ \tilde{\rho}_x = \tilde{\sigma}$, and the restriction of the regularized system to the fixed point set of $\tilde{\sigma}$,
\[
\mathrm{Fix}(\tilde{\sigma})
=
\{(\xi_0, \xi_1, \xi_2, 0, \eta_0, \eta_1, \eta_2, 0)\},
\]
corresponds to the Moser regularization of the planar CR3BP.\ Restricted to this subset, the remaining anti-symplectic symmetry reduces to
\begin{equation*}
	\tilde{\rho} : 
	(\xi_0,\xi_1,\xi_2,0,\eta_0,\eta_1,\eta_2,0)
	\mapsto
	(\xi_0,-\xi_1,\xi_2,0,-\eta_0,\eta_1,-\eta_2,0).
\end{equation*}
While the constraint function $g_1$ is invariant under each of the symmetries \eqref{symmetries_moser}, the constraint function $g_2$ is invariant under $\tilde{\sigma}$ and anti-invariant under both $\tilde{\rho}_x$ and $\tilde{\rho}_{xz}$.\ Inspecting the signs in~\eqref{eqn:ode-moser}, however, both $\tilde{\rho}_x$ and $\tilde{\rho}_{xz}$ are still anti-symplectic symmetries, and $\tilde{\sigma}$ is a symplectic symmetry of the equations~\eqref{eqn:ode-moser}.\ Therefore, classical techniques with respect to such symmetries can be applied:\ If $\tilde{\gamma}(s) =   (\boldsymbol{\xi}(s),\boldsymbol{\eta}(s))$ is a solution of \eqref{eqn:ode-moser}, then
$$ \tilde{\rho}_x \left( \tilde{\gamma}(s) \right),\qquad \tilde{\rho}_{xz} \left( \tilde{\gamma}(s) \right),\qquad \tilde{\sigma} \left( \tilde{\gamma}(s) \right), $$
are also solutions of \eqref{eqn:ode-moser}.\ Furthermore, the regularized flow $\Phi^{s}$ on the constraint surface $T^* S^3$ satisfies
\begin{align}\label{identity_1}
	\tilde{\sigma} \circ \Phi^{s} = \Phi^{s} \circ \tilde{\sigma},
\end{align}
and for the anti-symplectic symmetries $\tilde{\rho}_i$, $i \in \{x,xz\}$, it holds that
\begin{align}\label{identity_2}
	\tilde{\rho}_i \circ \Phi^{-s} = \Phi^{s} \circ \tilde{\rho}_i.
\end{align}
In the special case when a periodic orbit $\tilde{\gamma}(s)$ is invariant under one or two of the anti-symplectic symmetries, then, with respect to the fixed point sets of $\tilde{\rho}_x$ and $\tilde{\rho}_{xz}$,
\begin{equation} \label{eqn:fixedptset}
	\text{Fix} (\tilde{\rho}_x) = \{ (\xi_0,0,\xi_2,\xi_3,0,\eta_1,0,0) \},\qquad \text{Fix} (\tilde{\rho}_{xz}) = \{ (\xi_0,0,\xi_2,0,0,\eta_1,0,\eta_3) \}, 
\end{equation}
these constraints lead to the following four classes of symmetric periodic orbits with the period $S$:
\begin{itemize}
	\item[(I)] A periodic orbit $\tilde{\gamma}(s)$ is called \textit{simple symmetric with respect to $\tilde{\rho}_x$} if
	\begin{align}\label{bound_cond_1}
		\xi_1(0) = \eta_0(0) = \eta_2(0) = \eta_3(0) = 0,\quad \text{ and } \quad \xi_1(S/2) = \eta_0(S/2) = \eta_2(S/2) = \eta_3(S/2) = 0.
	\end{align}
	\item[(II)] A periodic orbit $\tilde{\gamma}(s)$ is called \textit{simple symmetric with respect to $\tilde{\rho}_{xz}$} if
	\begin{align}\label{bound_cond_2}
		\xi_1(0) = \xi_3(0) = \eta_0(0) = \eta_2(0) = 0,\quad \text{ and } \quad \xi_1(S/2) = \xi_3(S/2) = \eta_0(S/2) = \eta_2(S/2) = 0.
	\end{align}
	\item[(III)] A periodic orbit $\tilde{\gamma}(s)$ is called \textit{doubly symmetric with respect to $\tilde{\rho}_x$-$\tilde{\rho}_{xz}$} if
	\begin{align}\label{bound_cond_3}
		\xi_1(0) = \eta_0(0) = \eta_2(0) = \eta_3(0) = 0,\quad \text{ and } \quad \xi_1(S/4) = \xi_3(S/4) = \eta_0(S/4) = \eta_2(S/4) = 0.
	\end{align}
	\item[(IV)] A periodic orbit $\tilde{\gamma}(s)$ is called \textit{doubly symmetric with respect to $\tilde{\rho}_{xz}$-$\tilde{\rho}_x$} if
	\begin{align}\label{bound_cond_4}
		\xi_1(0) = \xi_3(0) = \eta_0(0) = \eta_2(0) = 0,\quad \text{ and } \quad \xi_1(S/4) = \eta_0(S/4) = \eta_2(S/4) = \eta_3(S/4) = 0.
	\end{align}
\end{itemize}
\noindent
Geometrically, the conditions \eqref{bound_cond_1} and \eqref{bound_cond_2} characterize orbits that meet the same symmetry section orthogonally at half period intervals, while \eqref{bound_cond_3} and \eqref{bound_cond_4} characterize orbits that alternate between two distinct symmetry sections at quarter period intervals.\ Furthermore, classical techniques \cite{bray_goudas} that simplify the analysis of symmetric solutions by reducing the integration of the equations to the half period (for simple symmetric solutions), or to the quarter period (for doubly symmetric solutions) can be applied:\ Let $V(s) := d \Phi^{s}$ denote the regularized linearized flow that solves the variational equations with respect to~\eqref{eqn:ode-moser} along a $S$-periodic solution.\ Then, by using the identities \eqref{identity_1} and \eqref{identity_2}, the monodromy matrix $M := V(S) \in \mathbb{R}^{8 \times 8}$ can be obtained by exploiting the symmetries:
\begin{equation}
	\begin{aligned}
		\text{(i) simple symmetric w.r.t. $\tilde{\rho}_x$:} & \quad M = \Sigma_x V^{-1}(S/2) \Sigma_x V(S/2),\\ 
		\text{(ii) simple symmetric w.r.t. $\tilde{\rho}_{xz}$:} & \quad M = \Sigma_{xz} V^{-1}(S/2) \Sigma_{xz} V(S/2), \\
		\text{(iii) doubly symmetric w.r.t. $\tilde{\rho}_x$ and $\tilde{\rho}_{xz}$:} & \quad M = \left[ \Sigma_x V^{-1}(S/4) \Sigma_{xz} V(S/4) \right]^2, \\
		\text{(iv) doubly symmetric w.r.t. $\tilde{\rho}_{xz}$ and $\tilde{\rho}_x$:} & \quad M = \left[ \Sigma_{xz} V^{-1}(S/4) \Sigma_x V(S/4) \right]^2,
	\end{aligned}
	\label{monodromy_symmetry_moser}
\end{equation}
where
\[
\Sigma_x = \mathrm{diag}\{1,-1,1,1,-1,1,-1,-1\}, \qquad
\Sigma_{xz} = \mathrm{diag}\{1,-1,1,-1,-1,1,-1,1\}.
\]

\subsection{Floquet Multipliers and Stability Indices}
\label{sec:floquet}

The restriction of the monodromy matrix $M \in \mathbb{R}^{8\times 8}$ to the tangent space of the constraint manifold $T^*S^3$ is symplectic.\ Therefore, the eigenvalues of $M$ occur in reciprocal pairs and are of the form 
\[
\{1, 1, 1, 1, \lambda_1, \lambda_1^{-1}, \lambda_2, \lambda_2^{-1}\},
\]
where the four eigenvalues equal to $1$ appear due to the conserved quantities $Q_c$, $g_1$, and~$g_2$, as well as the time invariance.\ We observe that, when switching from the original unregularized system to the Moser regularized framework, the two constraint functions $g_1$ and $g_2$ produce two additional eigenvalues equal to $1$.\ The remaining non-trivial eigenvalues, which are the \textit{Floquet multipliers}, match with those computed in the original unregularized system.\ Therefore, for each pair of Floquet multipliers the classical notion of the \textit{stability index} can be used, that is defined by
$$ s_k := \frac{1}{2} \left( \lambda_k + \lambda_k^{-1} \right),\qquad k=1,2. $$
The classifications of the stability are identical to the ones in the original unregularized case:\
\begin{itemize}[noitemsep]
	\item \textit{Elliptic:}\ $-1 < s_k = \cos \theta < 1$.\ It is equivalent to $\lambda_k = e^{i \theta}$.
	\item \textit{Positive hyperbolic:}\ $s_k > 1$.\ It is equivalent to $\lambda_k \in \mathbb{R}_{>0} \setminus \{1\}$.
	\item \textit{Negative hyperbolic:}\ $s_k < -1$.\ It is equivalent to $\lambda_k \in \mathbb{R}_{<0} \setminus \{-1\}$.
	\item \textit{Complex instability:}\ Case of quadruples $\{ \lambda, \lambda^{-1}, \overline{\lambda}, \overline{\lambda}^{-1} \}$.
\end{itemize}
The stable periodic solutions are those that simultaneously satisfy $|s_1| < 1$ and $|s_2| < 1$.\ Furthermore, in view of the additional two trivial eigenvalues that are equal to 1, the modified version of the classical formulas for the computations of the stability indices \cite{bray_goudas,broucke_stab} is now given by
\begin{align*}
	s_1 = - \frac{1}{4} \left( \alpha + \sqrt{\beta} \right),\quad s_2 = - \frac{1}{4} \left( \alpha - \sqrt{\beta} \right),
\end{align*}
where
$$ \alpha = 4 - \textnormal{tr}(M),\quad \beta  = 2 \textnormal{tr}(M^2) - \alpha^2. $$

\begin{remark}
In the numerical implementation, explicitly forming $V^{-1}$ in~\eqref{monodromy_symmetry_moser} may be inaccurate, especially for highly unstable orbits.\ In such cases, the following equivalent generalized eigenvalue formulation can be used to avoid computing the inverse directly.\ For example, if
$
	M = \Sigma V^{-1} \Sigma V,
$
then its eigenvalues $\lambda$ (including the trivial ones) can be computed from the generalized eigenvalue problem
\[
	Az = \lambda Bz, \qquad
	A = \Sigma V, \quad B = V\Sigma,
\]
which can be numerically preferable to forming $V^{-1}$ explicitly.\ For doubly symmetric orbits, where
$
	M = \left[ \Sigma_i V^{-1} \Sigma_j V \right]^2,
$
we first compute the generalized eigenvalues $\nu$ of
\[
	Az = \nu Bz, \qquad
	A = \Sigma_j V, \quad B = V \Sigma_i,
\]
and obtain the eigenvalues of $M$ by $\lambda=\nu^2$.\
\end{remark}

\subsection{Vertical Self-Resonant Bifurcations}

Owing to the symplectic symmetry $\sigma$, and its extension $\tilde{\sigma}$ to the Moser regularized system, the linear stability of a planar periodic orbit decomposes into an ``in-plane'' and a ``vertical'' (or ``out-of-plane'') component~\cite{henon_1973}.\ We denote the corresponding stability indices by $s_p$ and $s_v$, respectively.\ In the context of vertical bifurcations from planar periodic orbits, $s_v=1$ corresponds to an out-of-plane bifurcation giving rise to spatial periodic orbits of the same period, while $s_v=-1$ corresponds to a period-doubling bifurcation.\ Higher-order vertical bifurcations occur at \textit{Vertical Self-Resonant} (VSR) orbits, for which the vertical Floquet multipliers are roots of unity.\ Equivalently, the vertical stability index satisfies the self-resonance condition
\begin{align}
	\label{sv}
	s_v = \cos\left(2\pi\frac{d}{k}\right),
\end{align}
where $d,k\in\mathbb{N}$ ($k\geq 3$), $d/k$ is irreducible, and $d < k$, with the numerator chosen to be either $d$ or $k - d$ so that $2\pi d/k$ matches the rotation angle determined by the Krein sign (see Appendix~\ref{appendix_a}).\ A bifurcation point satisfying \eqref{sv} is referred to as a $d$:$k$ \textit{Vertical Self-Resonant (VSR) orbit}, according to its self-resonance ratio $d$:$k$, with $k$ denoting the \textit{multiplicity} of the resonance.\ Let $T$ denote the period of a $d$:$k$ VSR orbit.\ According to \cite{robin_markellos_1980,robin_markellos_1983}, its $k$-fold covering, of period $kT$, gives rise to two distinct families of spatial periodic orbits.\ Their symmetry properties depend on the parity of the multiplicity $k$:
\begin{itemize}[noitemsep]
	\item[a)] \textit{Odd multiplicity $k$.}\ The two bifurcating families consist of simple symmetric orbits:\ one family is $\rho_x$-symmetric, while the other is $\rho_{xz}$-symmetric.\ By the $\sigma$-symmetry, each of these branches has a symmetric counterpart.\
	\item[b)] \textit{Even multiplicity $k$.}\ Both bifurcating families consist of doubly symmetric orbits with respect to $\rho_x$ and $\rho_{xz}$.\
\end{itemize}
The same vertical bifurcation mechanism applies in the Moser regularized system, with the symmetries $\rho_x$, $\rho_{xz}$, and $\sigma$ replaced by their counterparts in the Moser regularized system $\tilde{\rho}_x$, $\tilde{\rho}_{xz}$, and $\tilde{\sigma}$, respectively.\

\subsection{Numerical Continuation and Newton Correction}
We compute families of symmetric periodic orbits using pseudo-arclength continuation~\cite[Section~10.2]{kuznetsov} in the unregularized system, allowing continuation through birth--death (or fold) bifurcations.\ As an orbit approaches collision, we switch to Moser regularization and use natural continuation in the Jacobi energy parameter $C=-2c$.\ The correction scheme in the regularized system, described below, depends on the symmetry type of the orbit.\ Instead of enforcing full periodicity, we impose boundary conditions on a fraction of the orbit depending on the type of symmetry at initial and final states with respect to the fixed point sets~\eqref{eqn:fixedptset} of $\tilde{\rho}_x$ and $\tilde{\rho}_{xz}$.\ Moreover, rather than using a Poincar\'e return map, we treat the period as an additional unknown, since intersections with prescribed sections are generally continuous only locally along a family.\ 

\medskip
\noindent\emph{Case 1: Initial point on Fix($\tilde{\rho}_{xz}$).}\ To formulate a Newton correction scheme we take a reduced set of free variables
\begin{equation}
\label{eqn:free_vector_case1}
X = (\xi_2, \eta_3, s),
\end{equation}
where $s$ denotes the integration time corresponding to the half period (for simple symmetric orbits), or to the quarter period (for doubly symmetric orbits).\ Given these free variables, we lift the coordinates to the appropriate fixed point set and energy surface.\ In the $\tilde{\rho}_{xz}$-symmetric case, $\xi_0$ is determined by
\[
\xi_0 = \pm \sqrt{1 - \xi_2^2},
\]
with the sign chosen to be consistent with the previously computed solution. The remaining coordinate~$\eta_1$ is then obtained by numerically solving the energy condition,
\[
K_c(\boldsymbol{\xi}, \boldsymbol{\eta}) = f(\boldsymbol{\xi}, \boldsymbol{\eta})|\boldsymbol{\eta}| - \mu = 0,
\]
for the root near the $\eta_1$ value of the previously computed solution.\ We then integrate the lifted initial condition for time $s$ under the regularized flow $\Phi^s$~\eqref{eqn:ode-moser} and measure the deviation from the target fixed point set.\ This defines the defect map
\[
F(X) = Y\bigl(\Phi^s(\boldsymbol{\xi}, \boldsymbol{\eta})\bigr),
\]
where $Y$ extracts the components that must vanish on the fixed point set:
\begin{equation}
\label{eqn:final_coords}
Y =
\begin{cases}
	(\xi_1,\eta_2,\eta_3), & \text{for $\tilde{\rho}_x$-symmetry}, \\
	(\xi_1,\xi_3,\eta_2), & \text{for $\tilde{\rho}_{xz}$-symmetry}.
\end{cases}
\end{equation}
The condition $F(X)=0$ defines a square $3\times 3$ correction scheme.\ The correction is performed using a Newton iteration
\[
X_{i+1} = X_i - DF(X_i)^{-1} F(X_i),
\qquad i=0,1,2,\dots,
\]
until convergence, i.e.,\ $\|F(X_i)\| < \varepsilon$ for a prescribed tolerance $\varepsilon>0$. The converged solution yields the corrected initial condition and the corresponding period.\ Note that only three components are enforced in~$Y$, since the remaining condition follows from the constraint $g_2(\boldsymbol{\xi},\boldsymbol{\eta})=\boldsymbol{\xi}\cdot\boldsymbol{\eta}=0$, which implies $\eta_0=0$ provided $\xi_0 \neq 0$ at the endpoint; this is verified a posteriori after convergence.

The Jacobian $DF$ of $F$ is obtained by combining the variational flow with the differential of the lifting.\ Writing $Y = (y_1,y_2,y_3)$, the chain rule yields, for $i=1,2,3$,
\begin{align}
	\frac{\partial y_i}{\partial \xi_2}
	&=
	v_{y_i, \xi_0} \frac{\partial \xi_0}{\partial \xi_2}
	+
	v_{y_i, \xi_2}
	+
	v_{y_i, \eta_1} \frac{\partial \eta_1}{\partial \xi_2}, \label{eqn:jacobian_case1_1}\\
	\frac{\partial y_i}{\partial \eta_3}
	&=
	v_{y_i, \eta_3}
	+
	v_{y_i, \eta_1} \frac{\partial \eta_1}{\partial \eta_3}, \label{eqn:jacobian_case1_2}\\
	\frac{\partial y_i}{\partial s}
	&=
	\dot{y_i}\left(\Phi^s(\boldsymbol{\xi},\boldsymbol{\eta})\right). \nonumber
\end{align}
Here $v_{y_i,x}$ denotes the entry of the variational matrix $V(s)\in\mathbb{R}^{8\times 8}$ corresponding to the $x$-component at time $0$ and the $y_i$-component at time $s$, while $\dot{y}_i$ is given by~\eqref{eqn:ode-moser}.\ The lifting derivatives follow from explicit differentiation and the implicit function theorem:
\[
\frac{\partial \xi_0}{\partial \xi_2}
=
\mp \frac{\xi_2}{(1-\xi_2^2)^{1/2}},
\quad
\frac{\partial \eta_1}{\partial \xi_2}
=
-\frac{\partial K_c/\partial \xi_0 \cdot \partial \xi_0/\partial \xi_2 + \partial K_c/\partial \xi_2}{\partial K_c/\partial \eta_1},
\quad
\frac{\partial \eta_1}{\partial \eta_3}
=
-\frac{\partial K_c/\partial \eta_3}{\partial K_c/\partial \eta_1}.
\]

\medskip
\noindent\emph{Case 2: Initial point on Fix($\tilde{\rho}_x$).}\ In the $\tilde{\rho}_x$-symmetric case, we take
\begin{equation}
	\label{eqn:free_vector_case2}
X = (\xi_2, \xi_3, s),
\end{equation}
as reduced set of free variables.\ For the remaining variables one uses
\[
\xi_0 = \pm \sqrt{1 - \xi_2^2 - \xi_3^2},
\]
and determines $\eta_1$ from the energy constraint. The Jacobian entries become
\begin{align}
	\frac{\partial y}{\partial \xi_i}
	&=
	v_{y, \xi_0} \frac{\partial \xi_0}{\partial \xi_i}
	+
	v_{y, \xi_i}
	+
	v_{y, \eta_1} \frac{\partial \eta_1}{\partial \xi_i}, 
	\qquad
	i = 2, 3, \label{eqn:jacobian_case2}\\
	\frac{\partial y}{\partial s}
	&=
	\dot{y}\left(\Phi^s(\boldsymbol{\xi},\boldsymbol{\eta})\right),\nonumber
\end{align}
with
\[
\frac{\partial \xi_0}{\partial \xi_i}
=
\mp \frac{\xi_i}{(1-\xi_2^2-\xi_3^2)^{1/2}},
\quad
\frac{\partial \eta_1}{\partial \xi_i}
=
-\frac{\partial K_c/\partial \xi_0 \cdot \partial \xi_0/\partial \xi_i + \partial K_c/\partial \xi_i}{\partial K_c/\partial \eta_1},
\quad
i = 2, 3.
\]

\medskip
\begin{remark}
	The correction schemes assume $\xi_0\neq 0$ at both endpoints of the segment:\ at the initial point to ensure the parametrization of $\xi_0$ to be regular, and at the final point so that $\eta_0=0$ is implied.\ In instances where $\xi_0 = 0$ occurred, we avoided the singularity by choosing the symmetry point with $\xi_0\neq 0$ as the initial point and temporarily correcting the full period orbit instead of the half period segment.\
\end{remark}

\subsection{Multiple Shooting Correction}
\label{sec:3.6}
For efficient continuation of high-period branches containing orbits that approach or pass through collision, especially those arising from period-three and period-four bifurcations, we use a multiple shooting correction in the regularized variables.\ The continuation parameter remains the Jacobi energy parameter $C = -2c$.\ Our implementation follows the standard multiple shooting formulation for periodic orbit correction, for instance, see~\cite{bury}.\

Let $s$ denote the integration time of the symmetry-reduced portion of the orbit.\ Thus, $s$ is either a half period for simple symmetric orbits or a quarter period for doubly symmetric orbits.\ Given an approximate orbit segment, we choose $N$ shooting nodes
$
x_0,x_1,\ldots,x_{N-1}
$
along this segment and set
$
\Delta s = \frac{s}{N}.
$
The first node is constrained to lie on the fixed point set of the prescribed symmetry and on the energy surface.\ Hence, only the free coordinates of $x_0$ are included as unknowns as in~\eqref{eqn:free_vector_case1} and~\eqref{eqn:free_vector_case2}, namely $u_0=(\xi_2, \eta_3)$ for $x_0\in \text{Fix}(\tilde{\rho}_{xz})$ and $u_0=(\xi_2, \xi_3)$ for $x_0\in \text{Fix}(\tilde{\rho}_{x})$.\ The intermediate nodes $x_1, \ldots, x_{N-1}$ are treated as full phase-space unknowns, and the integration time $s$ is also included.\ Thus, the Newton unknown is of the form
\[
X =
\left( u_0, x_1, \ldots, x_{N-1}, s \right) \in \mathbb{R}^{3 + 8 (N - 1)}.
\]

The constraint map consists of two types of equations.\ First, continuity is imposed between consecutive shooting nodes:
\[
F_i(X) = \Phi^{\Delta s}(x_i) - x_{i+1} = 0,
\qquad i = 0, \ldots, N-2.
\]
For the last segment, the corresponding symmetry boundary condition is imposed at its endpoint:
\[
F_{N - 1}(X) = Y\!\left(\Phi^{\Delta s}(x_{N-1})\right)=0,
\]
where $Y$, defined in~\eqref{eqn:final_coords}, extracts the components that vanish on the target fixed point set.\ Together, these equations define a square Newton system.\

The Jacobian is formulated blockwise from the variational flow along each segment.\ For the matching condition $F_i = 0$ $(i=1,\cdots,N-2)$, the nonzero derivatives are
\[
	d_{x_i}F_i=d\Phi^{\Delta s}(x_i),\qquad
	d_{x_{i+1}}F_i=-I,\qquad
	d_sF_i=\frac{1}{N}X_Q\!\left(\Phi^{\Delta s}(x_i)\right),
\]
where $X_Q$ denotes the regularized vector field.\ For the first matching condition $F_0 = 0$, the derivative with respect to $u_0$ is composed with the differential of the lifting map:
\[
	d_{u_0}F_0
	=
	d\Phi^{\Delta s}(x_0)\,d_{u_0}x_0.
\]
The resulting derivatives are given in the same way as in the single shooting case~\eqref{eqn:jacobian_case1_1}--\eqref{eqn:jacobian_case1_2}, and~\eqref{eqn:jacobian_case2}, with $v_{y_i,x}$ now referring to entries of $d\Phi^{\Delta s}(x_0)$.\ The remaining derivatives of $F_0$ (w.r.t.\ $x_1$ and $s$) are obtained in the same way as above.\ Finally, for the endpoint condition $F_{N - 1} = 0$, the corresponding derivatives are projected onto the components selected by $Y$:
\[
	d_{x_{N-1}}F_{N-1}
	=
	\pi_Y \, d\Phi^{\Delta s}(x_{N-1}),\qquad
	d_sF_{N-1}
	=
	\frac{1}{N} \pi_Y \, X_Q\!\left(\Phi^{\Delta s}(x_{N-1})\right).
\]

\section{Symmetric Periodic Collision Orbits of the CR3BP}
\label{sec:4}

We use the term \emph{periodic collision orbit} for an orbit which converges, in physical position space, to one of the primaries both in backward time, $t\to t_1$, and in forward time, $t\to t_2$, where $- \infty < t_1 < t_2 < + \infty$, and which becomes periodic after regularization.\ Such orbits appear along continuations of periodic orbit families, and many examples in the planar CR3BP are provided in the work of Broucke~\cite{broucke}.\ The collisions occur also for the planar $L_1$ and $L_2$ Lyapunov families, and a computer-assisted proof for the collision $L_1$ Lyapunov orbit in the Earth--Moon CR3BP was given in~\cite{capinski_kepley_mireles_james}.\ The properties of the intersection of the stable and unstable manifolds of this collision orbit were used in~\cite{capinski_pasiut} to prove existence of oscillatory collision approaches.\ 

Periodic collision orbits can play an organizing role in networks of periodic orbit families.\ In the Hill problem, the rectilinear vertical collision family is one of the basic families participating in connections between planar and spatial periodic orbit families~\cite{aydin_batkhin}.\ Bifurcations from vertical collisions give rise to spatial halo and butterfly families which persist into CR3BP models~\cite{joung_2026_bifurcations}.\ Thus, locating collision members can reveal branch connections which are obscured by the singularity.\

In Moser regularized coordinates, collision with the light primary
$\mathbf m = (1-\mu, 0, 0)$
corresponds to the North Pole
$\boldsymbol{\xi} = (1, 0, 0, 0)$
of the regularized momentum sphere $\boldsymbol{\xi} \in S^3$.\ This simplifies the expression for the regularized energy, where from~\eqref{eqn:energy_f} and \eqref{eqn:energy_Q} we get
\[
f(\boldsymbol{\xi}, \boldsymbol{\eta}) = 1,
\qquad
Q_c(\boldsymbol{\xi}, \boldsymbol{\eta}) = \frac{|\boldsymbol{\eta}|^2}{2\mu}.
\]
Together with the energy condition
$Q_c = \frac{\mu}{2}$
and the constraint $g_2 = \boldsymbol{\xi}\cdot\boldsymbol{\eta} = 0$ from~\eqref{eqn:constraints_g}, the collision locus over $\mathbf{m}$ on the regularized energy surface is the sphere
\[
\mathcal{C}_{\mathbf m} = \{
	(1, 0, 0, 0, 0, \eta_1, \eta_2, \eta_3) \, | \, \eta_1^2 + \eta_2^2 + \eta_3^2 = \mu^2
\} \simeq S^2.
\]
Its intersection with the fixed point sets of the symmetries determine the possible symmetric collision configurations:
\begin{equation*}
	\begin{aligned}
		&\mathcal{C}_{\mathbf m} \cap \text{Fix}(\tilde{\rho}_{xz})
		= \{
			(1, 0, 0, 0, 0, \eta_1, 0, \eta_3) \, | \, \eta_1^2 + \eta_3^2 = \mu^2
		\} \simeq S^1,
		\\
		&\mathcal{C}_{\mathbf m} \cap \text{Fix}(\tilde{\rho}_{x})
		= \{
			(1, 0, 0, 0, 0, \eta_1, 0, 0) \, | \, \eta_1 = \pm\mu
		\} \simeq S^0.
	\end{aligned}
\end{equation*}
This leads to the following description.
\begin{Proposition}[Symmetric periodic collisions]
	\label{prop:symmetric-collision-configurations}
	Let
	$\tilde{\gamma}(0)=(\boldsymbol{\xi}(0), \boldsymbol{\eta}(0))$, where $\boldsymbol{\xi}(0), \boldsymbol{\eta}(0) \in \mathbb{R}^4$,
	be a collision point of a periodic collision orbit $\tilde{\gamma}(s)$ of the Moser regularized CR3BP~\eqref{eqn:ode-moser}.\ Then,
	\[
	\boldsymbol{\xi}(0) = (1, 0, 0, 0),
	\qquad
	\boldsymbol{\eta}(0) = (0, \eta_1, \eta_2, \eta_3),
	\qquad
	\eta_1^2 + \eta_2^2 + \eta_3^2 = \mu^2.
	\]
	Moreover:
	\begin{enumerate}[(1)]
		\item Let $\tilde{\gamma}$ be planar, i.e.,\ $\eta_3=0$, and $\tilde{\rho}$-symmetric.\ Then,
		$
		\tilde{\gamma}(0)\in \operatorname{Fix}(\tilde{\rho})
		$
		if and only if $\eta_2=0$ and $\eta_1=\pm\mu$.
		
		\item Let $\tilde{\gamma}$ be spatial and $\tilde{\rho}_{xz}$-symmetric.\ Then,
		$
		\tilde{\gamma}(0)\in \operatorname{Fix}(\tilde{\rho}_{xz})
		$
		if and only if $\eta_2 = 0$ and
		$
		\eta_3 = \pm\sqrt{\mu^2 - \eta_1^2}.
		$
		
		\item Let $\tilde{\gamma}$ be spatial and $\tilde{\rho}_{x}$-symmetric.\ Then,
		$
		\tilde{\gamma}(0)\in \operatorname{Fix}(\tilde{\rho}_{x})
		$
		if and only if $\eta_2=\eta_3=0$ and $\eta_1=\pm\mu$.\ In particular, the symmetric periodic collision orbit is planar.\
	\end{enumerate}
\end{Proposition}
Thus, a genuinely spatial collision cannot occur at a $\tilde{\rho}_x$-symmetric point.\ If it occurs at a symmetric point, it can only occur at a $\tilde{\rho}_{xz}$-symmetric point.\ A $\tilde{\rho}_{x}$-symmetric orbit can still be a collision orbit in the doubly symmetric case, where the collision occurs at the other symmetric point, namely the $\tilde{\rho}_{xz}$-symmetric point.\ From the perspective of numerical search, planar symmetric collision orbits can be located by varying only the energy parameter, while for spatial $\tilde{\rho}_{xz}$-symmetric collision orbits we use $\eta_1$ and the energy as free variables, with $\eta_3$ determined by
$\eta_3=\pm\sqrt{\mu^2-\eta_1^2}.$

\begin{figure}[t!]
	\vspace{1em}
	\centering
	\includegraphics[width=1.0\textwidth]{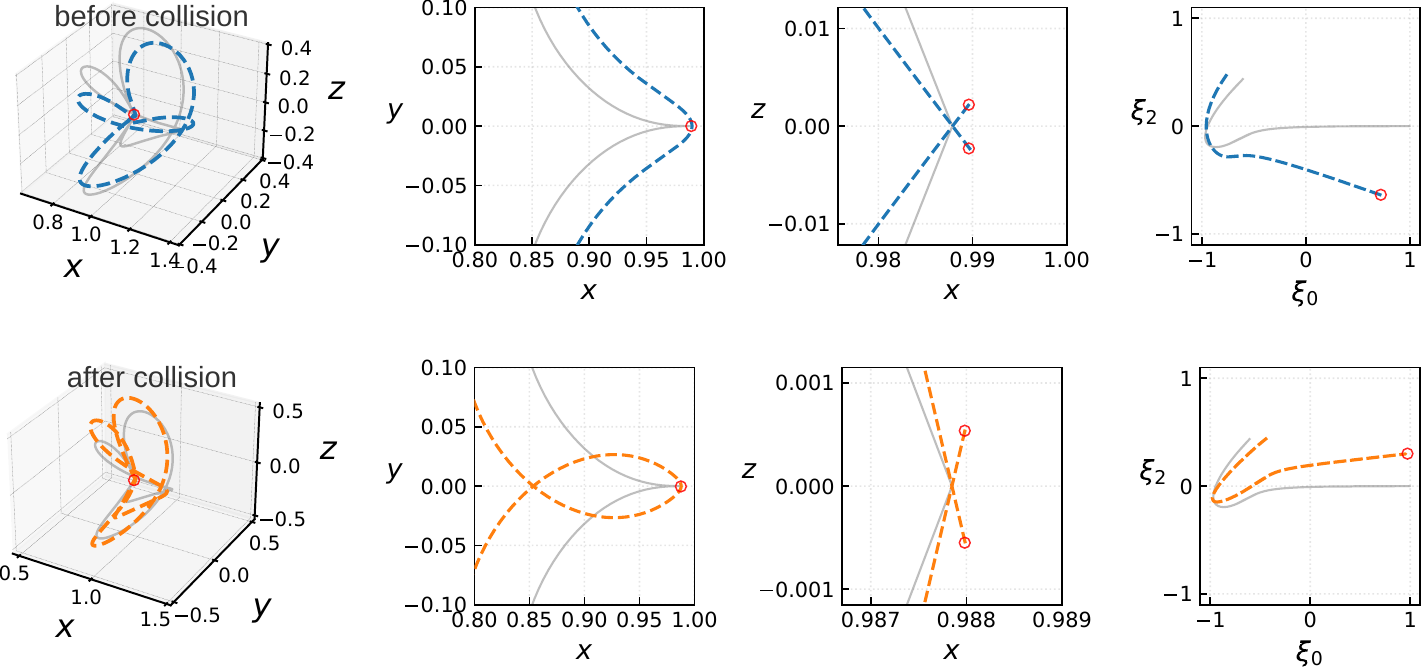}
	\caption{Local behavior of a spatial family continued through collision with the Moon.\ The example shown is the doubly symmetric collision orbit No.~6 in Table~\ref{tbl:periodic-collision-orbits}, generated from a vertical period-doubling bifurcation of the $L_2$ Lyapunov family.\ The gray solid curve denotes the collision orbit, while nearby members of the family are shown in dashed before (top, blue) and after (bottom, orange) collision.\ The red circle indicates the $\rho_{xz}$-symmetric point.\ In physical coordinates, the projections show $x$ and $z$ coordinates bouncing off $1 - \mu$ and $0$, quadratically at collision, respectively, while $p_y$ changes sign.\ In Moser regularized coordinates, the $\tilde{\rho}_{xz}$-symmetric point moves along the circle $\xi_0^2+\xi_2^2=1$ and crosses the collision locus at $(\xi_0,\xi_2)=(1,0)$.}
	\label{fig:collision-evolution}
\end{figure}

We also record the behavior of the orbits in physical coordinates as a family passes through collision.\ Along a $\tilde{\rho}_{xz}$-symmetric spatial family passing through a collision, the symmetric point moves along the circle
\[
\xi_1 = \xi_3 = 0,
\qquad
\xi_0^2 + \xi_2^2 = 1.
\]
Therefore, crossing of the collision locus $\xi_0=1$ is detected by a sign change of $\xi_2$.\ Parametrizing this circle locally near $(\xi_0, \xi_2) = (1, 0)$ by $\epsilon = \xi_2$, so that 
$1 - \xi_0 = 1 - (1 - \xi_2^2)^{1/2} = \frac{1}{2}\epsilon^2 + O(\epsilon^4),$
we have the following behavior in physical coordinates:
\begin{equation*}
	\begin{aligned}
		p_y - (1 - \mu) &= -\frac{\xi_2}{1 - \xi_0} = -\frac{1 + \xi_0}{\xi_2} = -\frac{2}{\epsilon} + O(\epsilon),\\
		x - (1 - \mu) &= (1 - \xi_0) \eta_1 = \frac{\eta_1^*}{2}\epsilon^2 + O(\epsilon^3), \\
		z &= (1 - \xi_0) \eta_3 = \frac{\eta_3^*}{2}\epsilon^2 + O(\epsilon^3),
	\end{aligned}
\end{equation*}
where $\eta_1^*, \eta_3^*$ denote the values at collision.\ Thus, the values of $x - (1 - \mu)$ and $z$ have quadratic contact with the collision point as the family crosses collision, while $p_y$ has a simple pole and changes sign.\ Figure~\ref{fig:collision-evolution} illustrates this behavior.\

We derive a numerical scheme to locate the periodic collision orbits.\ Let
$\mathcal{P}_{c, \pm}$ denote the Poincar\'e map of the Moser regularized flow for a chosen section, restricted to crossings in the direction $\pm$, and $\mathcal{P}^k_{c, \pm}$ its $k$th iterate.\ The energy parameter $c$ is included explicitly, and the number of iterates $k$ is chosen depending on the orbit so that the endpoint of $\mathcal{P}_{c, \pm}^k$ lies on a prescribed symmetry fixed point set.\ For the Poincar\'e section, we either use $\xi_1 = 0$ or $\xi_3 = 0$.\ The section $\xi_1 = 0$ applies to both planar and spatial collision orbits and is used in Section~\ref{sec:5.2}.\ The section $\xi_3 = 0$ can only be applied to out-of-plane orbits satisfying the boundary condition~\eqref{bound_cond_2}, and is used for vertical orbits in Section~\ref{sec:nrho-networks}.\

For a spatial $\tilde{\rho}_{xz}$-symmetric periodic collision orbit at the light primary $\mathbf m=(1-\mu,0,0)$, the initial condition is written as
\[
u_0(\eta_1)
=
(1, 0, 0, 0, 0, \eta_1, 0, \eta_3(\eta_1)),
\qquad
\eta_3(\eta_1) = \pm \sqrt{\mu^2 - \eta_1^2},
\]
where the sign of $\eta_3$ is chosen from the initial guess.\ We then solve
\[
F(\eta_1,c)
=
\pi_Y \circ \mathcal{P}_{c, \pm}^k
\left( u_0(\eta_1) \right)
=
0.
\]
Here $\pi_Y$ projects to the two coordinates that vanish on the target fixed point set.\ In the case where the Poincar\'e section is $\xi_1 = 0$, we use
\[
Y=(\xi_3, \eta_2)
\quad
\text{for } \operatorname{Fix}(\tilde{\rho}_{xz}),
\qquad
Y=(\eta_2, \eta_3)
\quad
\text{for } \operatorname{Fix}(\tilde{\rho}_{x}).
\]
In the case where the section is $\xi_3 = 0$ (only for target $\operatorname{Fix}(\tilde{\rho}_{xz})$), we use $Y = (\xi_1, \eta_2)$.
Thus, $F(\eta_1, c)=0$ imposes the two nontrivial boundary conditions defining the target fixed point set.\ The only remaining condition is $\eta_0 = 0$.\ This follows once we check $\xi_0 \neq 0$, since the constraint $g_2 = \boldsymbol{\xi}\cdot\boldsymbol{\eta} = 0$ from \eqref{eqn:constraints_g} is reduced to $\xi_0 \eta_0 = 0$ from the imposed boundary conditions.\

\begin{lemma}
	\label{lem:collision_newton}
	Assume that $u_0(\eta_1)$ is defined as above and that $\mathcal{P}^k_{c,\pm}$ is well-defined.\ Suppose that
	$
	F(\eta_1,c)=0,
	$
	and that the endpoint
	$
	u_k=\mathcal{P}_{c,\pm}^k(u_0(\eta_1))
	$
	satisfies $\xi_0(u_k)\neq 0$.\
	\begin{enumerate}[(1)]
		\item If the target fixed point set is $\operatorname{Fix}(\tilde{\rho}_{xz})$, then the segment is a half-chord of a $\tilde{\rho}_{xz}$-symmetric orbit.\ 
		\item If the target fixed point set is $\operatorname{Fix}(\tilde{\rho}_{x})$, then the segment is a quarter-chord of a doubly symmetric orbit.\
	\end{enumerate}
	In both cases, the orbit is a periodic collision orbit of the Moser regularized flow.\
\end{lemma}
If the non-degeneracy condition 
$
\det D_{(\eta_1, c)} F \neq 0
$
is met, then the periodic collision orbit is locally isolated, which can be identified by applying Newton correction to $F(\eta_1, c) = 0$.\ 
\begin{remark}
	In the planar case, Lemma~\ref{lem:collision_newton} reduces to a 1-dimensional zero finding problem 
	\[
	F(c)
	=
	\pi_{\eta_2} \circ \mathcal{P}_{c, \pm}^k
	\left( u_0 \right)
	=
	0
	\]
	where $u_0$ is fixed at $(\xi_0, \xi_1, \xi_2, \eta_0, \eta_1, \eta_2) = (1, 0, 0, 0, \pm \mu, 0)$.\
\end{remark}
\begin{remark}
	The correction scheme in Lemma~\ref{lem:collision_newton} is a single shooting formulation.\ One may instead use a multiple shooting version by adding intermediate nodes and imposing the corresponding continuity conditions, as in Section~\ref{sec:3.6}.\ We use the single shooting formulation here due to simplicity and because it is sufficient for the computations in Section \ref{sec:5}.\
\end{remark}

\section{Numerical Results in the Earth--Moon CR3BP}
\label{sec:5}

In this section, we apply the Moser regularization procedure described above to investigate bifurcation networks of periodic orbit families in the Earth--Moon CR3BP.\ We first study higher-period vertical bifurcations of the planar $L_1$ and $L_2$ Lyapunov families and continue the resulting spatial families, including through collisions with the primaries.\ We then extend the analysis to period-multiplying bifurcations of the $L_1$ and $L_2$ Near-Rectilinear Halo Orbit (NRHO) families and their connections to planar prograde and retrograde orbit families.\ The mass ratio assumed for the Earth--Moon system is $\mu = 1.215058560962404 \times 10^{-2}$, following the JPL Three-Body Periodic Orbit Catalog~\cite{JPL_periodic_orbits}.\ The numerical integrations are performed using the adaptive Taylor integrator provided by the \texttt{heyoka.py} library~\cite{biscani_izzo}, with the Taylor order set to 20 and the error tolerance set to machine epsilon.\ All computations are carried out in double precision, except for the collision orbits Nos.~2--4 of Table~\ref{tbl:periodic-collision-orbits} and for the Conley--Zehnder index computation of orbits with high instability in Tables~\ref{tbl:B2}--\ref{tbl:B7}, which are computed in quadruple precision.\ A numerically computed orbit is regarded as periodic when the norm of the vector of components that should vanish by the symmetry constraints is below the tolerance of $10^{-10}$.\

\subsection{Overview of the Obtained Bifurcation Networks}

\begin{figure}[t!]
	\centering
	\includegraphics[width=1.0\textwidth]{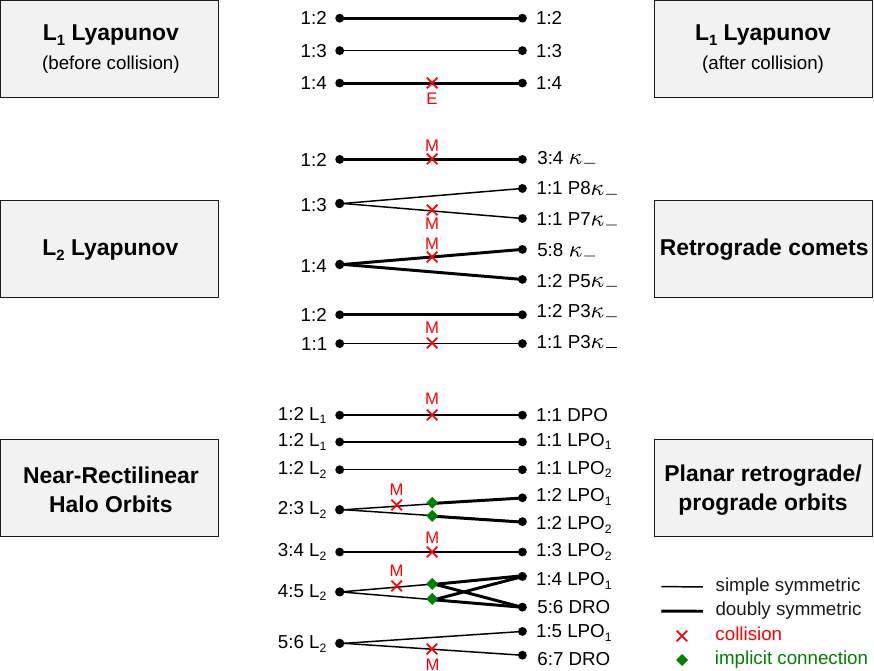}
	\caption{Schematic overview of the bifurcation networks obtained in the Earth--Moon CR3BP.\ The upper part shows the internal bridges between multiple covers of the $L_1$ Lyapunov family; the middle part shows the bridges from the $L_2$ Lyapunov family to planar families associated with the retrograde comet family; and the lower part shows the connections from the $L_1$ and $L_2$ NRHO families to planar prograde and retrograde orbit families.\ The labels $d$:$k$ indicate the corresponding self-resonances.\ Collisions with the Earth and the Moon are indicated by a cross with labels `E' and `M', respectively.\ }
	\label{fig:summary}
\end{figure}

We begin with an overview of the bifurcation network obtained from the numerical continuations, consisting of spatial periodic orbit families linking several well-known planar and spatial families originating from different asymptotic regimes of the CR3BP.\ Figure~\ref{fig:summary} summarizes the main results and serves as a roadmap for the detailed discussion in the following subsections.\ The network separates into three parts:
\begin{itemize}
	\item[(1)] The spatial families bifurcating from the 1:2, 1:3, and 1:4 VSR $L_1$ Lyapunov orbits form \textit{internal bridges}:\ after continuation, one family from each resonance reconnects to another VSR orbit of the corresponding multiple cover of the same $L_1$ Lyapunov family.\ More precisely, these bridges connect $L_1$ Lyapunov orbits with same multiplicities lying before and after the collision of the planar family with the Earth.\ Therefore, the spatial families provide global connections between two regions of the same underlying planar $L_1$ Lyapunov family separated by a collision orbit.\
	\item[(2)] For the $L_2$ Lyapunov orbits, the spatial families generated from its vertical self-resonances do not return to the $L_2$ Lyapunov family, but instead form bridge connections to planar families associated with the retrograde comet family $\kappa_-$ that has been recently studied in~\cite{aydin_comet}.\ Several of these connections pass through collision with the Moon.\
	\item[(3)] Families generated from period-multiplying bifurcations of the NRHOs form a structured network of connections to planar prograde and retrograde families.\ These connections occur either directly or implicitly.\ In the latter case, the bifurcated family connects to a doubly symmetric family at a point of a symmetry-breaking pitchfork bifurcation.\ The results align with the recent work~\cite{aydin_batkhin}, where it is shown for the Hill problem that the rectilinear vertical collision, prograde and retrograde orbit families form a structured bifurcation network organized by their covering multiplicities.\ In particular, there are connections between the $k$-covered rectilinear vertical collision family, the $(k-1)$-covered prograde family, and the $(k+1)$-covered retrograde family, for $k = 4,5,6$.\ At lower multiplicities the double and triple covers of the rectilinear collision family are connected to the simple and double covers of the prograde family, respectively.\ Since the continuation of the rectilinear vertical collision family from the Hill problem to the CR3BP connects to the $L_1$ and $L_2$ NRHO families~\cite{joung_2026_bifurcations}, the present results provide a corresponding structure for lunar orbits of the Earth--Moon CR3BP.\
\end{itemize}

\subsection{Vertical Bifurcations from $L_1$ and $L_2$ Lyapunov orbits}
\label{sec:5.2}
\begin{figure}[t!]
	\centering
	\includegraphics[width=1.0\textwidth]{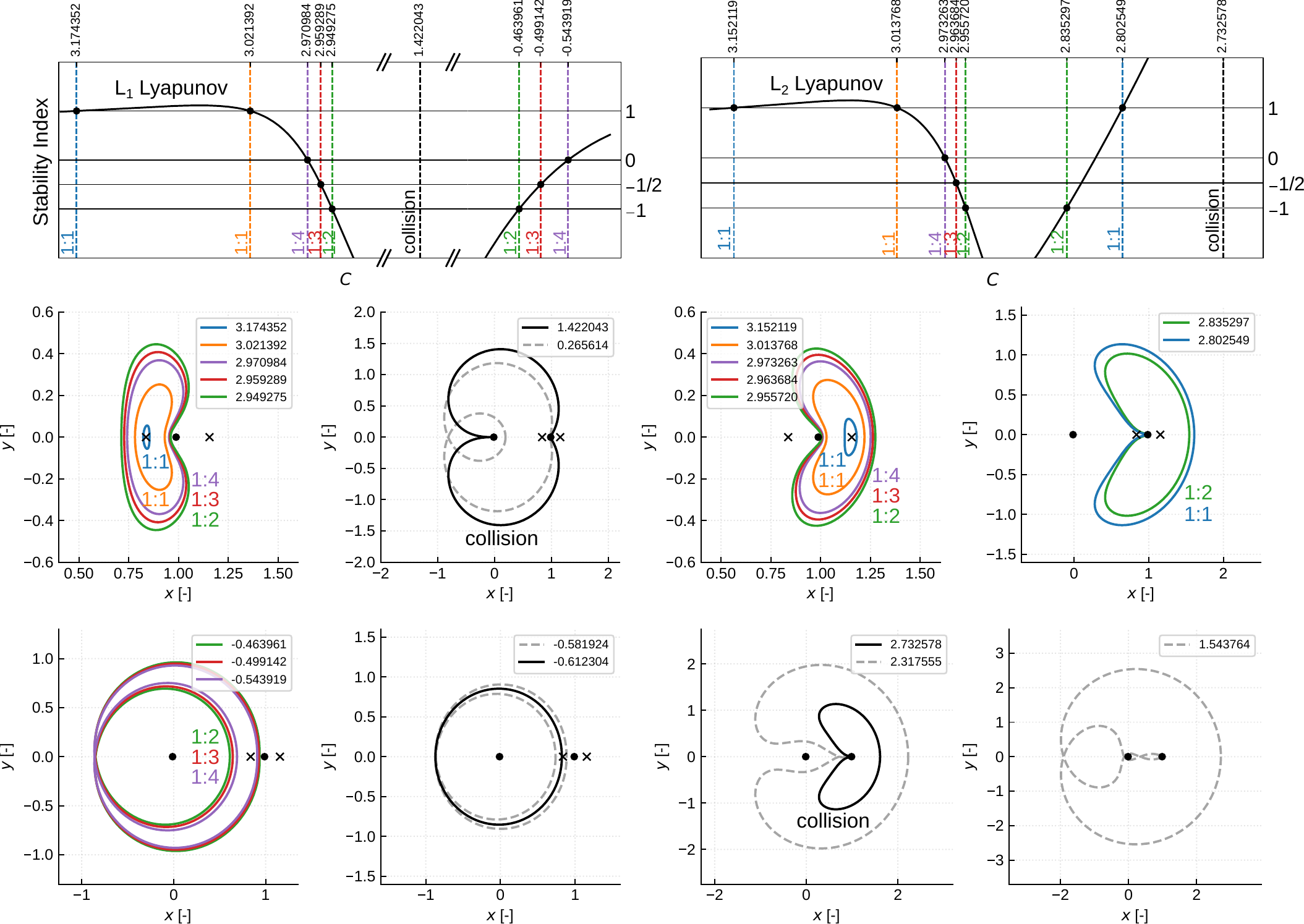}
	\caption{Vertical stability diagrams and orbit plots of the $L_1$ and $L_2$ planar Lyapunov families, split correspondingly into the left and right block.\ On the vertical stability diagrams, the dashed lines indicate Jacobi constants at which self-resonances or collisions occur.\ The corresponding orbit plots below are arranged in two rows, ordered from top left to top right, then bottom left to bottom right, labeled by their Jacobi constants.\ The axis units [-] indicate Earth--Moon distance units.\ }
	\label{fig:L1-L2}
\end{figure}

The first initial guess for the differential correction process for generating the planar Lyapunov orbit families is provided by the in-plane frequency that is associated to the linear behavior of the flow around the $L_1$ and $L_2$ libration points \cite{broucke,szebehely}.\ Near $L_1$ and $L_2$, respectively, the orbits are approximately elliptical and surround the corresponding libration point in the retrograde sense.\ As the families are continued away from the equilibria, their amplitudes increase, and the orbits undergo collision with the primaries and a sequence of stability transitions and vertical self-resonances.\ Figure~\ref{fig:L1-L2} shows the evolution of the $L_1$ and $L_2$ Lyapunov families together with their vertical stability indices and representative orbits, whose numerical data is provided in Table~\ref{tbl:B1}.\ Along the $L_1$ Lyapunov family, continuation leads to a collision with the Earth at $C \approx 1.422043$, while the $L_2$ Lyapunov family reaches collision with the Moon at $C \approx 2.732578$.\ The $L_1$ Lyapunov family terminates at a period-doubling bifurcation point of the planar retrograde family around the Earth (family $A_1$ in the notation of Broucke~\cite{broucke}) at $C \approx -0.612303$.\ In contrast, after passing through collision with the Moon, the $L_2$ Lyapunov family can be continued further into the positive hyperbolic regime, with increasing stability indices.\

The VSR orbits detected along the $L_1$ and $L_2$ Lyapunov families up to
multiplicity four are marked in Figure~\ref{fig:L1-L2} by their self-resonance
ratios $d$:$k$.\ The first two 1:1 VSR orbits along the Lyapunov families correspond to classical bifurcations that have been extensively studied.\ The first gives rise to the halo families~\cite{farquhar_kamel,breakwell_brown,howell_breakwell}, while the second gives rise to the axial families~\cite{doedel}.\ The latter form bridges to the vertical Lyapunov families.\ These bifurcations provide useful reference points in the stability diagrams, but they are not the main focus of the present work.\ Instead, we concentrate on the higher-period and near-collision vertical bifurcations.\ The VSR bifurcations considered naturally divide into two classes.\ The first consists of the period-doubling, period-tripling, and period-quadrupling bifurcations associated with the 1:2, 1:3, and 1:4 VSR orbits of the $L_1$ and $L_2$ Lyapunov families.\ These bifurcations are discussed in Sections~\ref{sec:5.2.1} and~\ref{sec:5.2.2}, respectively.\ The second class consists of the third 1:1 and the second 1:2 VSR orbits along the $L_2$ Lyapunov family shortly before its collision with the Moon.\ These near-collision bifurcations are specific to the $L_2$ family and do not occur in the limiting Hill problem~\cite{henon_1969}.\ These bifurcations are discussed in Section~\ref{sec:5.2.2}.\

Before presenting the individual bifurcations, we briefly summarize the graphical conventions used in the bifurcation diagrams below, following those in~\cite{aydin_batkhin,aydin_comet,aydin_cz,aydin_dro,moreno_aydin}:\ The vertices are bifurcation points and the edges are the orbit families, labeled with Conley--Zehnder index (see Appendix~\ref{appendix_a} for a brief description of the Conley--Zehnder index and its behavior at bifurcations).\ The diagrams are oriented from bottom to top in the direction of decreasing Jacobi constant (increasing Hamiltonian energy).\ Thin and thick edges distinguish families of simple and doubly symmetric orbits, respectively.\ At bifurcation points involving planar orbits, the corresponding planar families are represented locally by short vertical black edges.\ Spatial families generated from vertical bifurcations are shown in red or blue:\ if the bifurcated family initially has the same Conley--Zehnder index as the underlying orbit before the bifurcation, we refer to it as the \textit{red family}; otherwise, we refer to it as the \textit{blue family}.\ Symmetry-related families obtained through the $\sigma$-symmetry are shown in the same colour with dashed edges.\ Collisions with Earth and Moon are indicated by a cross labeled with `E' and `M', respectively.\ These diagrams provide a convenient representation of the bifurcation structure, where the labeling of Conley-Zehnder indices are used to verify invariance of the Euler characteristic of local Floer homology before and after a bifurcation.\ Note that bad orbits are indicated by an upper star on their Conley--Zehnder index.\ 

\subsubsection{Internal Bridges of the $L_1$ Lyapunov Family}
\label{sec:5.2.1}

\begin{figure}[t!]
    \centering
	\includegraphics[width=\textwidth]{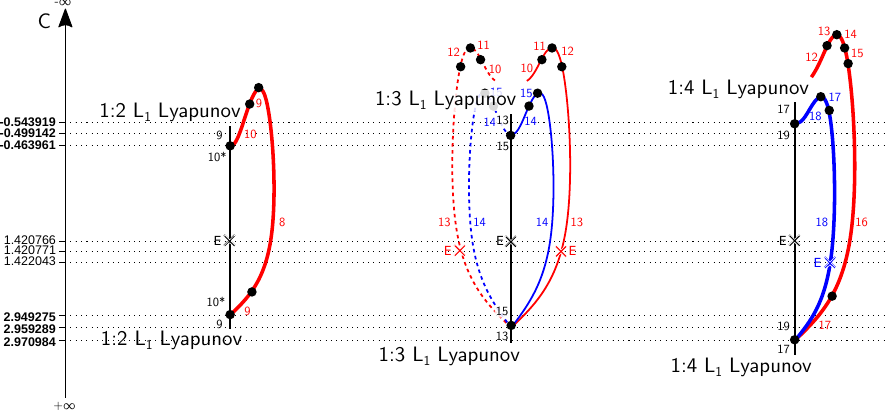}
    \caption{Bifurcation diagrams for the period multiplying vertical bifurcations of the $L_1$ Lyapunov family. The diagrams are labeled by the Conley--Zehnder indices, and the collisions with the Earth are marked by a cross. Note that there is a second bifurcated branch for each of the vertical bifurcations occurring at $C\approx -0.499142$ and $C\approx -0.543919$, which are not indicated in this diagram.}
    \label{fig:bif-L1}
\end{figure}

\begin{figure}[t!]
	\includegraphics[width=\textwidth]{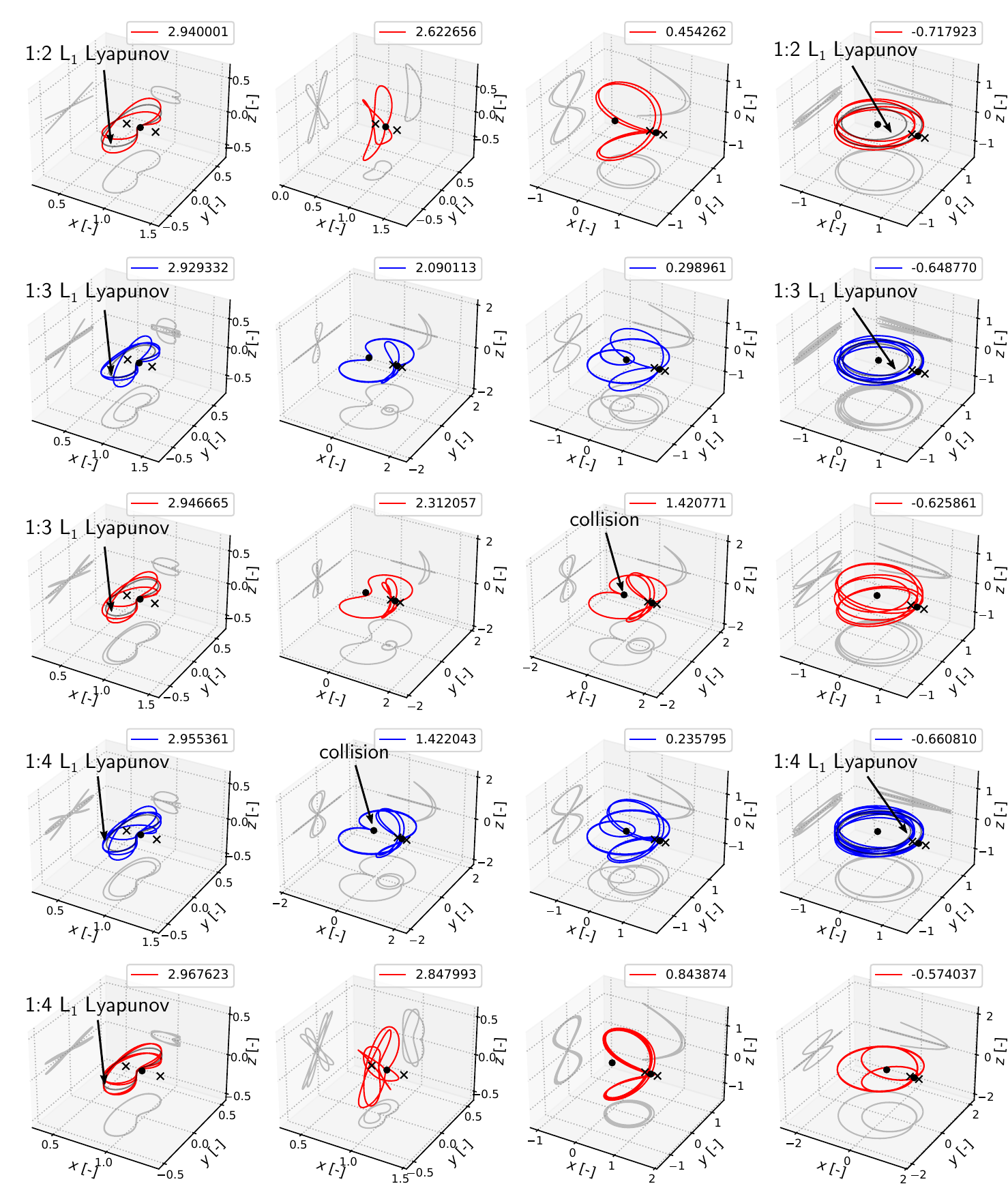}
	\caption{Representative orbits of the families emerging from the VSR bifurcations of the $L_1$ Lyapunov family.}
	\label{fig:orbits-L1}
\end{figure}

The bifurcation graphs for the spatial families branching out from the 1:2, 1:3, and 1:4 $L_1$ Lyapunov orbits are collected in Figure~\ref{fig:bif-L1} with representative orbits shown in Figure~\ref{fig:orbits-L1}.\
The numerical data are provided in Tables~\ref{tbl:B2},~\ref{tbl:B3}, and~\ref{tbl:B4}, respectively.\ For each self-resonance, one of the bifurcated families forms an \textit{internal bridge} that reconnects to another $L_1$ Lyapunov orbit with the same vertical self-resonance.\ 

The first internal bridge originates from the 1:2 VSR orbit at $C \approx 2.949275$.\ At this period-doubling bifurcation, the double cover of the planar $L_1$ Lyapunov orbit changes its Conley--Zehnder index from $9$ to $10$.\ The bifurcated doubly symmetric red family initially has Conley--Zehnder index $9$ and forms a bridge back to the double cover of the $L_1$ Lyapunov family at $C \approx -0.463961$, where the Conley--Zehnder index of the planar family jumps from $10$ back to $9$.\ Along continuation, the spatial family undergoes several stability transitions including a birth--death bifurcation.\ 

\begin{figure}[t!]
	\centering
	\includegraphics[width=\textwidth]{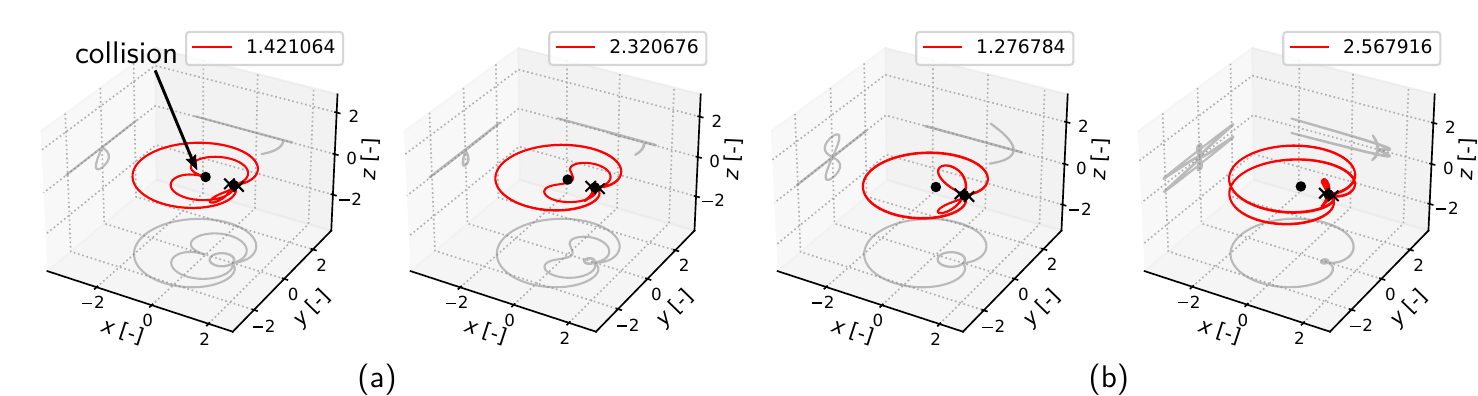}
	\caption{Representative orbits obtained by continuing the red families beyond what is indicated in the bifurcation diagrams of Figure~\ref{fig:bif-L1}. Subfigures (a) and (b) correspond respectively to the branches from the period-tripling and period-quadrupling VSR bifurcations of the $L_1$ Lyapunov family.}
	\label{fig:bif-graph-L1-unstable}
\end{figure}

At the 1:3 VSR orbit at $C \approx 2.959289$, the triple cover of the planar $L_1$ Lyapunov orbit undergoes a Conley--Zehnder index jump from $13$ to $15$, giving rise to two simple symmetric spatial families with different symmetries, as well as their $\sigma$-symmetric counterparts.\ The red family starts with Conley--Zehnder index 13 and consists of $\rho_{xz}$-symmetric orbits, while the blue family starts with Conley--Zehnder index 14 and consists of $\rho_{x}$-symmetric orbits.\ The blue family forms an internal bridge.\ Along its continuation, it passes through a birth--death bifurcation at $C \approx -0.723964$ and reconnects to another 1:3 $L_1$ Lyapunov orbit at $C\approx-0.499142$, where the Conley--Zehnder index of the planar family jumps from $15$ back to $13$.\ The red family passes through collision with the Earth at $C \approx 1.420771$, and goes through a birth--death bifurcation at $C\approx -0.862931$.\ Beyond what is indicated in Figure~\ref{fig:bif-L1}, we were able to continue the red family up to approximately $C \approx 2.320676$, where the orbits remain highly unstable along the continuation.\ Representative orbits from this segment are shown in Figure~\ref{fig:bif-graph-L1-unstable}.\ No connection to a planar family was detected in this continuation range, but we report that the family passes through another collision with the Earth near $C \approx 1.421064$.\ 

At the 1:4 VSR orbit at $C \approx 2.970984$, the Conley--Zehnder index of the planar family jumps from $17$ to $19$, and two doubly symmetric families are generated.\ The blue family passes through collision with the Earth at $C \approx 1.420766$, then forms an internal bridge by reconnecting to another 1:4 $L_1$ Lyapunov orbit at $C \approx -0.543919$.\ The orbits of the red family go through a series of bifurcations at high energy, with Conley--Zehnder index dropping sequentially from 16 to 12, before continuing in decreasing energy direction with high instability.\ Beyond what is indicated in Figure~\ref{fig:bif-L1}, we were able to continue the red family up to approximately $C \approx 2.567916$.\ We did not detect a connection to a planar family in the computed range.\ Representative orbits from this highly unstable segment are shown in Figure~\ref{fig:bif-graph-L1-unstable}.\

\subsubsection{Bridges from the $L_2$ Lyapunov Family to Retrograde Comet Orbits}
\label{sec:5.2.2}

\begin{figure}[t!]
    \centering
	\includegraphics[width=\textwidth]{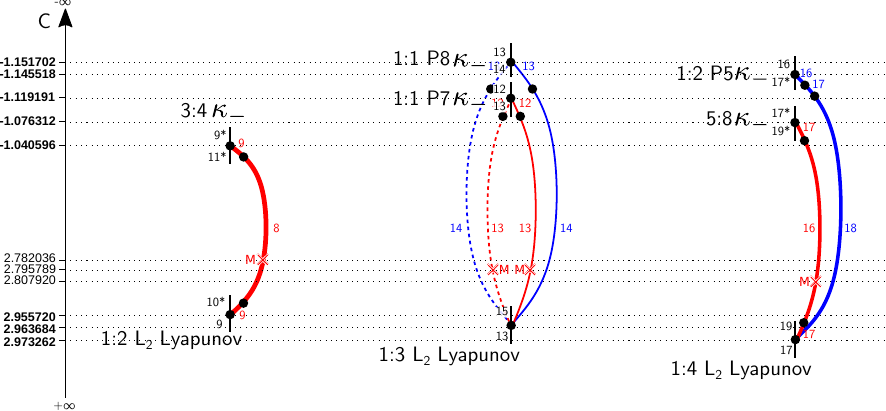}
    \caption{Bifurcation diagrams for the period multiplying vertical bifurcations of the $L_2$ Lyapunov family. The diagrams are labeled by the Conley--Zehnder indices, and the collisions with the Moon are marked by a cross. Note that there is a second bifurcated branch for each of the vertical bifurcations occurring at $C\approx -1.040596$, $C\approx -1.076312$, $C\approx -1.119191$, $C\approx -1.145518$, and $C\approx -1.151702$, which are not indicated in this diagram.}
    \label{fig:bif-L2}
\end{figure}

\begin{figure}[t!]
	\includegraphics[width=\textwidth]{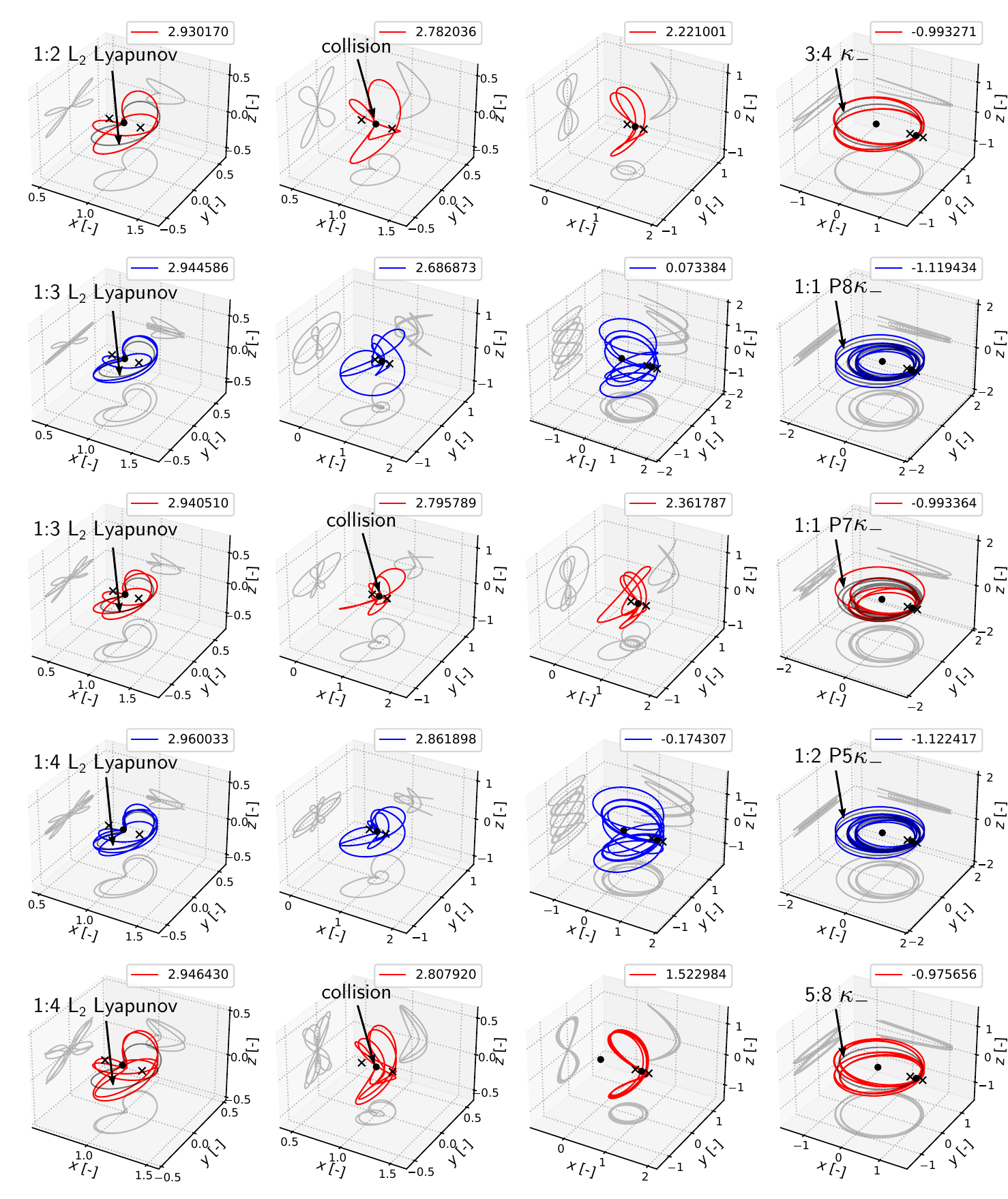}
	\caption{Representative orbits of the families emerging from the VSR bifurcations of the $L_2$ Lyapunov family.}
	\label{fig:orbits-L2}
\end{figure}

\begin{figure}[t!]
    \centering
	\includegraphics[width=.55\textwidth]{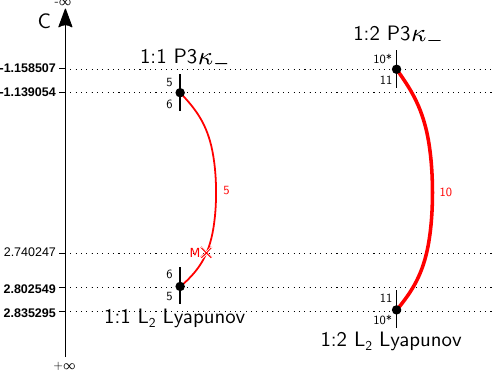}
    \caption{Bifurcation diagrams for the vertical pitchfork (left) and period-doubling (right) bifurcations of the $L_2$ Lyapunov family occurring near collision. The diagrams are labeled by the Conley--Zehnder indices, and collision with the Moon is marked by a cross.}
    \label{fig:bif-L2-second}
\end{figure}

\begin{figure}[t!]
	\includegraphics[width=\textwidth]{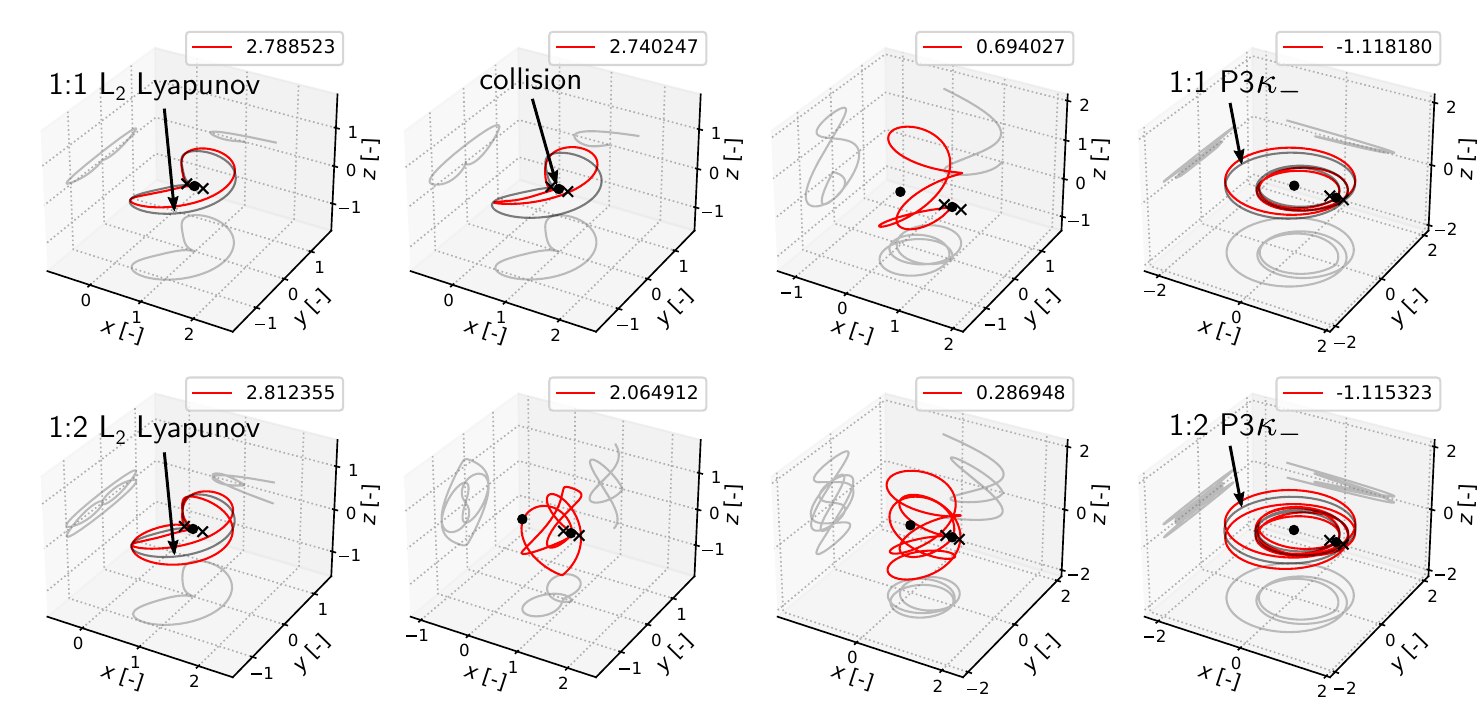}
	\caption{Representative orbits of the families emerging from the near-collision vertical bifurcations of the $L_2$ Lyapunov family.}
	\label{fig:orbits-L2-second}
\end{figure}

In contrast to the internal bridges arising from the $L_1$ Lyapunov family, the spatial families generated from the $L_2$ Lyapunov family do not reconnect to another region of the same planar family.\ Instead, their continuation reveals spatial bridges to the planar retrograde comet family $\kappa_-$ that has been recently studied in~\cite{aydin_comet}, or to planar families that bifurcate from multiple covers of the $\kappa_-$ comet family.\ For convenience, we denote by P$m$$\kappa_-$ a planar family bifurcating from the $m$-fold cover of the retrograde comet family $\kappa_-$.\ The resulting bifurcation diagrams are shown in Figures~\ref{fig:bif-L2} and~\ref{fig:bif-L2-second}, where Figure~\ref{fig:bif-L2} describes the families generated from the first $1$:$2$, $1$:$3$, and $1$:$4$ self-resonances, and Figure~\ref{fig:bif-L2-second} from the $1$:$1$ and $1$:$2$ self-resonances which are specific to the $L_2$ Lyapunov family and occur near collision.\ The representative orbits are shown in Figures~\ref{fig:orbits-L2} and~\ref{fig:orbits-L2-second}, respectively, and the numerical data are provided in Tables~\ref{tbl:B5}--\ref{tbl:B8}.\ Unlike the branches generated from the $L_1$ Lyapunov family, these spatial bridges encounter no birth--death bifurcation and are therefore monotone with respect to the Jacobi constant.\ Several of these spatial bridges pass through collision with the Moon.\ 

The first period-doubling bifurcation of the $L_2$ Lyapunov family occurs at the $1$:$2$ VSR orbit at $C \approx 2.955720$.\ The resulting doubly symmetric spatial family forms a bridge from the double cover of an $L_2$ Lyapunov orbit to the $4$-fold cover of the retrograde comet orbit $\kappa_-$ at $C \approx -1.040596$.\ The bifurcated family initially consists of elliptic--positive hyperbolic orbits with Conley--Zehnder index $9$.\ At $C \approx 2.939704$, it undergoes a stability transition to the positive--positive hyperbolic regime, as the Conley--Zehnder index jumps to $8$.\ The family passes through collision with the Moon at $C \approx 2.782036$.\ After a further bifurcation at $C \approx -1.019989$, where the index returns to $9$, the family connects with $3$:$4$ self-resonant $\kappa_-$ at $C \approx -1.040596$.\ 

At the $1$:$3$ VSR orbit at $C \approx 2.963684$, the Conley--Zehnder index of the $3$-fold cover of the $L_2$ Lyapunov orbit jumps from $13$ to $15$, and two simple symmetric spatial families are generated, together with their $\sigma$-symmetric counterparts.\ The red family is $\rho_{xz}$-symmetric, while the blue family is $\rho_{x}$-symmetric.\ Both families form connections to planar families associated with retrograde comet family $\kappa_-$.\ Namely, the families connect to vertical pitchfork bifurcations of the planar families P$7$$\kappa_-$ and P$8$$\kappa_-$, which originate from planar bifurcations of the $7$-fold and $8$-fold cover of the retrograde comet family $\kappa_-$, respectively.\ The red family passes through collision with the Moon at $C \approx 2.795789$.\

At the $1$:$4$ VSR orbit at $C \approx 2.973262$, the Conley--Zehnder index of the $4$-fold cover of the $L_2$ Lyapunov orbit jumps from $17$ to $19$, and two doubly symmetric spatial families are generated.\ The red family goes through collision with the Moon at $C\approx 2.807920$, then terminates at the $5$:$8$ VSR retrograde comet orbit $\kappa_-$ at $C \approx -1.076312$, whereas the blue family connects to the double cover of the planar family P$5$$\kappa_-$ at $C \approx -1.145518$.\ 

The $L_2$ Lyapunov family contains two further vertical bifurcations shortly before its collision with the Moon at $C \approx 2.732578$ (see Figure~\ref{fig:L1-L2}), namely at the second 1:2 VSR orbit at $C \approx 2.835295$ and at the third $1$:$1$ VSR orbit at $C \approx 2.802549$.\ The resulting bifurcation diagrams are given in Figure~\ref{fig:bif-L2-second}, and representative orbits are shown in Figure~\ref{fig:orbits-L2-second}.\ The vertical pitchfork bifurcation at the $1$:$1$ VSR orbit gives rise to a $\rho_{xz}$-symmetric family which continues with constant Conley--Zehnder index $5$.\ It passes through collision with the Moon at $C \approx 2.740247$ and connects to the vertical pitchfork bifurcation of P$3$$\kappa_-$ at $C \approx -1.139054$.\ The spatial family generated from the near-collision $1$:$2$ VSR orbit is doubly symmetric and has constant Conley--Zehnder index $10$ along the continuation.\ This family connects to a period-doubling bifurcation of P$3$$\kappa_-$ at $C \approx -1.158507$.\

Combining the results of Sections~\ref{sec:5.2.1} and~\ref{sec:5.2.2}, a total of 9 symmetric periodic collision orbits are identified along the continued families, including the two from the planar Lyapunov families.\ The collisions are identified by verifying a sign change in the $\xi_2$-coordinate along the continuation of the $\tilde{\rho}_{xz}$-symmetric point in the Moser regularization (see Proposition~\ref{prop:symmetric-collision-configurations}).\ To precisely locate these collision orbits, we additionally apply the Newton correction scheme of Lemma~\ref{lem:collision_newton} with $\xi_1 = 0$ as the Poincar\'e section. To compute derivatives of the Poincar\'e map with respect to energy $c$, we augment the ODE with $\dot{c} = 0$ and solve the variational equations of the augmented system.\ The resulting numerical data is shown in Table~\ref{tbl:periodic-collision-orbits}.\ For orbits Nos.~2--4, convergence of the correction scheme was not achieved in double precision arithmetic, and was instead obtained using extended precision arithmetic.\ We report that the Jacobi constants of the collision orbits are clustered around similar values, depending on whether they originate from the $L_1$ or $L_2$ Lyapunov family.\ 

\begin{table}[t]
\setlength{\tabcolsep}{5pt}
\centering
\small
\caption{Symmetric periodic collision orbits found along families related to the planar $L_1$ and $L_2$ Lyapunov families, computed using the correction scheme of Lemma~\ref{lem:collision_newton}. Here, $C$ is the Jacobi constant, $k$ is the number of iterates of the Poincar\'e map $\mathcal{P}_{c,\pm}$ on the section $\xi_1=0$, and ``dir.'' denotes the crossing direction. The column $\eta_3$ gives the sign of $\eta_3=\pm\sqrt{m_P^2-\eta_1^2}$ in the case of spatial collision orbits, where $m_P = \mu$ for collisions with the Moon, and $m_P = 1 - \mu$ for collisions with the Earth.\ }
\label{tbl:periodic-collision-orbits}
\begin{tabular}{cclccccccc}
\hline
No. & Type & Origin & Collision & Symmetry & $C$ & $\eta_1$ & $\eta_3$ & $k$ & dir. \\
\hline
1 & planar  & $L_1$ Lyapunov family 	& Earth & $\tilde{\rho}$ 				 	   & $1.4220430992$ & $-(1 - \mu)$ & $0$ & $2$ & $+$ \\
2 & spatial & $1\!:\!3$ $L_1$ Lyapunov 	& Earth & $\tilde{\rho}_{xz}$ 				   & $1.4207713748$ & $-0.9877845625$ & $-$ & $5$ & $+$ \\
3 & spatial & $1\!:\!3$ $L_1$ Lyapunov 	& Earth & $\tilde{\rho}_{xz}$ 				   & $1.4210644051$ & $-0.9878494140$ & $-$ & $4$ & $-$ \\
4 & spatial & $1\!:\!4$ $L_1$ Lyapunov 	& Earth & $\tilde{\rho}_{xz}$-$\tilde{\rho}_x$ & $1.4207660622$ & $-0.9877842881$ & $+$ & $3$ & $-$ \\
5 & planar  & $L_2$ Lyapunov family 	& Moon  & $\tilde{\rho}$ 				 	   & $2.7325778911$ & $\mu$ & $0$ & $1$ & $-$ \\
6 & spatial & $1\!:\!2$ $L_2$ Lyapunov 	& Moon  & $\tilde{\rho}_{xz}$-$\tilde{\rho}_x$ & $2.7820369902$ & $\phantom{+}0.0046062245$ & $-$ & $1$ & $+$ \\
7 & spatial & $1\!:\!3$ $L_2$ Lyapunov 	& Moon  & $\tilde{\rho}_{xz}$ 				   & $2.7957888525$ & $\phantom{+}0.0042680552$ & $+$ & $3$ & $-$ \\
8 & spatial & $1\!:\!4$ $L_2$ Lyapunov 	& Moon  & $\tilde{\rho}_{xz}$-$\tilde{\rho}_x$ & $2.8079202799$ & $\phantom{+}0.0040093290$ & $-$ & $2$ & $+$ \\
9 & spatial & $1\!:\!1$ $L_2$ Lyapunov 	& Moon  & $\tilde{\rho}_{xz}$ 				   & $2.7402468921$ & $\phantom{+}0.0106918628$ & $+$ & $1$ & $-$ \\
\hline 
\end{tabular}
\end{table}

\subsection{Bifurcation Network connecting NRHOs with Planar Families}
\label{sec:nrho-networks}

\begin{figure}[!t]
    \centering
	\includegraphics[width=1.0\textwidth]{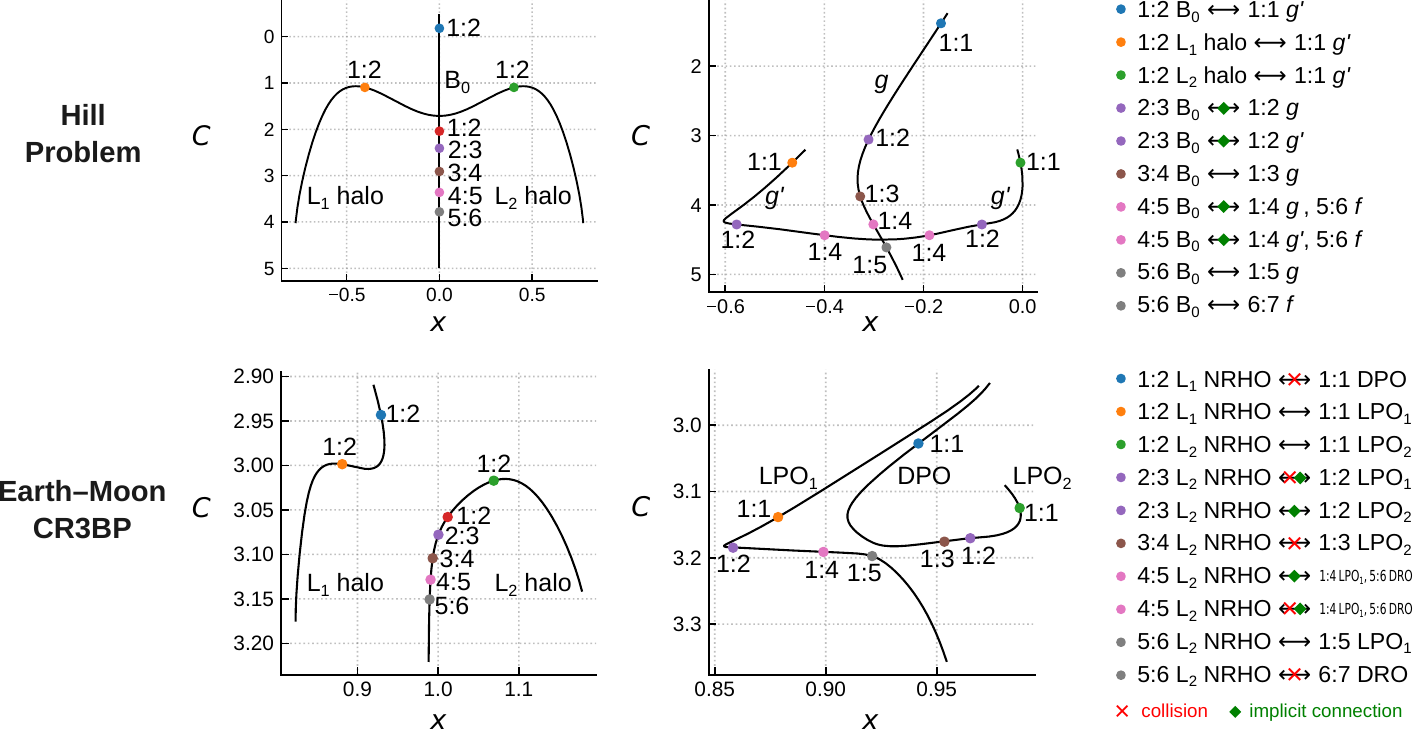}
	\caption{Comparison of the $(x, C)$-bifurcation graphs of the $L_1$ and $L_2$ halo families (left) and planar retrograde and prograde orbit families (right), for the Hill problem (top) and the Earth--Moon CR3BP (bottom).\ The markers on the bifurcation graphs indicate self-resonant points involved in connections through bifurcated branches, summarized on the right.\ Branches passing through periodic collision orbits with the Moon are indicated with a red cross.}
    \label{fig:network}
\end{figure}

In this section, we describe the bifurcation network involving Near-Rectilinear Halo Orbits (NRHOs) and planar prograde and retrograde orbit families in the Earth--Moon CR3BP.\ The motivation comes from the corresponding structure in the Hill problem, which is a limiting approximation of the CR3BP describing the dynamics in a neighbourhood of the light primary when the mass ratio tends to zero.\ In~\cite{aydin_batkhin}, it was shown that the rectilinear vertical collision family $\mathcal{B}_0^{\pm}$, the planar prograde families $g$ and $g'$, and the planar retrograde family $f$ are systematically linked through spatial bifurcated families, with the connections organized by their covering multiplicities and symmetry properties.\ The two prograde families are themselves related, with $g'$ emerging from $g$ at a symmetry-breaking pitchfork bifurcation~\cite{henon_1969}.\ The perturbation of these families from the Hill problem to the CR3BP provides the natural setting to study the corresponding bifurcations.\ In particular, the Hill problem possesses additional symmetries (see~\cite{aydin_symmetries}) that are not preserved under perturbation to the CR3BP, resulting in the breaking of the associated pitchfork bifurcations.\ The rectilinear vertical collision family extends to the $L_1$ and $L_2$ halo and NRHO families~\cite{joung_2026_bifurcations}, while the planar prograde families $g$ and $g'$ extend to the families denoted, following~\cite{restrepo_russell}, by DPO, $LPO_1$, and $LPO_2$ (see Figure~\ref{fig:network}).\ Here, $LPO_1$ (Low Prograde Orbit) denotes the family of small prograde periodic orbits around the light primary, while $LPO_2$ forms a second planar prograde family connected to the DPO (Distant Prograde Orbit) family through a birth--death bifurcation.\ The planar retrograde family extends to the familiar Distant Retrograde Orbit (DRO) family.\ We demonstrate below that the systematic bifurcation structure from the Hill problem persists in the Earth--Moon CR3BP, as illustrated in Figure~\ref{fig:network}.\

The $L_1$ and $L_2$ halo families are the spatial families that bifurcate out-of-plane from the corresponding planar Lyapunov families~\cite{breakwell_brown,howell}.\ Along continuation, both halo families undergo a birth--death bifurcation, after which the orbits acquire the characteristic near-rectilinear geometry.\ This part of the family is commonly referred to as the NRHO family~\cite{howell_breakwell}.\ At this transition, the halo orbits become elliptic--elliptic for a small range of Jacobi energy, until they become unstable after a subsequent period-doubling bifurcation.\ Further along the continuation, the $L_1$ and $L_2$ halo families transition to a family of highly vertical, near-collision orbits.\ For the $L_1$ halo family, the orbits undergo a second birth--death bifurcation and continue toward increasing energy, while the $L_2$ halo family transitions back to elliptic--elliptic type after another period-doubling bifurcation (see Figure~\ref{fig:network}).\ We notice that the further continuation of the $L_1$ halo family terminates at a planar orbit, at a vertical pitchfork bifurcation of the retrograde comet family~\cite{doedel,aydin_comet}.\ Bifurcations of halo and NRHO families are of great interest in mission design~\cite{campbell,zimovan_spreen}.\ In recent studies, a class of families generated by period-multiplying bifurcations of the $L_2$ NRHOs has been referred to as tulip-shaped orbits~\cite{koblick_kelly,zhang} due to their characteristic geometry.\ These families provide a range of spatial configurations with varying inclination, altitudes, and ground-tracks.\ From a dynamical point of view, a natural question is which parts of the bifurcation network observed in the Hill problem persist in the Earth-Moon system, and how the corresponding connections are modified upon perturbation.\

The connections in the Earth--Moon CR3BP are summarized in Figure~\ref{fig:network}.\ Families generated from period-multiplying bifurcations of the $L_1$ and $L_2$ NRHO families up to bifurcations of the sixth cover are computed.\ A single family is generated from each of the four period-doubling bifurcations, while two families are generated from each of the four higher-period bifurcations.\ In total, this gives twelve bifurcated families.\ Each family is continued until it either reaches a planar periodic orbit, or connects to a doubly symmetric family which itself connects to a planar orbit.\ We refer to the latter case as an \emph{implicit connection}.\ Among the twelve families, ten are found to form either direct or implicit connections to a planar family.\ For the remaining two families, generated from the $1$:$2$ and $3$:$4$ self-resonances of the $L_2$ NRHO family, no explicit connections were found.\ The numerical data for the self-resonant orbits corresponding to endpoints of these connections, as well as the critical orbits involved in the implicit connections, are provided in Table~\ref{tbl:nrho-critical-orbits}.\

The results exhibit the following structure with respect to the covering multiplicities.\ The two $1$:$2$ bifurcations of the $L_1$ NRHO family connect to vertical pitchfork bifurcations of the DPO and $LPO_1$ families, respectively.\ For the $L_2$ NRHO family, the bifurcations associated with $k$-covered $L_2$ NRHOs are connected to $(k-1)$-covered planar prograde families for $k=2,...,6$, while for $k=5,6$ additional connections to the corresponding $(k+1)$-covered DROs are obtained.\ The connections from $k$-covered $L_2$ NRHOs for $k=3, 5$ are implicit, while the other connections are direct.\ This structure of the bifurcation network aligns with the structure observed in the Hill problem~\cite{aydin_batkhin}.\ 

\begin{figure}[t]
    \centering
    \includegraphics[width=0.6\textwidth]{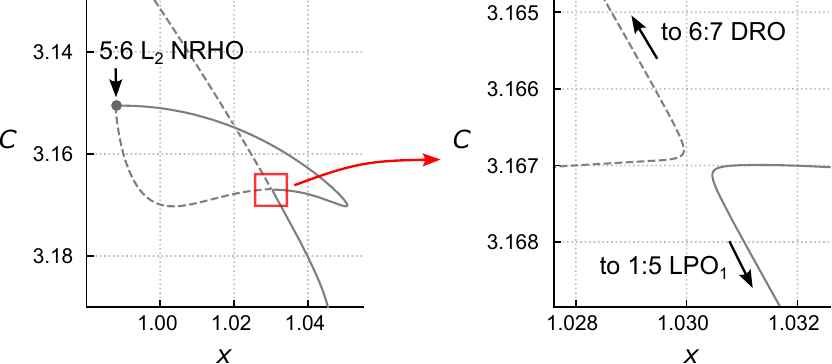}\hspace{2em}
    \caption{Branches generated from the 5:6 NRHO in the $(x,C)$ bifurcation diagram.\ The zoomed out view (left) shows the two branches approaching closely, while the zoomed in view (right) shows them as separated branches with no direct connection.\ The branch with its symmetric point at apoapsis (dashed) continues toward decreasing $C$ and connects to the 6:7 DRO, whereas the branch with its symmetric point at periapsis (solid) continues toward increasing $C$ and connects to the 1:5 $LPO_1$.\ }
    \label{fig:near-bif}
\end{figure}

We report that the two branches generated from the $5$:$6$ $L_2$ NRHO approach each other closely in the $(x,C)$ bifurcation diagram (Figure~\ref{fig:near-bif}), while a zoomed in view shows the two branches as distinct, with no direct connection between them.\ The branch whose symmetric point lies at apoapsis turns toward decreasing Jacobi constant and connects to the 6:7 DRO, whereas the branch whose symmetric point lies at periapsis turns toward increasing Jacobi constant and connects to the 1:5 $LPO_1$.\

\begin{figure}[!t]
    \centering
	\includegraphics[width=1.0\textwidth]{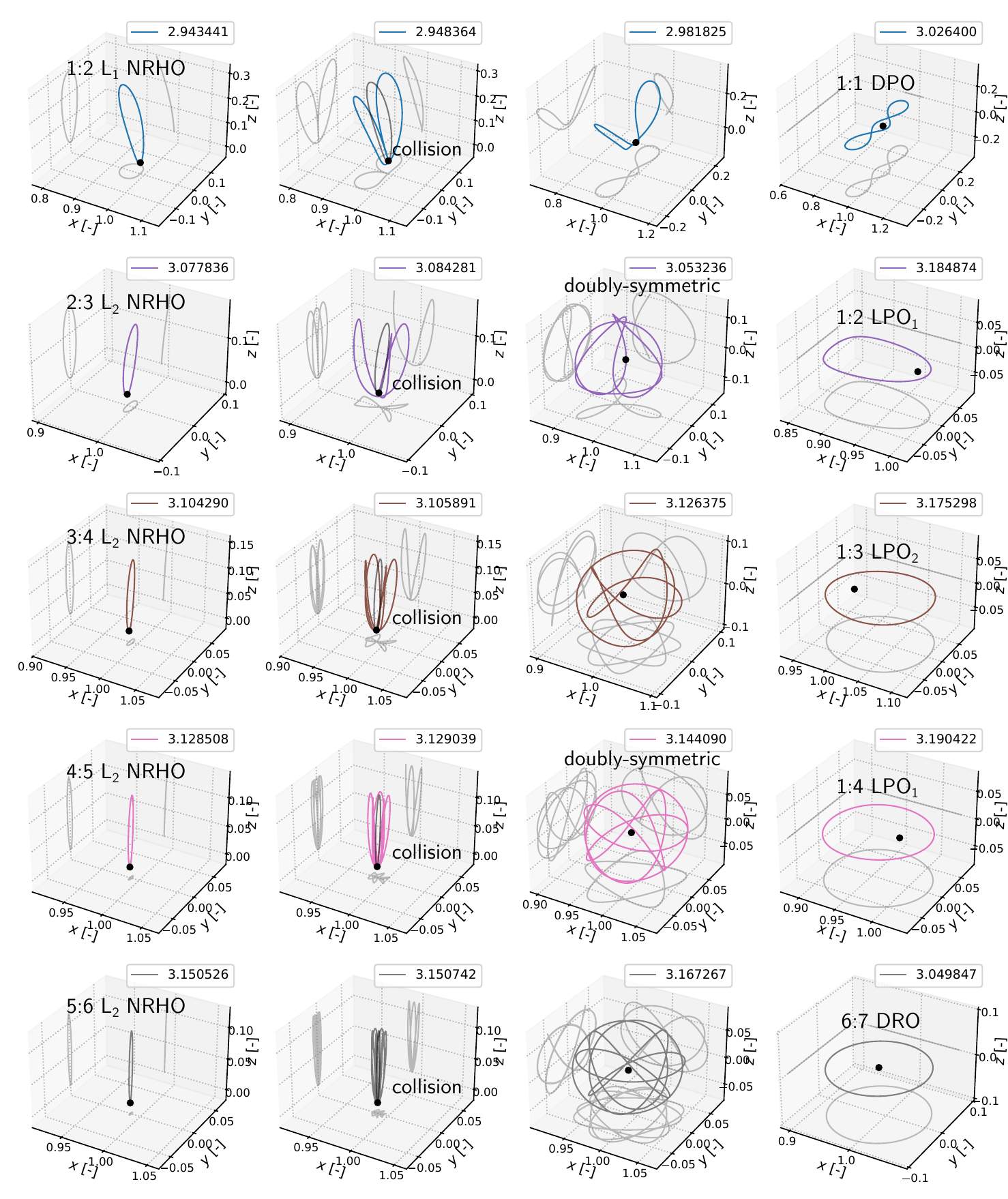}
	\caption{Orbits along the five connecting families from multiple covers of NRHO to planar orbits which go through collision with the Moon. The five orbits on the second column show periodic collision orbits with the Moon which are summarized in Table~\ref{tbl:halo-periodic-collision-orbits}.\ The legends indicate the Jacobi constants of the orbits.}
    \label{fig:network-orbits}
\end{figure}

Five of the connecting branches pass through collision with the Moon.\ We apply Lemma~\ref{lem:collision_newton} to precisely locate the collision orbits, with $\xi_3 = 0$ as the Poincar\'e section. These collision orbits are shown in Figure~\ref{fig:network-orbits}, and their numerical data, including stability indices, are given in Table~\ref{tbl:halo-periodic-collision-orbits}.\ Compared with the periodic collision orbits related to the planar Lyapunov families (Table~\ref{tbl:periodic-collision-orbits}), these examples exhibit better stability as indicated by their small stability indices.\ 

\begin{table}[!ht]
\setlength{\tabcolsep}{5pt}
\centering
\small
\caption{Periodic collision orbits appearing along continuations of period-multiplying bifurcations of NRHOs.\ Here $C$ is the Jacobi constant, $k$ is the number of iterates of the Poincar\'e map $\mathcal{P}_{c,\pm}$ on the section $\xi_3=0$, and ``dir.'' denotes the crossing direction.\ The column $\eta_3$ gives the sign of $\eta_3=\pm\sqrt{\mu^2-\eta_1^2}$.}
\label{tbl:halo-periodic-collision-orbits}
\begin{tabular}{cccccccc}
\hline
Origin & Collision & $C$ & $\eta_1$ & $\eta_3$ & $k$ & dir. & Stability index \\
\hline
1:2 $L_1$ NRHO & Moon & $2.94836431$ & $\phantom{+}0.0028624125$ & $-$ & 1 & $-$ & $s_1\!=\!1.136,\, s_2\!=\!-5.254$ \\
2:3 $L_2$ NRHO & Moon & $3.08428108$ & $          -0.0020506482$ & $-$ & 2 & $+$ & $s_1\!=\!1.013,\, s_2\!=\!-0.2713$ \\
3:4 $L_2$ NRHO & Moon & $3.10589063$ & $          -0.0015068301$ & $-$ & 2 & $-$ & $s_1\!=\!0.9988,\, s_2\!=\!-0.6447$ \\
4:5 $L_2$ NRHO & Moon & $3.12903856$ & $          -0.0010117363$ & $-$ & 3 & $+$ & $s_1\!=\!1.000,\, s_2\!=\!-0.8880$ \\
5:6 $L_2$ NRHO & Moon & $3.15074188$ & $          -0.0007239258$ & $-$ & 3 & $-$ & $s_1\!=\!0.999,\, s_2\!=\!-0.9906$ \\
\hline 
\end{tabular}
\end{table}

\section{Conclusion}
In this work, we developed a computational framework based on Moser regularization for the study of symmetric periodic orbit families and periodic collision orbits in the spatial CR3BP.\ The method was applied to investigate the bifurcation networks of lunar periodic orbits in the Earth--Moon system.\ The framework combines symmetry reduced numerical continuation, Newton and multiple shooting corrections, and the computation of Floquet multipliers and stability indices entirely within Moser regularized coordinates.\ We also established a characterization and numerical detection scheme for symmetric periodic collision orbits.\ A central feature of this approach is that the symplectic structure is preserved under the Moser regularization, allowing the stability and bifurcation properties of periodic orbits to be studied consistently, even as the corresponding trajectories approach or pass through collision with a primary.\

The application to the Earth--Moon CR3BP demonstrated bifurcation networks connecting several classical lunar periodic orbit families, as summarized in Figure~\ref{fig:summary}.\ Higher-period vertical bifurcations of the planar $L_1$ and $L_2$ Lyapunov families, with covering multiplicities up to four, generate spatial bridges connecting different regions of the $L_1$ Lyapunov family across its collision with the Earth and connecting the $L_2$ Lyapunov family to planar retrograde comet orbit families, respectively.\ Furthermore, period-multiplying bifurcations of the Near-Rectilinear Halo Orbit (NRHO) families, with covering multiplicities up to six, give rise to a systematic network of connections to planar prograde and retrograde orbit families, reflecting the characteristic organization by covering multiplicities and symmetries as established in the limiting Hill problem~\cite{aydin_batkhin}.\ An important feature of the resulting bifurcation connections is that several spatial families pass through periodic collisions with the Earth or the Moon, showing that collision need not represent an endpoint of a periodic orbit family but may instead occur along an orbit family connecting different regions.\ Together with Conley--Zehnder indices and related symplectic invariants, the regularized framework provides a systematic approach to studying the global organization of periodic orbit families.\

The present study focused on the geometric and bifurcation structure of periodic orbit families, rather than their direct applications to trajectory design.\ The computed families and their connections may nevertheless provide useful starting points for such applications.\ The developed framework suggests several directions for further research.\ The characterization and numerical detection of periodic collision orbits suggest a systematic classification of these orbits and the families containing them.\ Extending the numerical methods to asymmetric periodic orbits and higher-fidelity dynamical models may broaden its applicability.\ Moreover, higher-order variational equations and normal form methods would allow a more detailed local analysis of bifurcations.\ 

\begin{appendix}
\section{Conley--Zehnder Index and Krein Signature}
\label{appendix_a}

The Conley--Zehnder index~\cite{conley_zehnder} associates to Hamiltonian periodic solutions an integer invariant that intuitively measures the winding number of the linearized flow along the orbit.\ In this appendix we briefly discuss its definition~\cite{conley_zehnder,salamon_zehnder} and related invariants, together with the method of computation, which is based on~\cite{moreno_aydin}.\

\subsection*{The Conley--Zehnder index}
First, consider a path of symplectic matrices $\Psi \colon [0,T] \to Sp(2n)$, for $T > 0$, that starts at the identity and is non-degenerate, meaning that $\Psi(T)$ has no eigenvalue equal to $1$.\ Then define the \textit{Maslov cycle} $V$ by
$$ V = \{ M \in Sp(2n) \colon \det (M - \mathrm{id}) = 0 \}, $$
which is the set of symplectic matrices with eigenvalue $1$.\ The Maslov cycle $V$ splits $Sp(2n)$ into two components, namely
$$V_+ = \{ M \in Sp(2n) \colon \det (M - \mathrm{id}) > 0 \},\quad V_- = \{ M \in Sp(2n) \colon \det (M - \mathrm{id}) < 0 \}, $$
to which we choose correspondingly the base points
$$ B_+ = \text{diag}\{ -1,-1,...,-1,-1 \},\quad B_- = \text{diag}\{ 2, 1/2, -1,-1,...,-1,-1 \}. $$
Any symplectic matrix $M \in Sp(2n)$ has a polar decomposition \cite[Chapter 1]{long}, $M = US$, where $U$ is unitary and $S$ is symmetric positive definite.\ The unitary part is extracted by the retraction
$$ \rho \colon Sp(2n) \to U(n) \subset Sp(2n),\quad M \mapsto (MM^T)^{-1/2}M, $$
where $U(n)$ is identified with the matrices $\left(\begin{smallmatrix} A & -B \\ B & A \end{smallmatrix}\right)$, whose complex determinant is $\det_{\mathbb{C}} = \det(A + iB)$.\ Since the path $\Psi$ is non-degenerate, we can extend it to a path $\tilde{\Psi} \colon [0,2T] \to Sp(2n)$ that ends either at $B_+$ or at $B_-$, such that the extension does not cross the Maslov cycle $V$.\ Then, the \textit{Conley--Zehnder index} of the path $\Psi$ is defined by
$$ \mu_{CZ} (\Psi) = \deg \left( t \mapsto \underset{\mathbb{C}}{\det} (\rho \circ \tilde{\Psi}(t))^2 \right), $$
which counts the number of half turns of the complex determinant around the unit circle.\ The index is invariant under continuous deformations of the path that keep its initial point at the identity and its endpoint non-degenerate.

To connect this construction to non-degenerate Hamiltonian periodic orbits, we choose a symplectic frame along the orbit, i.e., a symplectic basis of the tangent space of the phase space at each point of the orbit.\ The matrix representation of the linearized flow in this frame then gives a path of symplectic matrices, to which the above construction is applied.\ Since the Hamiltonian of the CR3BP is time-independent, the linearized flow always has the eigenvalue $1$ twice, due to the flow direction and the conservation of energy.\ We therefore use a frame that is transverse to the Hamiltonian vector field on the energy level set, and the resulting index is called the \textit{transverse Conley--Zehnder index}.
For our computations, we use the frame $\langle U_1, V_1, U_2, V_2 \rangle$ of~\cite{moreno_aydin}.\ In the coordinates $(x,y,z,p_x,p_y,p_z)$, let
$$
	J_0(x,y,z,p_x,p_y,p_z) = (-y,\, x,\, 0,\, p_y,\, -p_x,\, 0), 
	\qquad 
	K_0(x,y,z,p_x,p_y,p_z) = (-p_y,\, p_x,\, 0,\, -y,\, x,\, 0).
$$
Then
$$
	U_1 = \frac{J_0 \nabla H}{\omega(J_0\nabla H, K_0 \nabla H)}, 
	\qquad
	V_1 = K_0 \nabla H, 
	\qquad
	U_2 = \frac{P\partial_{p_z}}{\omega(P\partial_{p_z}, P\partial_z)}, 
	\qquad
	V_2 = P\partial_z.
$$
Here $P$ denotes the symplectic projection onto the $\omega$-orthogonal complement of $\langle Z, X_H, U_1, V_1\rangle$, where $Z = \nabla H/|\nabla H|^2$.\ The index does not depend on the choice of frames provided the frame extends to a capping disk of the orbit.\ This is the case for the frame chosen here, which is well defined unless the Hamiltonian vector field is purely vertical, i.e., $X_H \in \langle \partial_z, \partial_{p_z}\rangle$.

The index is computed as follows.\ We integrate the orbit together with its variational equations and express the linearized flow in the above frame to obtain $\Psi$.\ We then extend $\Psi(T)$ to $B_\pm$ without crossing $V$ as in~\cite{moreno_aydin}.\ Finally, we compute the degree by summing the increments of $\arg \det_{\mathbb{C}}(\rho\circ\tilde{\Psi})$.\ For orbits passing near a collision, we integrate the variational equations of~\eqref{eqn:ode-moser} in Moser regularized coordinates, then compose with the derivatives of the transformation from physical to regularized coordinates given by~\eqref{change_pos_mom} and~\eqref{eqn:moser-transform} to obtain $\Psi(t)$.
For planar orbits, $\nabla H$ has vanishing $z$ and $p_z$ components, so that $U_2 = \partial_{p_z}$ and $V_2 = \partial_z$.\ By virtue of the symmetry $\sigma$, the matrices $\Psi(t)$ split into planar and spatial blocks $\Psi = \Psi^p \oplus \Psi^s$, and the index splits accordingly: $\mu_{CZ} = \mu_{CZ}^p + \mu_{CZ}^s$.

\subsection*{Local Floer homology and Euler characteristic}

Locally near a family of non-degenerate periodic orbits, the Conley--Zehnder index provides the grading of a topological bifurcation invariant that is based on holomorphic curves~\cite{ginzburg}, the \textit{local Floer homology}.\ A useful related invariant is its \textit{Euler characteristic}, the alternating sum of the ranks of the homology groups.\ As the periodic orbits are the generators of the local homology groups graded by the Conley--Zehnder index, the index provides information on how different families are related to each other at branch points.\ To define good and bad orbits, let $\gamma$ be a single-turn periodic orbit and let $\gamma^k$ be its $k$-fold cover.\ In the sense of SFT~\cite{Eliashberg2010}, $\gamma^k$ is referred to as a \textit{bad orbit} if the underlying single-turn orbit $\gamma$ has exactly one pair of negative hyperbolic Floquet multipliers and $k$ is even.\ Otherwise, $\gamma^k$ is referred to as a \textit{good orbit}.\ Bad orbits are ignored in the local Floer homology, and only good orbits count.\ The invariance of the Euler characteristic at a bifurcation then reads
$$ \sum_{\gamma \text{: good}}^{\text{before}} (-1)^{\mu_{CZ}(\gamma)} = \sum_{\gamma \text{: good}}^{\text{after}} (-1)^{\mu_{CZ}(\gamma)}, $$
where the sums run over the periodic orbits near the bifurcating orbit before and after the bifurcation~\cite{moreno_aydin}.\

\subsection*{Krein signature and index change}

Furthermore, the \textit{Krein signature} of elliptic Floquet multipliers is computed to verify index changes as a family goes through a bifurcation, as proposed in~\cite{aydin_cz,aydin_batkhin,aydin_dro,aydin_comet}.\ The Krein signature of elliptic eigenvalues of a symplectic matrix is defined as follows~\cite[Appendix~29]{arnold_avez}: Let $M \in Sp(2n)$, and let $\lambda = e^{i \theta}$ and $\overline{\lambda} = e^{-i \theta}$ be a pair of elliptic eigenvalues of $M$ with $\lambda \neq \pm 1$.\ On the real 2-dimensional plane $P := \langle \text{Re}(v), \text{Im}(v) \rangle$, where $v$ is an eigenvector of~$M$ with eigenvalue $\lambda$, a quadratic form is defined by
\begin{align}\label{quadratic_form}
	Q(\xi,\xi) := \omega \left( \xi, M \xi \right),\quad \xi \in P.
\end{align}
This quadratic form is non-degenerate, and either positive or negative definite.\ The \textit{Krein signature} associated to the elliptic pair $\{\lambda,\overline{\lambda}\}$ is defined as the sign of $Q$, which is a symplectic invariant and independent of the choice of the eigenvector.\ If the Krein signature is positive, the rotation angle is $\vartheta := \theta \in (0,\pi)$, and if the Krein signature is negative, the rotation angle is $\vartheta := 2\pi - \theta \in (\pi,2\pi)$.\ The main result by Krein was that two pairs of elliptic Floquet multipliers with different signatures can collide and move off the unit circle and transition to complex instability, while a collision of pairs with the same signature cannot~\cite[Theorem A 29.9]{arnold_avez}.\ 

In this paper, the rotation angle is used to determine the value of the integer $d$ in the self-resonance condition~\eqref{sv}, and to verify the index change at bifurcation points.\ Namely, if the elliptic pair becomes positive hyperbolic as $\vartheta \to 0$ (resp. $\vartheta \to 2\pi$), the index decreases (resp. increases) by $1$.\ Conversely, if a positive hyperbolic pair becomes elliptic with $\vartheta$ emerging from $0$ (resp. $2\pi$), the index increases (resp. decreases) by $1$.\ If the pair goes through $1$ and remains elliptic, the index increases by $2$ if $\vartheta$ increases and decreases by $2$ if $\vartheta$ decreases.\ All bifurcation results and index changes in Section~\ref{sec:5} are consistent with the behavior of the Krein signatures.\

\section{Tables of Numerical Data}
\label{appendix_b}

In this appendix we collect tables of data containing initial conditions, stability indices, and Conley--Zehnder indices of selected periodic orbits from the families explored in this work.\ The column $T$ gives the full period of each orbit in the unregularized system, whereas the numerical correction was performed by integrating over a half or quarter period ($T/2$ or $T/4$), depending on the orbit's symmetry.\ In the elliptic region, the sign $\pm$ next to the stability index corresponds to the Krein signature.\ For degenerate orbits, the Conley--Zehnder indices on either side of the degeneracy are indicated by ``$a\shortto b$''.\

\begin{table}[!h]
\setlength{\tabcolsep}{3pt}
\centering
\small
\caption{Numerical data of vertical self-resonant orbits of the $L_1$ and $L_2$ Lyapunov families, listed in the first and second blocks, respectively.\ The last row of each block corresponds to the endpoint of the respective continuation.\ In the first block, the endpoint is a period-doubling bifurcation of a retrograde orbit around the Earth, and the index change in the final column refers to its double cover.\ The final column gives the total Conley--Zehnder index $\mu_{CZ} = \mu_{CZ}^p + \mu_{CZ}^s$, together with its planar and spatial components.\ }
	\label{tbl:B1}
\begin{tabular}{ccccccc}
\hline
$C$ & $x$ & $p_y$ & $T$ & Stability index & $\mu_{CZ}$ / $\mu_{CZ}^p$ / $\mu_{CZ}^s$ \\
\hline
$ 3.17435195$ & $0.82339090$ & $ 0.94971730$ & 2.742994 & $s_p\!=\! 1181,\, s_v\!\approx\! 1.000 $      & $3\shortto 4$ / $2$ / $1\shortto 2$ \\
$ 3.02139213$ & $0.78157393$ & $ 1.22477173$ & 3.949999 & $s_p\!=\!200.3,\, s_v\!\approx\! 1.000 $      & $4\shortto 5$ / $2$ / $2\shortto 3$ \\
$ 2.97098397$ & $0.74200952$ & $ 1.28806600$ & 5.019915 & $s_p\!=\!90.02,\, s_v\!\approx\! 0.000\, (+)$ & $5$ / $2$ / $3$ \\
$ 2.95928890$ & $0.72729740$ & $ 1.30591388$ & 5.339252 & $s_p\!=\!75.59,\, s_v\!\approx\!-0.500\, (+)$ & $5$ / $2$ / $3$ \\
$ 2.94927519$ & $0.71282811$ & $ 1.32306715$ & 5.618251 & $s_p\!=\!66.46,\, s_v\!\approx\!-1.000 $      & $5$ / $2$ / $3$ \\
$-0.46396130$ & $0.61064403$ & $-1.40766805$ & 5.434121 & $s_p\!=\!2.422,\, s_v\!\approx\!-1.000 $      & $5$ / $2$ / $3$ \\
$-0.49914210$ & $0.64259230$ & $-1.35740052$ & 5.460252 & $s_p\!=\!1.926,\, s_v\!\approx\!-0.500\, (+)$ & $5$ / $2$ / $3$ \\
$-0.54391851$ & $0.68883417$ & $-1.29059595$ & 5.497832 & $s_p\!=\!1.451,\, s_v\!\approx\! 0.000\, (+)$ & $5$ / $2$ / $3$ \\
$-0.61230376$ & $0.82226979$ & $-1.12784447$ & 5.564367 & $s_p\!\approx\!1.000,\, s_v\!=\!0.516\, (+)$  & $6\shortto 5$ / $3\shortto 2$ / $3$  \\
\hline
$3.15211890$ & $1.18089856$ & $1.02504231$ & 3.415531 & $s_p\!=\!606.1,\, s_v\!\approx\! 1.000 $      & $3\shortto 4$ / $2$ / $1\shortto 2$ \\
$3.01376751$ & $1.21997758$ & $0.79248347$ & 4.310509 & $s_p\!=\!127.9,\, s_v\!\approx\! 1.000 $      & $4\shortto 5$ / $2$ / $2\shortto 3$ \\
$2.97326315$ & $1.25096005$ & $0.75281111$ & 5.165816 & $s_p\!=\!73.95,\, s_v\!\approx\! 0.000\, (+)$ & $5$ / $2$ / $3$ \\
$2.96368427$ & $1.26295814$ & $0.74416382$ & 5.451914 & $s_p\!=\!65.39,\, s_v\!\approx\!-0.500\, (+)$ & $5$ / $2$ / $3$ \\
$2.95572007$ & $1.27473051$ & $0.73696910$ & 5.713490 & $s_p\!=\!59.57,\, s_v\!\approx\!-1.000 $      & $5$ / $2$ / $3$ \\
$2.83529660$ & $1.54428542$ & $0.61554036$ & 8.770529 & $s_p\!=\!107.0,\, s_v\!\approx\!-1.000 $      & $5$ / $2$ / $3$ \\
$2.80254910$ & $1.61072962$ & $0.58685572$ & 9.100501 & $s_p\!=\!145.0,\, s_v\!\approx\! 1.000 $      & $5\shortto 6$ / $2$ / $3\shortto 4$ \\
$0.80895979$ & $2.96147346$ & $0.02241137$ & 12.390909 & $s_p\!=\!5061,\, s_v\!=\!376.3$              & $6$ / $2$ / $4$ \\
\hline
\end{tabular}
\end{table}

\begin{table}[!h]
\setlength{\tabcolsep}{3pt}
\centering
\small
\caption{Numerical data of the family emerging from the period-doubling bifurcations of the $L_1$ Lyapunov family.}
	\label{tbl:B2}
\begin{tabular}{ccccccc}
\hline
$C$ & $x$ & $p_y$ & $p_z$ & $T$ & Stability index & $\mu_{CZ}$ \\
\hline
$ 2.94927519$ & $0.71282811$ & $ 1.32306715$ & $0.$         & 11.236502 & $s_1\!=\!8832,\, s_2\!\approx\!1.000$ & $9\shortto 10$ \\
$ 2.94883824$ & $0.71299767$ & $ 1.32243458$ & $0.03277716$ & 11.236230 & $s_1\!=\!8918,\, s_2\!=\!0.822\, (-)$ & $9$ \\
$ 2.93133031$ & $0.71985293$ & $ 1.29715392$ & $0.20805731$ & 11.227156 & $s_1\!=\!1.3\cdot 10^{4},\, s_2\!=\!0.025\, (+)$ & $9$ \\
$ 2.92807004$ & $0.72114058$ & $ 1.29246160$ & $0.22576683$ & 11.225996 & $s_1\!=\!1.4\cdot 10^{4},\, s_2\!\approx\!1.000$ & $9\shortto 8$ \\
$ 2.32891577$ & $0.82880938$ & $ 0.72917622$ & $0.92207326$ & 12.173813 & $s_1\!=\!9.3\cdot 10^{6},\, s_2\!=\!3.3\cdot 10^{4}$ & $8$ \\
$ 0.45426234$ & $0.87746664$ & $-0.39372574$ & $1.06805919$ & 12.531276 & $s_1\!=\!1.1\cdot 10^{7},\, s_2\!=\!2.5\cdot 10^{5}$ & $8$ \\
$-0.75957225$ & $0.84901283$ & $-1.11729737$ & $0.28857260$ & 12.214579 & $s_1\!=\!6.662,\, s_2\!\approx\!1.000$ & $8\shortto 9$ \\
$-0.75951055$ & $0.84749423$ & $-1.11923188$ & $0.28488481$ & 12.203701 & $s_1\!=\!5.976,\, s_2\!=\!0.119\, (+)$ & $9$ \\
$-0.75409690$ & $0.83257873$ & $-1.13634774$ & $0.25685710$ & 12.101499 & $s_1\!=\!2.176,\, s_2\!=\!0.326\, (-)$ & $9$ \\
$-0.75306427$ & $0.83084986$ & $-1.13819392$ & $0.25431929$ & 12.090126 & $s_1\!=\!1.591,\, s_2\!\approx\!1.000$ & $9\shortto 10$ \\
$-0.71792307$ & $0.79091367$ & $-1.17900200$ & $0.21563021$ & 11.846163 & $s_1,s_2\!=\!2.059\pm2.449i$ & $10$ \\
$-0.46594016$ & $0.61180413$ & $-1.40578909$ & $0.02447175$ & 10.874914 & $s_1\!=\!9.182,\, s_2\!=\!1.289$ & $10$ \\
$-0.46396130$ & $0.61064403$ & $-1.40766805$ & $0.$         & 10.868241 & $s_1\!=\!10.73,\, s_2\!\approx\!1.000$ & $10\shortto 9$ \\
\hline
\end{tabular}
\end{table}

\begin{table}[!h]
\setlength{\tabcolsep}{3pt}
\centering
\small
\caption{Numerical data of the two families emerging from the period-tripling bifurcations of the $L_1$ Lyapunov family.}
	\label{tbl:B3}
\begin{tabular}{ccccccc}
\hline
$C$ & $x$ & $p_y$ & $p_z$ & $T$ & Stability index & $\mu_{CZ}$ \\
\hline
$2.95928993$ & $0.72729880$ & $1.30591221$ & $0.$ & 16.017669 & $s_1\!=\!1.7\cdot 10^{6},\, s_2\!\approx\!1.000$ & $13\shortto 15$ \\
$2.88971670$ & $0.77544400$ & $1.16790688$ & $0.42485477$ & 19.024564 & $s_1\!=\!1.7\cdot 10^{6},\, s_2\!=\!212.1$ & $14$ \\
$1.21599329$ & $0.98247724$ & $-1.51930337$ & $0.01600346$ & 18.462203 & $s_1\!=\!4.1\cdot 10^{8},\, s_2\!=\!5.5\cdot 10^{6}$ & $14$ \\
$-0.72396486$ & $0.95092702$ & $-1.12949420$ & $0.09841999$ & 17.868097 & $s_1\!=\!3.658,\, s_2\!\approx\!1.000$ & $14\shortto 15$ \\
$-0.72392358$ & $0.95094221$ & $-1.12950112$ & $0.09909963$ & 17.864369 & $s_1\!=\!3.337,\, s_2\!=\!0.272\, (+)$ & $15$ \\
$-0.72187456$ & $0.95077082$ & $-1.12830087$ & $0.10233603$ & 17.830248 & $s_1\!=\!1.903,\, s_2\!=\!-3.447$ & $15$ \\
$-0.71613573$ & $0.95006315$ & $-1.12462556$ & $0.10406320$ & 17.770965 & $s_1\!\approx\!1.000,\, s_2\!=\!-6.430$ & $15\shortto 14$ \\
$-0.61370768$ & $0.93824048$ & $-1.07542430$ & $0.08853823$ & 17.062641 & $s_1\!=\!0.369\,(-),\, s_2\!=\!-15.76$ & $14$ \\
$-0.51256306$ & $0.92728041$ & $-1.04130584$ & $0.03708307$ & 16.461985 & $s_1\!=\!10.31,\, s_2\!=\!1.006$ & $14$ \\
$-0.49914210$ & $0.92564498$ & $-1.03735651$ & $0.$ & 16.380755 & $s_1\!=\!22.78,\, s_2\!\approx\!1.000$ & $15\shortto 13$ \\
\hline
$C$ & $x$ & $z$ & $p_y$ & $T$ & Stability index & $\mu_{CZ}$ \\
\hline
$2.95928932$ & $0.95666838$ & $0.$         & $0.07657270$ & 16.017721 & $s_1\!=\!1.7\cdot 10^{6},\, s_2\!\approx\!1.000$ & $13\shortto 15$ \\
$2.95210007$ & $0.95951773$ & $0.00865723$ & $0.05282265$ & 16.225781 & $s_1\!=\!1.8\cdot 10^{6},\, s_2\!=\!0.568\, (-)$ & $13$ \\
$2.30039507$ & $0.98755176$ & $0.00245258$ & $-2.25087284$ & 19.318824 & $s_1\!=\!2.4\cdot 10^{9},\, s_2\!=\!-1.3\cdot 10^{6}$ & $13$ \\
$0.39757053$ & $0.98760523$ & $0.00364819$ & $-2.04558619$ & 18.059885 & $s_1\!=\!6.1\cdot 10^{8},\, s_2\!=\!-1.0\cdot 10^{7}$ & $13$ \\
$-0.81339207$ & $0.99602941$ & $0.01592630$ & $-1.26719444$ & 18.598472 & $s_1\!=\!123.1,\, s_2\!=\!-0.118\, (+)$ & $13$ \\
$-0.81546761$ & $0.99601993$ & $0.01588274$ & $-1.26838540$ & 18.613816 & $s_1\!=\!129.0,\, s_2\!\approx\!1.000$ & $13\shortto 12$ \\
$-0.85335174$ & $0.99510301$ & $0.01660967$ & $-1.27328429$ & 18.975595 & $s_1\!=\!233.9,\, s_2\!=\!31.97$ & $12$ \\
$-0.86293060$ & $0.99325437$ & $0.01740567$ & $-1.27558340$ & 19.277960 & $s_1\!=\!220.5,\, s_2\!\approx\!1.000$ & $12\shortto 11$ \\
$-0.86291841$ & $0.99317373$ & $0.01740597$ & $-1.27603953$ & 19.287845 & $s_1\!=\!209.3,\, s_2\!=\!-0.612\, (-)$ & $11$ \\
$-0.86284357$ & $0.99303762$ & $0.01740098$ & $-1.27687310$ & 19.304163 & $s_1\!=\!183.5,\, s_2\!=\!0.107\, (+)$ & $11$ \\
$-0.86282319$ & $0.99301335$ & $0.01739937$ & $-1.27703017$ & 19.307028 & $s_1\!=\!177.7,\, s_2\!\approx\!1.000$ & $11\shortto 10$ \\
$-0.85770975$ & $0.99154248$ & $0.01690891$ & $-1.29167856$ & 19.466071 & $s_1,s_2\!=\!-406.4\pm704.7i$ & $10$ \\
\hline
\end{tabular}
\end{table}

\begin{table}[!h]
\setlength{\tabcolsep}{3pt}
\centering
\small
\caption{Numerical data of the two families emerging from the period-quadrupling bifurcations of the $L_1$ Lyapunov family.}
	\label{tbl:B4}
\begin{tabular}{ccccccc}
\hline
$C$ & $x$ & $z$ & $p_y$ & $T$ & Stability index & $\mu_{CZ}$ \\
\hline
$2.97098310$ & $0.95228647$ & $0.$         & $0.13512913$ & 20.079752 & $s_1\!=\!5.3\cdot 10^{8},\, s_2\!\approx\!1.000$ & $17\shortto 19$ \\
$2.91593414$ & $0.66238922$ & $0.03851421$ & $1.38431062$ & 23.233942 & $s_1\!=\!3.4\cdot 10^{8},\, s_2\!=\!425.5$ & $18$ \\
$1.83985256$ & $0.01660373$ & $0.00039932$ & $8.19514114$ & 25.709072 & $s_1\!=\!1.8\cdot 10^{12},\, s_2\!=\!4.5\cdot 10^{8}$ & $18$ \\
$0.53035296$ & $0.10729169$ & $0.00139714$ & $-3.89878601$ & 23.742565 & $s_1\!=\!8.7\cdot 10^{10},\, s_2\!=\!6.6\cdot 10^{8}$ & $18$ \\
$-0.70212168$ & $0.75494902$ & $0.02882468$ & $-1.23234784$ & 23.835822 & $s_1\!=\!143.5,\, s_2\!\approx\!1.000$ & $18\shortto 17$ \\
$-0.70455841$ & $0.75743019$ & $0.02972873$ & $-1.22957357$ & 23.819104 & $s_1\!=\!66.04,\, s_2\!=\!-0.471\, (-)$ & $17$ \\
$-0.70739968$ & $0.76144680$ & $0.03223349$ & $-1.22479646$ & 23.750794 & $s_1\!\approx\!1.000,\, s_2\!=\!-0.960\, (+)$ & $17\shortto 18$ \\
$-0.62808752$ & $0.72209138$ & $0.04003913$ & $-1.25950719$ & 22.799628 & $s_1,s_2\!=\!-3.887\pm0.4136i$ & $18$ \\
$-0.55177356$ & $0.69122609$ & $0.01717726$ & $-1.28852041$ & 22.070577 & $s_1\!=\!14.74,\, s_2\!=\!1.058$ & $18$ \\
$-0.54391848$ & $0.68883414$ & $0.$         & $-1.29059600$ & 21.991327 & $s_1\!=\!19.60,\, s_2\!\approx\!1.000$ & $19\shortto 17$ \\
\hline
$2.97098310$ & $0.95228647$ & $0.$         & $0.13512913$ & 20.079752 & $s_1\!=\!5.3\cdot 10^{8},\, s_2\!\approx\!1.000$ & $17\shortto 19$ \\
$2.96762321$ & $0.95387241$ & $0.00729687$ & $0.12522343$ & 20.121111 & $s_1\!=\!5.6\cdot 10^{8},\, s_2\!=\!0.629\, (-)$ & $17$ \\
$2.95558559$ & $0.95927420$ & $0.01388705$ & $0.08570923$ & 20.272678 & $s_1\!=\!7.1\cdot 10^{8},\, s_2\!=\!0.173\, (+)$ & $17$ \\
$2.95449555$ & $0.95973981$ & $0.01420885$ & $0.08180137$ & 20.286648 & $s_1\!=\!7.3\cdot 10^{8},\, s_2\!\approx\!1.000$ & $17\shortto 16$ \\
$0.24164861$ & $0.98767088$ & $0.00390993$ & $-1.99874039$ & 25.114058 & $s_1\!=\!3.5\cdot 10^{10},\, s_2\!=\!1.2\cdot 10^{8}$ & $16$ \\
$-0.86265429$ & $0.99318820$ & $0.01787859$ & $-1.26875320$ & 26.048410 & $s_1\!=\!1.2\cdot 10^{6},\, s_2\!\approx\!1.000$ & $16\shortto 15$ \\
$-0.86269167$ & $0.99317106$ & $0.01787706$ & $-1.26887669$ & 26.053378 & $s_1\!=\!1.2\cdot 10^{6},\, s_2\!=\!-0.432\,(-)$ & $15$ \\
$-0.86276736$ & $0.99313292$ & $0.01787321$ & $-1.26915585$ & 26.064379 & $s_1\!=\!1.2\cdot 10^{6},\, s_2\!=\!-0.650\,(+)$ & $15$ \\
$-0.86281408$ & $0.99310641$ & $0.01787018$ & $-1.26935349$ & 26.071986 & $s_1\!=\!1.3\cdot 10^{6},\, s_2\!\approx\!1.000$ & $15\shortto 14$ \\
$-0.86299494$ & $0.99295410$ & $0.01784729$ & $-1.27054524$ & 26.115157 & $s_1\!=\!1.5\cdot 10^{6},\, s_2\!=\!13.51$ & $14$ \\
$-0.86304118$ & $0.99282430$ & $0.01782067$ & $-1.27163645$ & 26.151403 & $s_1\!=\!1.7\cdot 10^{6},\, s_2\!\approx\!1.000$ & $14\shortto 13$ \\
$-0.86303760$ & $0.99278752$ & $0.01781197$ & $-1.27195837$ & 26.161609 & $s_1\!=\!1.8\cdot 10^{6},\, s_2\!=\!-0.920\,(+)$ & $13$ \\
$-0.86303276$ & $0.99276777$ & $0.01780710$ & $-1.27213348$ & 26.167078 & $s_1\!=\!1.8\cdot 10^{6},\, s_2\!\approx\!1.000$ & $13\shortto 12$ \\
$-0.75561606$ & $0.98832833$ & $0.01110322$ & $-1.43935914$ & 28.393138 & $s_1\!=\!9.4\cdot 10^{12},\, s_2\!=\!1.6\cdot 10^{12}$ & $12$ \\
\hline
\end{tabular}
\end{table}

\begin{table}[t!]
\setlength{\tabcolsep}{3pt}
\centering
\small
\caption{Numerical data of the family emerging from the period-doubling bifurcations of the $L_2$ Lyapunov family.}
	\label{tbl:B5}
\begin{tabular}{ccccccc}
\hline
$C$ & $x$ & $p_y$ & $p_z$ & $T$ & Stability index & $\mu_{CZ}$ \\
\hline
$2.95572007$ & $1.27473051$ & $0.73696910$ & $0.$         & 11.426980 & $s_1\!=\!7096,\, s_2\!\approx\!1.000$ & $9\shortto 10$ \\
$2.95347280$ & $1.27354431$ & $0.73670328$ & $0.04451647$ & 11.411586 & $s_1\!=\!7646,\, s_2\!=\!-0.012\,(-)$ & $9$ \\
$2.94448410$ & $1.26872931$ & $0.73568412$ & $0.10012592$ & 11.352114 & $s_1\!=\!1.0\cdot 10^{4},\, s_2\!=\!-0.632\,(+)$ & $9$ \\
$2.93970471$ & $1.26612783$ & $0.73517137$ & $0.11991892$ & 11.322192 & $s_1\!=\!1.2\cdot 10^{4},\, s_2\!\approx\!1.000$ & $9\shortto 8$ \\
$2.82239314$ & $1.20731220$ & $0.72339619$ & $0.36318569$ & 11.194610 & $s_1\!=\!1.1\cdot 10^{5},\, s_2\!=\!770.8$ & $8$ \\
$2.10838148$ & $1.12754538$ & $0.46769598$ & $0.79690960$ & 12.273667 & $s_1\!=\!2.6\cdot 10^{7},\, s_2\!=\!1.7\cdot 10^{5}$ & $8$ \\
$-1.01998825$ & $1.03594071$ & $-1.07635681$ & $0.14737180$ & 11.917080 & $s_1\!=\!53.17,\, s_2\!\approx\!1.000$ & $8\shortto 9$ \\
$-1.02316693$ & $1.03340143$ & $-1.08704778$ & $0.13373539$ & 11.870048 & $s_1\!=\!61.32,\, s_2\!=\!-0.536\,(+)$ & $9$ \\
$-1.03930262$ & $1.02416750$ & $-1.13513827$ & $0.03464973$ & 11.668335 & $s_1\!=\!51.31,\, s_2\!=\!0.931\,(-)$ & $9$ \\
$-1.04059623$ & $1.02363253$ & $-1.13853958$ & $0.$         & 11.654822 & $s_1\!=\!44.96,\, s_2\!\approx\!1.000$ & $11\shortto 9$ \\
\hline
\end{tabular}
\end{table}

\begin{table}[!h]
\setlength{\tabcolsep}{3pt}
\centering
\small
\caption{Numerical data of the two families emerging from the period-tripling bifurcations of the $L_2$ Lyapunov family.}
	\label{tbl:B6}
\begin{tabular}{ccccccc}
\hline
$C$ & $x$ & $p_y$ & $p_z$ & $T$ & Stability index & $\mu_{CZ}$ \\
\hline
$2.96368427$ & $1.26295814$ & $0.74416382$ & $0.$         & 16.355742 & $s_1\!=\!1.1\cdot 10^{6},\, s_2\!\approx\!1.000$ & $13\shortto 15$ \\
$2.93427990$ & $1.23074459$ & $0.75088127$ & $0.19956354$ & 18.070808 & $s_1\!=\!1.0\cdot 10^{6},\, s_2\!=\!49.95$ & $14$ \\
$2.78625505$ & $1.15096027$ & $0.74947459$ & $0.47421794$ & 23.175878 & $s_1\!=\!3.8\cdot 10^{6},\, s_2\!=\!821.7$ & $14$ \\
$-1.11942374$ & $1.11189456$ & $-0.96231317$ & $0.08348652$ & 25.484470 & $s_1\!=\!211.7,\, s_2\!\approx\!1.000$ & $14\shortto 13$ \\
$-1.11954678$ & $1.11197544$ & $-0.96225398$ & $0.08323936$ & 25.484755 & $s_1\!=\!210.4,\, s_2\!=\!0.092\,(-)$ & $13$ \\
$-1.14431940$ & $1.13089803$ & $-0.94768323$ & $0.03281975$ & 25.549155 & $s_1\!=\!18.30,\, s_2\!=\!-107.3$ & $13$ \\
$-1.15170260$ & $1.13722598$ & $-0.94282218$ & $0.$         & 25.569706 & $s_1\!\approx\!1.000,\, s_2\!=\!-137.2$ & $14\shortto 13$ \\
\hline
$C$ & $x$ & $z$ & $p_y$ & $T$ & Stability index & $\mu_{CZ}$ \\
\hline
$2.96368427$ & $1.26295814$ & $0.$         & $0.74416382$ & 16.355742 & $s_1\!=\!1.1\cdot 10^{6},\, s_2\!\approx\!1.000$ & $13\shortto 15$ \\
$2.95535991$ & $1.27355275$ & $0.02765679$ & $0.73731526$ & 16.581767 & $s_1\!=\!1.2\cdot 10^{6},\, s_2\!=\!-0.328\,(-)$ & $13$ \\
$2.93930460$ & $1.29887458$ & $0.04943956$ & $0.72367755$ & 17.079103 & $s_1\!=\!1.5\cdot 10^{6},\, s_2\!=\!-30.40$ & $13$ \\
$2.36178671$ & $1.04001642$ & $0.96591268$ & $0.68193685$ & 21.405627 & $s_1\!=\!1.4\cdot 10^{8},\, s_2\!=\!-7.9\cdot 10^{4}$ & $13$ \\
$-1.08924087$ & $-1.38311942$ & $0.23373209$ & $0.72236379$ & 22.072496 & $s_1\!\approx\!1.000,\, s_2\!=\!-6.983$ & $13\shortto 12$ \\
$-1.10336572$ & $-1.34409955$ & $0.15614114$ & $0.75199038$ & 22.158735 & $s_1\!=\!-0.129\,(-),\, s_2\!=\!-4.274$ & $12$ \\
$-1.11884999$ & $-1.29868052$ & $0.02116219$ & $0.78734951$ & 22.151254 & $s_1\!=\!0.911\,(-),\, s_2\!=\!0.448\,(-)$ & $12$ \\
$-1.11919077$ & $-1.29762551$ & $0.$         & $0.78818282$ & 22.150501 & $s_1\!\approx\!1.000,\, s_2\!=\!0.463\,(-)$ & $13\shortto 12$ \\
\hline
\end{tabular}
\end{table}

\begin{table}[!h]
\setlength{\tabcolsep}{3pt}
\centering
\small
\caption{Numerical data of the two families emerging from the period-quadrupling bifurcations of the $L_2$ Lyapunov family.}
	\label{tbl:B7}
\begin{tabular}{ccccccc}
\hline
$C$ & $x$ & $z$ & $p_y$ & $T$ & Stability index & $\mu_{CZ}$ \\
\hline
$2.97326221$ & $1.25096112$ & $0.$         & $0.75281026$ & 20.663370 & $s_1\!=\!2.4\cdot 10^{8},\, s_2\!\approx\!1.000$ & $17\shortto 19$ \\
$2.92474163$ & $1.32168260$ & $0.08004396$ & $0.71221552$ & 23.844574 & $s_1\!=\!1.9\cdot 10^{8},\, s_2\!=\!622.2$ & $18$ \\
$-0.60849952$ & $-1.10344276$ & $0.92078411$ & $0.69110626$ & 31.559419 & $s_1\!=\!6.5\cdot 10^{6},\, s_2\!=\!5.7\cdot 10^{4}$ & $18$ \\
$-1.12027838$ & $-1.41058308$ & $0.19079080$ & $0.71543446$ & 31.847082 & $s_1\!=\!203.1,\, s_2\!\approx\!1.000$ & $18\shortto 17$ \\
$-1.12429356$ & $-1.41079239$ & $0.17429667$ & $0.71675896$ & 31.863331 & $s_1\!=\!149.1,\, s_2\!=\!-0.716\,(-)$ & $17$ \\
$-1.13618277$ & $-1.41001190$ & $0.11417013$ & $0.72144799$ & 31.915351 & $s_1\!=\!37.19,\, s_2\!=\!0.281\,(+)$ & $17$ \\
$-1.13937913$ & $-1.40940599$ & $0.09225235$ & $0.72292784$ & 31.929940 & $s_1\!=\!19.22,\, s_2\!\approx\!1.000$ & $17\shortto 16$ \\
$-1.14521304$ & $-1.40781485$ & $0.02041560$ & $0.72590144$ & 31.956734 & $s_1\!=\!3.120,\, s_2\!=\!1.403$ & $16$ \\
$-1.14551777$ & $-1.40771351$ & $0.$         & $0.72606708$ & 31.958132 & $s_1\!=\!3.326,\, s_2\!\approx\!1.000$ & $17\shortto 16$ \\
\hline
$2.97326221$ & $1.25096112$ & $0.$         & $0.75281026$ & 20.663370 & $s_1\!=\!2.4\cdot 10^{8},\, s_2\!\approx\!1.000$ & $17\shortto 19$ \\
$2.97046511$ & $1.00945141$ & $0.00440323$ & $2.05094817$ & 20.663619 & $s_1\!=\!2.6\cdot 10^{8},\, s_2\!=\!0.616\,(-)$ & $17$ \\
$2.96278655$ & $1.00595178$ & $0.00741309$ & $2.11589373$ & 20.665443 & $s_1\!=\!3.3\cdot 10^{8},\, s_2\!=\!-0.827\,(+)$ & $17$ \\
$2.96003526$ & $1.00478086$ & $0.00789721$ & $2.14189684$ & 20.666666 & $s_1\!=\!3.6\cdot 10^{8},\, s_2\!\approx\!1.000$ & $17\shortto 16$ \\
$2.87216971$ & $0.98847586$ & $0.00118778$ & $5.25180440$ & 21.228504 & $s_1\!=\!2.0\cdot 10^{9},\, s_2\!=\!1.7\cdot 10^{4}$ & $16$ \\
$1.52298360$ & $0.98799116$ & $0.00194002$ & $-2.74316409$ & 24.901548 & $s_1\!=\!4.4\cdot 10^{12},\, s_2\!=\!3.4\cdot 10^{8}$ & $16$ \\
$-1.04203881$ & $1.01531244$ & $0.03026040$ & $-1.12703013$ & 24.332414 & $s_1\!=\!143.2,\, s_2\!\approx\!1.000$ & $16\shortto 17$ \\
$-1.07167418$ & $1.03440079$ & $0.01188150$ & $-1.09519293$ & 24.203598 & $s_1\!=\!88.39,\, s_2\!=\!0.994\,(-)$ & $17$ \\
$-1.07631193$ & $1.03708327$ & $0.$         & $-1.09092804$ & 24.196991 & $s_1\!=\!70.58,\, s_2\!\approx\!1.000$ & $19\shortto 17$ \\
\hline
\end{tabular}
\end{table}

\begin{table}[!h]
\setlength{\tabcolsep}{3pt}
\centering
\small
\caption{Numerical data of the family emerging from the near-collision single-turn and period-doubling bifurcations of the $L_2$ Lyapunov family (first and second block, correspondingly).}
	\label{tbl:B8}
\begin{tabular}{ccccccc}
\hline
$C$ & $x$ & $z$ & $p_y$ & $T$ & Stability index & $\mu_{CZ}$ \\
\hline
$2.80254910$ & $1.61072962$ & $0.$         & $0.58685572$ & 9.100501 & $s_1\!=\!145.0,\, s_2\!\approx\! 1.000 $ & $5\shortto 6$ \\
$2.79037943$ & $1.59584464$ & $0.18991333$ & $0.58852964$ & 9.107888 & $s_1\!=\!151.3,\, s_2\!=\!-0.554\,(-)$ & $5$ \\
$2.55474032$ & $1.34430010$ & $0.79412841$ & $0.61130969$ & 9.235817 & $s_1\!=\!328.9,\, s_2\!=\!-27.13$ & $5$ \\
$-0.82707915$ & $-1.31846309$ & $0.73484152$ & $0.65476728$ & 9.479441 & $s_1\!=\!155.6,\, s_2\!=\!-72.64$ & $5$ \\
$-1.13902689$ & $-1.46207265$ & $0.00627665$ & $0.69413092$ & 9.482272 & $s_1\!=\!20.37,\, s_2\!=\!0.994\,(-)$ & $5$ \\
$-1.13905411$ & $-1.46206836$ & $0.$         & $0.69414252$ & 9.482258 & $s_1\!=\!20.37,\, s_2\!\approx\!1.000$ & $6\shortto 5$ \\
\hline
$2.83529544$ & $1.54428788$ & $0.$         & $0.61553929$ & 17.541085 & $s_1\!=\!2.3\cdot 10^{4},\, s_2\!\approx\!1.000$ & $10\shortto 11$ \\
$2.81235454$ & $1.51543202$ & $0.24965361$ & $0.61975863$ & 17.661490 & $s_1\!=\!2.1\cdot 10^{4},\, s_2\!=\!19.65$ & $10$ \\
$-0.60557577$ & $-1.12935926$ & $0.94967085$ & $0.67033705$ & 18.971519 & $s_1\!=\!9742,\, s_2\!=\!846.0$ & $10$ \\
$-1.15459880$ & $-1.46706367$ & $0.07365913$ & $0.69464442$ & 19.128149 & $s_1\!=\!759.4,\, s_2\!=\!5.156$ & $10$ \\
$-1.15850741$ & $-1.46834781$ & $0.$         & $0.69526000$ & 19.129681 & $s_1\!=\!737.4,\, s_2\!\approx\!1.000$ & $11\shortto 10$ \\
\hline
\end{tabular}
\end{table}

\begin{table}[!ht]
\setlength{\tabcolsep}{2pt}
\centering
\footnotesize
\caption{Numerical data for critical orbits appearing in the bifurcation networks involving NRHOs and planar retrograde and direct orbits.\ The first block lists the self-resonant NRHOs, the second block lists the planar endpoint orbits, and the third block lists the critical orbits involved in the implicit connections.}
\label{tbl:nrho-critical-orbits}
\begin{tabular}{cccccccc}
\hline
Name & $C$ & $x$ & $z$ & $p_y$ & $T$ & Stability index & $\mu_{CZ}$ \\
\hline
$1$:$2$ $L_1$ NRHO & $2.94344119$ & $0.92917310$ & $0.29137298$ & $1.01089978$ & $2.150939$ & $s_1\!=\!0.712\,(-),\, s_2\!\approx\!-1.000$ & $5\!\shortto\! 4$ \\
$1$:$2$ $L_1$ NRHO & $2.99863011$ & $0.88118523$ & $0.19370917$ & $1.10281265$ & $2.113270$ & $s_1\!=\!0.602\,(+),\, s_2\!\approx\!-1.000$ & $7\!\shortto\! 6$ \\
$1$:$2$ $L_2$ NRHO & $3.01714367$ & $1.06876977$ & $0.20136006$ & $0.88441353$ & $2.169092$ & $s_1\!=\!-0.087\,(+),\, s_2\!\approx\!-1.000$ & $7\!\shortto\! 6$ \\
$1$:$2$ $L_2$ NRHO & $3.05802218$ & $1.01182891$ & $0.17391370$ & $0.93189789$ & $1.374328$ & $s_1\!=\!0.769\,(+),\, s_2\!\approx\!-1.000$ & $8\!\shortto\! 7$ \\
$2$:$3$ $L_2$ NRHO & $3.07783646$ & $1.00032341$ & $0.15930058$ & $0.94994570$ & $1.190605$ & $s_1\!=\!0.827\,(+),\, s_2\!\approx\!-0.500\,(-)$ & $12\!\shortto\! 10$ \\
$3$:$4$ $L_2$ NRHO & $3.10429009$ & $0.99324249$ & $0.14077724$ & $0.96407022$ & $1.000283$ & $s_1\!=\!0.849\,(+),\, s_2\!\approx\!0.000\,(-)$ & $16\!\shortto\! 14$ \\
$4$:$5$ $L_2$ NRHO & $3.12850828$ & $0.99058077$ & $0.12607043$ & $0.97103146$ & $0.859534$ & $s_1\!=\!0.863\,(+),\, s_2\!\approx\!\cos(2\pi/5)\,(-)$ & $20\!\shortto\! 18$ \\
$5$:$6$ $L_2$ NRHO & $3.15052568$ & $0.98942208$ & $0.11462996$ & $0.97493289$ & $0.753105$ & $s_1\!=\!0.878\,(+),\, s_2\!\approx\!0.500\,(-)$ & $24\!\shortto\! 22$ \\
\hline
$1$:$1$ DPO   & $3.02640034$ & $1.01897928$ & $0.$ & $1.86076860$ & $5.124068$ & $s_p\!=\!884.7,\, s_v\!\approx\!1.000$ & $5\!\shortto\! 6$ \\
$1$:$1$ $LPO_1$ & $3.13629470$ & $0.87971912$ & $0.$ & $0.60118129$ & $3.364393$ & $s_p\!=\!-236.4,\, s_v\!\approx\!1.000$ & $6\!\shortto\! 7$ \\
$1$:$1$ $LPO_2$ & $3.12466614$ & $1.09261784$ & $0.$ & $1.39167705$ & $3.910194$ & $s_p\!=\!-107.0,\, s_v\!\approx\!1.000$ & $6\!\shortto\! 7$ \\
$1$:$2$ $LPO_1$ & $3.18487359$ & $0.85900709$ & $0.$ & $0.76137673$ & $2.009316$ & $s_p\!=\!-0.364\,(+),\, s_v\!\approx\!-1.000$ & $11\!\shortto\! 12$\\
$1$:$2$ $LPO_2$ & $3.17065216$ & $1.12351796$ & $0.$ & $1.22570003$ & $2.323728$ & $s_p\!=\!-0.519\,(+),\, s_v\!\approx\!-1.000$ & $11\!\shortto\! 12$\\
$1$:$3$ $LPO_2$ & $3.17529837$ & $1.10058572$ & $0.$ & $1.26514610$ & $1.698176$ & $s_p\!=\!0.393\,(+),\, s_v\!\approx\!-0.500\,(+)$ & $14\!\shortto\! 16$\\
$1$:$4$ $LPO_1$ & $3.19042228$ & $0.89415817$ & $0.$ & $0.67413150$ & $1.294848$ & $s_p\!=\!0.882\,(+),\, s_v\!\approx\!0.000\,(+)$ & $18\!\shortto\! 20$\\
$1$:$5$ $LPO_1$ & $3.19593933$ & $1.04523881$ & $0.$ & $1.47940121$ & $1.138585$ & $s_p\!=\!0.921\,(+),\, s_v\!\approx\!\cos(2\pi/5)\,(+)$ & $22\!\shortto\! 24$\\
$5$:$6$ DRO   & $3.02343980$ & $1.07820497$ & $0.$ & $0.60914918$ & $1.293528$ & $s_p\!=\!0.284\,(-),\, s_v\!\approx\!0.500\,(-)$ & $20\!\shortto\! 18$\\
$6$:$7$ DRO   & $3.04984732$ & $0.91148346$ & $0.$ & $1.39956554$ & $1.056877$ & $s_p\!=\!0.483\,(-),\, s_v\!\approx\!\cos(2\pi/7)\,(-)$ & $24\!\shortto\! 22$\\
\hline
{\tiny $2$:$3$ $L_2$ NRHO $\!\leftrightarrow\!$ $1$:$2$ $LPO_1$} & $3.05323592$ & $0.96037619$ & $-0.10979604$ & $0.64024519$ & $6.314176$ & $s_1\!=\!7.016,\, s_2\!\approx\!1.000$ & $11\!\shortto\! 12$ \\
{\tiny $2$:$3$ $L_2$ NRHO $\!\leftrightarrow\!$ $1$:$2$ $LPO_2$} & $3.07753697$ & $1.03304677$ & $0.11102138$ & $1.30150443$ & $5.896381$ & $s_1\!=\!28.44,\, s_2\!\approx\!1.000$ & $11\!\shortto\! 12$ \\
{\tiny $4$:$5$ $L_2$ NRHO $\!\leftrightarrow\!$ $1$:$4$ $LPO_1$} & $3.14409046$ & $0.94629451$ & $0.05594400$ & $0.54886773$ & $5.896286$ & $s_1\!=\!0.998\,(-),\, s_2\!\approx\!1.000$ & $18\!\shortto\! 19$ \\
{\tiny $4$:$5$ $L_2$ NRHO $\!\leftrightarrow\!$ $1$:$4$ $LPO_1$} & $3.14409790$ & $1.03268353$ & $0.05891572$ & $1.40483810$ & $5.896433$ & $s_1\!=\!1.002,\, s_2\!\approx\!1.000$ & $19\!\shortto\! 20$ \\
\hline
\end{tabular}
\end{table}

\end{appendix}

\clearpage

\noindent
\textbf{Acknowledgments.}\, 
CA acknowledges support by the Deutsche Forschungsgemeinschaft (DFG, German Research Foundation), Project-ID 541062288.\ CJ was supported by BK21 SNU Mathematical Sciences Division.\ The authors thank Professor Otto van Koert for assisting the preparation of this manuscript.\ CJ would like to thank researchers at Universität Heidelberg for their kind hospitality during his two visits.\

\bigskip
\noindent
\textbf{Data Availability.}\, 
All numerical results presented in Tables~\ref{tbl:periodic-collision-orbits}--\ref{tbl:nrho-critical-orbits} were computed using the numerical algorithms developed and implemented by the authors.\ Claude Opus 5.5 was used only at the final stages to assist in preparing the public repository and checking the reproducibility of the results.\ The source code and scripts for reproducing every table entry from the supplied initial data are publicly available at \url{https://github.com/ckjoung/moser}.\

\addcontentsline{toc}{section}{References}
\hypertarget{bibliolist}{}\printbibliography

@article{aydin_batkhin,
	author = {Aydin, C. and Batkhin, A.},
	journal = {Celest. Mech. Dyn. Astron.},
	title = {Studying network of symmetric periodic orbit families of the Hill problem via symplectic invariants},
	year = {2025},
    volume = {137},
    number = {12},
    doi = {10.1007/s10569-025-10241-7}}

@article {moser_1970_regularization,
    AUTHOR = {Moser, J.},
     TITLE = {Regularization of {K}epler's problem and the averaging method
              on a manifold},
   JOURNAL = {Comm. Pure Appl. Math.},
  FJOURNAL = {Communications on Pure and Applied Mathematics},
    VOLUME = {23},
      YEAR = {1970},
     PAGES = {609--636},
       DOI = {10.1002/cpa.3160230406},
}

@article {belbruno_2019_family,
    AUTHOR = {Belbruno, Edward and Frauenfelder, Urs and {van }Koert, Otto},
     TITLE = {A family of periodic orbits in the three-dimensional lunar
              problem},
   JOURNAL = {Celestial Mech. Dynam. Astronom.},
  FJOURNAL = {Celestial Mechanics \& Dynamical Astronomy. An International
              Journal of Space Dynamics},
    VOLUME = {131},
      YEAR = {2019},
    NUMBER = {2},
     PAGES = {Paper No. 7, 22},
       DOI = {10.1007/s10569-019-9882-8},
}

@article {joung_2026_bifurcations,
    AUTHOR = {Joung, Chankyu and Koh, Dayung and van Koert, Otto},
     TITLE = {Bifurcations of highly inclined near halo orbits using {M}oser
              regularization},
   JOURNAL = {Celestial Mech. Dynam. Astronom.},
  FJOURNAL = {Celestial Mechanics \& Dynamical Astronomy. An International
              Journal of Space Dynamics},
    VOLUME = {138},
      YEAR = {2026},
    NUMBER = {4},
     PAGES = {Paper No. 31, 42},
      ISSN = {0923-2958,1572-9478},
   MRCLASS = {70H12 (70F07 70G45)},
  MRNUMBER = {5090630},
       DOI = {10.1007/s10569-026-10306-1},
}

@book{szebehely,
	address = {New~York and London},
	author = {Szebehely, Victor},
	hyphenation = {english},
	language = {english},
	publisher = {Academic Press},
	title = {Theory of Orbits. The Restricted Problem of Three Bodies},
	year = {1967}}

@article{lamb_roberts,
	author = {Lamb, J.S.W. and Roberts, J.A.G.},
	journal = {Physica D},
	title = {Time-reversal symmetry in dynamical systems: A survey},
	year = {1998},
    volume = {112},
    pages = {1--39},
    doi = {10.1016/S0167-2789(97)00199-1}}

@article{bray_goudas,
	author = {Bray, T.A. and Goudas, C.L.},
	journal = {The Astronomical Journal},
	title = {Doubly Symmetric Orbits about the Collinear Lagrangian Points},
	year = {1967},
    volume = {72},
    number = {2},
    pages = {202--213}}

@article{robin_markellos_1980,
	author = {Robin, I.A. and Markellos, V.V.},
	journal = {Celest. Mech.},
	title = {Numerical determination of three-dimensional periodic orbits generated from vertical self-resonant satellite orbits},
	year = {1980},
    volume = {21},
    pages = {395--434},
    DOI = {10.1007/BF01231276}}

@article{robin_markellos_1983,
	author = {Robin, I.A. and Markellos, V.V.},
	journal = {International Astronomical Union Colloquium, Volume 74: Dynamical Trapping and Evolution in the Solar System},
	title = {Part IV - Periodic Orbits},
	year = {1983},
    pages = {213--224},
    DOI = {10.1017/S0252921100097104}}

@article{aydin_dro,
	author = {Aydin, C.},
	title = {Exploration of vertical self-resonant bifurcations from Distant Retrograde Orbits (DROs) in the Earth-Moon Circular Restricted Three-Body Problem (CR3BP)},
	year = {2026},
    url = {https://arxiv.org/abs/2508.17286}}

@article{aydin_comet,
	author = {Aydin, C.},
	title = {Comet-type periodic motions and their out-of-plane bifurcations in the Earth-Moon CR3BP: a computational symplectic analysis},
	year = {2026},
    url = {https://arxiv.org/abs/2602.16354}}

@article{aydin_cz,
	author = {Aydin, C.},
    journal = {Celest Mech Dyn Astron},
	title = {The Conley--Zehnder indices of the spatial Hill three-body problem},
	year = {2023},
    volume = {135},
    number = {32},
    DOI = {10.1007/s10569-023-10134-7}}

@article{aydin_symmetries,
	author = {Aydin, C.},
    journal = {Archiv der Mathematik},
	title = {The linear symmetries of Hill's lunar problem},
	year = {2023},
    volume = {120},
    DOI = {10.1007/s00013-022-01822-1}}

@article{aydin_contact,
	author = {Aydin, C.},
    journal = {Arnold Mathematical Journal},
	title = {Contact geometry of Hill's approximation in a spatial restricted four-body problem},
	year = {2025},
    volume = {11},
    DOI = {10.56994/ARMJ.011.001.005}}

@article{moreno_aydin,
	author = {Moreno, A. and Aydin, C. and van Koert, O. and Frauenfelder, U. and Koh, D.},
    journal = {J Astronaut Sci},
	title = {Bifurcation Graphs for the CR3BP via Symplectic Methods},
	year = {2024},
    volume = {71},
    number = {51},
    DOI = {10.1007/s40295-024-00462-7}}

@article{broucke,
	author = {Broucke, R.A.},
    journal = {JPL Technical Report 32-1168},
	title = {Periodic Orbits in the Restricted Three-Body Problem with Earth-Moon Masses},
	year = {1968}}

@article{broucke_stab,
	author = {Broucke, R.A.},
    journal = {AIAA Journal},
	title = {Stability of periodic orbits in the elliptic, restricted three-body problem},
	year = {1969},
    volume = {7},
    doi = {10.2514/3.5267}}

@Article{biscani_izzo,
  author  = {Biscani, F. and Izzo, D.},
  title   = {Revisiting high-order Taylor methods for astrodynamics and celestial mechanics},
  journal = {Monthly Notices of the Royal Astronomical Society},
  year    = {2021},
  volume  = {504},
  pages   = {2614--2628},
  doi     = {10.1093/mnras/stab1032},
  issue   = {2}}

@Article{farquhar_kamel,
  author  = {Farquhar, R.W. and Kamel, A.A.},
  title   = {Quasi-periodic orbits about the translunar libration point},
  journal = {Celestial Mechanics},
  year    = {1973},
  volume  = {7},
  pages   = {458--473},
  doi     = {10.1007/BF01227511}}

@Article{breakwell_brown,
  author  = {Breakwell, J.V. and Brown, J.V.},
  title   = {The ‘Halo’ family of 3-dimensional periodic orbits in the Earth-Moon restricted 3-body problem},
  journal = {Celestial Mechanics},
  year    = {1979},
  volume  = {20},
  pages   = {389--404},
  doi     = {10.1007/BF01230405}}

@Article{howell_breakwell,
  author  = {Howell, K.C. and Breakwell, J.V.},
  title   = {Almost rectilinear halo orbits},
  journal = {Celestial Mechanics},
  year    = {1984},
  volume  = {32},
  pages   = {29--52},
  doi     = {10.1007/BF01358402}}

@Article{doedel,
  author  = {Doedel, E.J. and Romanov, V.A. and Paffenroth, R.C. and Keller, H.B. and Dichmann, D.J. and Galán-Vioque, J. and Vanderbauwhede, A.},
  title   = {Elemental Periodic Orbits Associated with the Libration Points in the Circular Restricted 3-Body Problem},
  journal = {International Journal of Bifurcation and Chaos},
  year    = {2007},
  volume  = {17},
  number  = {08},
  pages   = {2625--2677},
  doi     = {10.1142/S0218127407018671}}

@book {arnold_avez,
    AUTHOR = {Arnold, V. I. and Avez, A.},
     TITLE = {Ergodic problems of classical mechanics},
      NOTE = {Translated from the French by A. Avez},
 PUBLISHER = {W. A. Benjamin, Inc., New York-Amsterdam},
      YEAR = {1968},
     PAGES = {ix+286},
   MRCLASS = {28.70 (70.00)},
  MRNUMBER = {232910},
}

@Article{salamon_zehnder,
  author  = {Salamon, D. and Zehnder, E.},
  title   = {Morse theory for periodic solutions of Hamiltonian systems and the Maslov index},
  journal = {Commun Pure Appl Math},
  year    = {1992},
  volume  = {45},
  number = {10},
  pages   = {1303--1360},
  doi     = {10.1002/cpa.3160451004}}

@Article{conley_zehnder,
  author  = {Conley, C. and Zehnder, E.},
  title   = {Morse-type index theory for flows and periodic solutions for Hamiltonian Equations},
  journal = {Commun Pure Appl Math},
  year    = {1984},
  volume  = {37},
  number = {2},
  pages   = {207--253},
  doi     = {10.1002/cpa.3160370204}}

@Article{ginzburg,
  author  = {Ginzburg, V.L.},
  title   = {The Conley conjecture},
  journal = {Annals of mathematics},
  year    = {2010},
  volume  = {172},
  number = {2},
  pages   = {1127--1180},
  doi     = {10.4007/annals.2010.172.1127}}

@book {long,
    AUTHOR = {Long, Yiming},
     TITLE = {Index theory for symplectic paths with applications},
    SERIES = {Progress in Mathematics},
    VOLUME = {207},
 PUBLISHER = {Birkh\"auser Verlag, Basel},
      YEAR = {2002},
     PAGES = {xxiv+380},
      ISBN = {3-7643-6647-8},
   MRCLASS = {37J45 (34C25 37-02 53D40 58-02 58E05 70H05)},
  MRNUMBER = {1898560},
MRREVIEWER = {Thomas\ Bartsch},
       DOI = {10.1007/978-3-0348-8175-3},
}

@incollection{Eliashberg2010,
  title={Introduction to symplectic field theory},
  author={Eliashberg, Yakov and Glvental, A and Hofer, Helmut},
  booktitle={Visions in Mathematics: GAFA 2000 Special Volume, Part II},
  pages={560--673},
  year={2000},
  publisher={Springer},
  doi="10.1007/978-3-0346-0425-3_4"
}

@online{JPL_periodic_orbits,
  author       = {{Jet Propulsion Laboratory}},
  title        = {JPL Three-Body Periodic Orbit Catalog},
  url          = {https://ssd.jpl.nasa.gov/tools/periodic_orbits.html},
  note         = {Accessed: 2026-05-26},
}

@article {bury,
    AUTHOR = {Bury, Luke and McMahon, Jay and Lo, Martin},
     TITLE = {A study of periodic orbits near {E}uropa},
   JOURNAL = {Celestial Mech. Dynam. Astronom.},
  FJOURNAL = {Celestial Mechanics \& Dynamical Astronomy. An International
              Journal of Space Dynamics},
    VOLUME = {134},
      YEAR = {2022},
    NUMBER = {3},
     PAGES = {Paper No. 27, 22},
      ISSN = {0923-2958,1572-9478},
   MRCLASS = {70K42 (70F15)},
  MRNUMBER = {4431077},
       DOI = {10.1007/s10569-022-10076-6},
}

@article{henon_1969,
  title={Numerical exploration of the restricted problem, V},
  author={H{\'e}non, Michel},
  journal={Astronomy and Astrophysics},
  volume={1},
  pages={223--238},
  year={1969}
}

@book{celletti,
  title={Stability and chaos in celestial mechanics},
  author={Celletti, Alessandra},
  year={2010},
  publisher={Springer Science \& Business Media}
}

@article{levi_civita,
  title={Traiettorie singolari ed urti nel problema ristretto dei tre corpi},
  author={Levi-Civita, T.},
  journal={Annali di Matematica, Serie III},
  volume={9},
  pages={1--32},
  year={1904},
  DOI={10.1007/BF02419867}
}

@article{kustaanheimo_stiefel,
  title={Perturbation theory of Kepler motion based on spinor regularization},
  author={Kustaanheimo, P. and Stiefel, E.},
  journal={J. Reine Angew. Math},
  volume={218},
  pages={204--219},
  year={1965}
}

@article{ligon_schaaf,
  title={On the Global Symmetry of the Classical Kepler Problem},
  author={Ligon, T. and Schaaf, M.},
  journal={Reports on Mathematical Physics},
  volume={9},
  pages={281--300},
  year={1976}
}

@article{birkhoff,
  title={The restricted problem of three bodies},
  author={Birkhoff, G.D.},
  journal={Rend. Circ. Mat. Palermo},
  volume={39},
  pages={265--334},
  year={1915},
  DOI={10.1007/BF03015982}
}

@article{waldvogel,
  title={Die Verallgemeinerung der Birkhoff-Regularisierung für das räumliche Dreikörperproblem},
  author={Waldvogel, J.},
  journal={Bulletin astronomique},
  pages={295--341},
  year={1967}
}

@article{kalantonis,
  title={Numerical Investigation for Periodic Orbits in the Hill Three-Body Problem},
  author={Kalantonis, V.S.},
  journal={Universe},
  volume={6},
  year={2020},
  DOI={10.3390/universe6060072}
}

@article{ragos_kalantonis,
  title={Numerical Investigation for 3D Branches of the Lyapunov Families in the Hill’s Problem with Radiation Pressure},
  author={Ragos, O. and Perdiou, A.E. and Perdios, E.A. and Kalantonis, V.S.},
  journal={Mathematics},
  volume={14},
  year={2026},
  DOI={10.3390/math14101625}
}

@article{schmidt,
  title={Families of Periodic Orbits in the Restricted Problem of Three Bodies Connecting Families of Direct and Retrograde Orbits},
  author={Schmidt, D.S.},
  journal={SIAM J. Appl. Math.},
  volume={22},
  year={1972},
  DOI={10.1137/0122004}
}

@article{palacian,
  title={Invariant tori of rectilinear type in the spatial three-body problem},
  author={Palacián, F. and Sayas, F. and Yanguas, P.},
  journal={Journal of Differential Equations},
  volume={399},
  year={2024},
  pages={82--180},
  DOI={10.1016/j.jde.2024.03.008}
}

@article{korneev,
  title={Low-Thrust Trajectory Optimization in Kustaanheimo–Stiefel Variables},
  author={Korneev, K. and Trofimov, S.P.},
  journal={Cosmic Research},
  volume={62},
  year={2024},
  pages={256--265},
  DOI={10.1134/S0010952524600288}
}

@article{oshima,
  title={Analysis of medium-energy transfers to the Moon},
  author={Oshima, K. and Topputo, F. and Campagnola, S. and Yanao, T.},
  journal={Celest Mech Dyn Astron},
  volume={127},
  year={2017},
  pages={285--300},
  DOI={10.1007/s10569-016-9727-7}
}

@article{lega,
  title={Detection of close encounters and resonances in three-body problems through Levi-Civita regularization},
  author={Lega, E. and Guzzo, M. and Froeschlé, C.},
  journal={Monthly Notices of the Royal Astronomical Society},
  volume={418},
  year={2011},
  pages={107--113},
  DOI={10.1111/j.1365-2966.2011.19467.x}
}

@article{cardin,
  title={Integrability of close encounters in the spatial restricted three-body problem},
  author={Cardin, Franco and Guzzo, Massimiliano},
  journal={Communications in Contemporary Mathematics},
  volume={24},
  number={06},
  pages={2150040},
  year={2022},
  publisher={World Scientific},
  DOI={10.1142/S0219199721500401}
}

@article {albers,
    AUTHOR = {Albers, Peter and Frauenfelder, Urs and van Koert, Otto and
              Paternain, Gabriel P.},
     TITLE = {Contact geometry of the restricted three-body problem},
   JOURNAL = {Comm. Pure Appl. Math.},
  FJOURNAL = {Communications on Pure and Applied Mathematics},
    VOLUME = {65},
      YEAR = {2012},
    NUMBER = {2},
     PAGES = {229--263},
      ISSN = {0010-3640,1097-0312},
   MRCLASS = {53D35 (70F07 70G45)},
  MRNUMBER = {2855545},
       DOI = {10.1002/cpa.21380},
}

@article{joung_koert,
  title={Computational symplectic topology and symmetric orbits in the restricted three-body problem},
  author={Joung, C. and van Koert, O.},
  journal={Nonlinearity},
  volume={38},
  year={2025},
  DOI={10.1088/1361-6544/ada7ba}
}

@article{sundman,
  title={Recherches sur le problème des trois corps},
  author={Sundman, K.F.},
  journal={Acta Societatis Scientificae Fennicae},
  volume={34},
  year={1907}
}

@article{kummer,
  title={On the regularization of the Kepler problem},
  author={Kummer, M.},
  journal={Commun. Math. Phys.},
  volume={84},
  pages={133--152},
  year={1982}
}

@book{frauenfelder_koert,
	author = {Frauenfelder, U. and van Koert, O.},
	language = {english},
	publisher = {Birkhäuser Cham},
	title = {The Restricted Three-Body Problem and Holomorphic Curves},
	year = {2018},
    DOI={10.1007/978-3-319-72278-8}}

@book{koon,
  author    = {Koon, W. S. and Lo, M. W. and Marsden, J. E. and Ross, S. D.},
  title     = {Dynamical Systems, the Three-Body Problem and Space Mission Design},
  year      = {2022},
  publisher = {Marsden Books},
  isbn      = {978-0-615-24095-4},
}

@article{howell_barden,
  title={Application of Dynamical Systems Theory to Trajectory Design for a Libration Point Mission},
  author={Howell, K.C. and Barden, B.T. and Lo, M.W.},
  journal={J of Astronaut Sci},
  volume={45},
  pages={161--178},
  year={1997},
  DOI={10.1007/BF03546374}
}

@article {capinski_kepley_mireles_james,
    AUTHOR = {Capi\'nski, Maciej J. and Kepley, Shane and Mireles James, J.
              D.},
     TITLE = {Computer assisted proofs for transverse collision and near
              collision orbits in the restricted three body problem},
   JOURNAL = {J. Differential Equations},
  FJOURNAL = {Journal of Differential Equations},
    VOLUME = {366},
      YEAR = {2023},
     PAGES = {132--191},
      ISSN = {0022-0396,1090-2732},
   MRCLASS = {70F07 (34C25 34C37 37C29 37J46)},
  MRNUMBER = {4576879},
       DOI = {10.1016/j.jde.2023.03.053},
}

@article {capinski_pasiut,
    AUTHOR = {Capi\'nski, Maciej J. and Pasiut, Aleksander},
     TITLE = {Oscillatory collision approach in the {E}arth-{M}oon
              restricted three body problem},
   JOURNAL = {Commun. Nonlinear Sci. Numer. Simul.},
  FJOURNAL = {Communications in Nonlinear Science and Numerical Simulation},
    VOLUME = {152},
      YEAR = {2026},
     PAGES = {Paper No. 109173, 24},
      ISSN = {1007-5704,1878-7274},
   MRCLASS = {37N05 (37C29 37J46 70F07)},
  MRNUMBER = {4941140},
MRREVIEWER = {Martha\ Alvarez-Ram\'irez},
       DOI = {10.1016/j.cnsns.2025.109173},
}

@book {roa,
    AUTHOR = {Roa, Javier},
     TITLE = {Regularization in orbital mechanics},
    SERIES = {De Gruyter Studies in Mathematical Physics},
    VOLUME = {42},
 PUBLISHER = {De Gruyter, Berlin},
      YEAR = {2017},
     PAGES = {xv+403},
   MRCLASS = {70M20 (35Q70 85A05)},
  MRNUMBER = {3793192},
       DOI = {10.1515/9783110559125},
}

@book {gomez_et_al,
    AUTHOR = {G\'omez, G. and Llibre, J. and Mart\'inez, R. and Sim\'o, C.},
     TITLE = {Dynamics and mission design near libration points. {V}ol.
              {I}.},
    SERIES = {World Scientific Monograph Series in Mathematics},
    VOLUME = {2},
 PUBLISHER = {World Scientific Publishing Co., Inc., River Edge, NJ},
      YEAR = {2001},
     PAGES = {xviii+443},
   MRCLASS = {70M20 (37N05 70-08 70F07 70F15)},
  MRNUMBER = {1867240},
MRREVIEWER = {Wang\ Sang\ Koon},
       DOI = {10.1142/9789812810632\_bmatter},
}

@article {howell,
    AUTHOR = {Howell, K. C.},
     TITLE = {Three-dimensional, periodic, ``halo''\ orbits},
   JOURNAL = {Celestial Mech.},
  FJOURNAL = {Celestial Mechanics},
    VOLUME = {32},
      YEAR = {1984},
    NUMBER = {1},
     PAGES = {53--71},
      ISSN = {0008-8714},
   MRCLASS = {70F07 (70-08)},
  MRNUMBER = {738148},
       DOI = {10.1007/BF01358403},
}

@article {zimovan_spreen,
    AUTHOR = {Zimovan-Spreen, Emily M. and Howell, Kathleen C. and Davis,
              Diane C.},
     TITLE = {Near rectilinear halo orbits and nearby higher-period
              dynamical structures: orbital stability and resonance
              properties},
   JOURNAL = {Celestial Mech. Dynam. Astronom.},
  FJOURNAL = {Celestial Mechanics \& Dynamical Astronomy. An International
              Journal of Space Dynamics},
    VOLUME = {132},
      YEAR = {2020},
    NUMBER = {5},
     PAGES = {Paper No. 28, 25},
      ISSN = {0923-2958,1572-9478},
   MRCLASS = {70M20},
  MRNUMBER = {4114268},
       DOI = {10.1007/s10569-020-09968-2},
}

@phdthesis{campbell,
  title={Bifurcations from families of periodic solutions in the circular restricted problem with application to trajectory design},
  author={Campbell, E. T.},
  year={1999},
  school={Purdue University}
}

@article{koblick_kelly,
  title={Novel three-body tulip-shaped orbit families for lunar missions},
  author={Koblick, Darin C. and Kelly, Patrick},
  journal={The Journal of the Astronautical Sciences},
  volume={72},
  number={4},
  pages={32},
  year={2025},
  publisher={Springer}
}

@article{zhang,
  title={Time-regularized bifurcation framework for constructing tulip orbits in the CRTBP},
  author={Zhang, Jiye and Jiang, Xinyu and Yuan, Yuan and Dai, Honghua},
  journal={Nonlinear Dynamics},
  volume={114},
  number={9},
  pages={620},
  year={2026},
  publisher={Springer}
}

@article{restrepo_russell,
  title={A database of planar axisymmetric periodic orbits for the Solar system},
  author={Restrepo, R.L. and Russell, R.P.},
  journal={Celest. Mech. Dyn. Astron.},
  volume={130},
  number={49},
  year={2018},
  DOI = {10.1007/s10569-018-9844-6},
}

@book {kuznetsov,
    AUTHOR = {Kuznetsov, Yuri A.},
     TITLE = {Elements of applied bifurcation theory},
    SERIES = {Applied Mathematical Sciences},
    VOLUME = {112},
   EDITION = {Fourth},
 PUBLISHER = {Springer, Cham},
      YEAR = {2023},
     PAGES = {xxvi+703},
   MRCLASS = {37Gxx (34C23)},
  MRNUMBER = {4592585},
       DOI = {10.1007/978-3-031-22007-4},
}

@article{henon_1973,
  title={Vertical stability of periodic orbits in the restricted problem. I. Equal masses},
  author={H{\'e}non, Michel},
  journal={Astronomy and Astrophysics, Vol. 28, p. 415 (1973)},
  volume={28},
  pages={415},
  year={1973}
}

\end{document}